\documentclass[11pt,a4paper]{article}
\usepackage[T1]{fontenc}
\usepackage{lmodern}
\usepackage{amsmath,amsthm,amssymb,mathtools,mathrsfs}
\usepackage{geometry,microtype,booktabs,array,enumitem,needspace}
\usepackage[dvipsnames]{xcolor}
\usepackage{tikz}
\usetikzlibrary{decorations.pathreplacing}
\usepackage{hyperref}
\usepackage[nameinlink,noabbrev]{cleveref}
\hypersetup{colorlinks=true,linkcolor=MidnightBlue,citecolor=OliveGreen,urlcolor=MidnightBlue,
pdftitle={Critical Ising and Temperley-Lieb-Jones States},pdfauthor={David E. Evans}}
\theoremstyle{plain}
\newtheorem{theorem}{Theorem}[section]
\newtheorem{proposition}[theorem]{Proposition}
\newtheorem{lemma}[theorem]{Lemma}
\newtheorem{corollary}[theorem]{Corollary}
\theoremstyle{definition}

\newtheorem{remark}[theorem]{Remark}
\newcommand{\C}{\mathcal C}

\newcommand{\DHR}{\mathrm{DHR}}
\newcommand{\Rep}{\mathrm{Rep}}
\newcommand{\Vect}{\mathrm{Vec}}
\newcommand{\Hom}{\mathrm{Hom}}
\newcommand{\End}{\mathrm{End}}
\newcommand{\Ad}{\mathrm{Ad}}
\newcommand{\ran}{\mathrm{ran}}
\newcommand{\sgn}{\operatorname{sgn}}
\newcommand{\id}{\mathrm{id}}
\newcommand{\one}{\mathbf 1}
\newcommand{\Z}{\mathbb Z}
\newcommand{\R}{\mathbb R}
\newcommand{\T}{\mathbb T}
\newcommand{\K}{\mathcal K}
\newcommand{\HH}{\mathcal H}
\newcommand{\EE}{\mathsf E}
\newcommand{\HS}{\mathcal L^2}

\newcommand{\Pf}{\operatorname{Pf}}

\title{\textbf{Critical Ising and Temperley--Lieb--Jones States:\\
Cone Haag Duality and the Trivial Algebra Object}}
\author{David E.~Evans}
\begin{document}
\maketitle

\begin{abstract}
We compute the algebra object of two critical Ising states in the operator algebraic framework of fusion category symmetry. The boundary subalgebra is the even part of the Pauli chain, and its DHR category is the toric code. We treat two models. The first is the classical Ising model in two dimensions of Araki--Evans and Evans--Lewis, through its row transfer matrix. The second is the Ising chain in a transverse field of Araki and Araki--Matsui. For the critical transfer matrix state and for the ground state of the critical chain, the algebra object is exactly the tensor unit. The calculation uses explicit sector representatives, through which the multiplicities become the equivalence criteria of Araki and Evans for even restrictions of Fock states. At criticality the half-line commutator is not Hilbert--Schmidt, as found by Araki--Evans and by Araki--Matsui, and this excludes the two bosonic sectors. Foit's twisted duality gives cone Haag duality for any pure quasifree state. At criticality the half-line algebras are factors and the half-line field subspaces are standard. The two states are related by a graded automorphism with exponentially decaying locality bounds. In a second part we treat the critical Temperley--Lieb--Jones chains at $\delta=2\cos(\pi/(k+2))$. For every $k\ge2$ the ground state selected by the vacuum block is pure and has trivial algebra object in the Drinfeld centre of the even part. In the critical quantum random-cluster model there are separating loops at all scales, and these force every charged translation coefficient to decay. At $k=2$ this recovers the critical Ising chain.
\end{abstract}

\tableofcontents

\section{Introduction}
\label{sec:intro}

A physical boundary inclusion $B\subseteq A$ associates categorical data to symmetric states on an infinite spin chain \cite{EvansJones2025}. A state $\phi$ whose restriction to $B$ is pure determines a connected algebra object $H_\phi\in\DHR(B)$. Topological states are those for which this algebra is Lagrangian. The framework also predicts gapless states, in the categorical sense that their algebra objects are not Morita equivalent to Lagrangian algebras. An explicit computation for the Ising model provides a test of these definitions on a familiar critical system.

We take
\[
 A^P=\bigotimes_{j\in\Z}M_2(\mathbb C),\qquad
 B=(A^P)^{\Lambda},\qquad
 \Lambda=\bigotimes_{j\in\Z}\Ad(\sigma_z^{(j)}).
\]
The inclusion is a physical boundary inclusion with symmetry category $\Vect(\Z_2)$ and
\[
 \DHR(B)\simeq\mathcal Z(\Rep(\Z_2)),
\]
the toric code category \cite{Jon24,EvansJones2025}. Its simple objects are $\one,e,m,f=em$, where $e,m$ are bosons and $f$ is a fermion. We choose $e$ to label the disorder sector and $m$ the spin-charge sector. The precise representatives are fixed in Section~\ref{sec:sectors}.

\subsection{The two Ising models}
\label{sec:models}

Two different Ising models give states on $A^P$. Both are written in terms of the Pauli matrices $\sigma_x^{(j)},\sigma_y^{(j)},\sigma_z^{(j)}$ at the sites $j\in\Z$, and $\Lambda$ is the global spin flip $\sigma_x^{(j)}\mapsto-\sigma_x^{(j)}$, $\sigma_z^{(j)}\mapsto\sigma_z^{(j)}$.

\paragraph{The two-dimensional classical model.}
Araki--Evans and Evans--Lewis \cite{ArakiEvans83,EvansLewis86} consider the classical Ising model on $\Z^2$ with spins $\xi_{j,k}=\pm1$ and energy
\[
 H(\xi)=-\sum_{(j,k)\in\Z^2}\bigl(J_1\,\xi_{j,k}\xi_{j,k+1}+J_2\,\xi_{j,k}\xi_{j+1,k}\bigr),
 \qquad J_1,J_2>0,
\]
where $k$ labels the rows. Put $K_i=\beta J_i$ and define $K_1^*$ by $\tanh K_1^*=e^{-2K_1}$. Identifying the spins of a row with the eigenvalues of the $\sigma_x^{(j)}$, the row-to-row transfer matrix is, up to a positive scalar,
\begin{equation}
 T=Y^{1/2}XY^{1/2},\qquad
 Y=\exp\Bigl(K_2\sum_j\sigma_x^{(j)}\sigma_x^{(j+1)}\Bigr),\qquad
 X=\exp\Bigl(K_1^*\sum_j\sigma_z^{(j)}\Bigr),
 \label{eq:transfer}
\end{equation}
on finite rows \cite[(2.4)]{EvansLewis86}. The normalized Perron--Frobenius vector states of the finite volume transfer matrices converge to a $\Lambda$-invariant state $\phi_\beta$ of $A^P$, given explicitly by Proposition~\ref{prop:transferconvention} and \eqref{eq:offcritical-symbol}. The critical inverse temperature $\beta_c$ is given by $K_1^*=K_2$, that is, $\sinh(2K_1)\sinh(2K_2)=1$. Araki and Evans proved that the GNS representation of $\phi_\beta$ is factorial for $\beta\le\beta_c$ and has two-dimensional centre for $\beta>\beta_c$ \cite[Theorem~1]{ArakiEvans83}. Their proof, through \cite[Theorem~2 and Lemmas~5.1, 5.2]{ArakiEvans83}, shows that $\phi_\beta$ is pure for $\beta\le\beta_c$ and is the average of two inequivalent pure states for $\beta>\beta_c$. Off criticality these pure states are ground states of the one-dimensional dynamics determined by the transfer matrix, with spectral gap $4(K_2-K_1^*)$ for $\beta>\beta_c$ and $2(K_1^*-K_2)$ for $\beta<\beta_c$ \cite[Proposition~6.2]{Araki86}. We write $\phi_c=\phi_{\beta_c}$.

\paragraph{The quantum chain.}
Araki and Araki--Matsui \cite{Araki84,ArakiMatsui85} determine the ground states of the $XY$ chain. Its Ising case is the transverse field Hamiltonian
\begin{equation}
 H_\lambda=-\frac12\sum_{j\in\Z}
 \bigl(\sigma_x^{(j)}\sigma_x^{(j+1)}+\lambda\,\sigma_z^{(j)}\bigr),
 \qquad\lambda\in\R,
 \label{eq:Hlambda}
\end{equation}
which is the case $J=1/4$, $\gamma=1$ of \cite[(1.1)]{ArakiMatsui85}. The infinite sum specifies an interaction and its dynamics, rather than an element of $A^P$. By \cite[Theorem~1 and Proposition~3.3]{ArakiMatsui85}, $H_\lambda$ has a unique ground state for $|\lambda|\ge1$, which is pure, while for $0<|\lambda|<1$ it has exactly two extremal ground states, which are pure and are interchanged by $\Lambda$. The one-particle energies are $\varepsilon_\lambda(\theta)=(1-2\lambda\cos\theta+\lambda^2)^{1/2}\ge\bigl|1-|\lambda|\bigr|$ (see the proof of Proposition~\ref{prop:groundstate}). The infimum $\bigl|1-|\lambda|\bigr|$ is attained at $\theta=0$ for $\lambda>0$ and at $\theta=\pi$ for $\lambda<0$, so the one-particle spectrum has a positive gap precisely when $|\lambda|\ne1$. The critical Hamiltonian is
\begin{equation}
 H_1=-\frac12\sum_{j\in\Z}
 \bigl(\sigma_x^{(j)}\sigma_x^{(j+1)}+\sigma_z^{(j)}\bigr),
 \label{eq:H1}
\end{equation}
and $\phi_\sharp$ denotes its unique ground state.

\paragraph{The two critical states.}
Thus $\phi_c$ belongs to the classical model at $\beta_c$ and $\phi_\sharp$ to the quantum chain at $\lambda=1$. Both are pure and $\Lambda$-invariant on $A^P$, and pure on $B$. They are distinct for every finite critical coupling. The quantum chain state is invariant under Kramers--Wannier duality, while the transfer matrix state is not (Section~\ref{sec:states}). The state $\phi_\sharp$ is the anisotropic limit $K_2=K_1^*\to0$ of the critical transfer matrix states (Remark~\ref{rem:anisotropic}).

\subsection{Results}

A state $\phi$ on $A^P$ whose restriction to $B$ is pure induces a DHR($B$) bimodule $H_\phi$, the internal end of $L^2(B,\phi|_B)$, which is a connected algebra object (Section~\ref{sec:invariant}). We analyse these algebra objects for the two critical states.

\begin{theorem}[Main result]
\label{thm:main}
For the two critical states,
\[
 H_{\phi_c}=H_{\phi_\sharp}=\one.
\]
Both restrictions to $B$ satisfy cone Haag duality with zero spread, and their half-line observable algebras are factors. There is a graded, translation covariant automorphism $\alpha$ of $A^P$, with exponentially decaying locality bounds for $\alpha$ and $\alpha^{-1}$, such that
\[
 \phi_\sharp=\phi_c\circ\alpha.
\]
The calculation of the algebra object and the duality theorem are independent of this conjugacy.
\end{theorem}

The sector calculation is the main point. Each multiplicity in $H_\phi$ records whether a DHR automorphism is implementable in the irreducible representation of $B$ defined by the state. With explicit sector representatives these multiplicities are exactly the equivalence criteria of Araki and Evans for restrictions of Fock states to the even algebra \cite[Lemmas~4.2, 4.3 and Theorem~4]{ArakiEvans83}. At criticality the half-line commutator is not Hilbert--Schmidt, by Araki--Evans for the transfer matrix state and by Araki--Matsui for the quantum chain, and both bosonic sectors are excluded. This determines the algebra object directly, without a classification of all Morita classes or an inference from failure of a particular conjugacy.

Cone Haag duality has a separate and more general explanation. Foit's theorem applies to every closed real one-particle subspace. Passing explicitly to the even vacuum sector yields duality for the half-line observable algebras of any pure quasifree state. Generic position is not needed for that argument. We nevertheless establish generic position for the two critical projections. It says that the half-line real subspaces are standard, so that the Fock vacuum is cyclic and separating for the corresponding field algebras. For $\phi_c$, the proof reduces to the impossibility of a meromorphic square root of a rational function with a pole of odd order. Standardness does not determine the type of the half-line algebras. At criticality we show separately, from the same half-line jump that excludes the bosonic sectors, that the half-line field and observable algebras are factors (Proposition~\ref{prop:factor}).

The critical quantum chain state is also explicit in the Temperley--Lieb presentation. Let $m_{2j-1},m_{2j}$ be the two Majorana operators at site $j$ (Section~\ref{sec:states}). Then
\[
 e_k=\frac12(1+i m_km_{k+1}),\qquad
 \phi_\sharp(e_k)=\frac12+\frac1\pi.
\]
Its Majorana covariance is $2/(\pi n)$ on odd nonzero lags and zero on even nonzero lags. These formulae identify the known critical Ising ground state within the categorical framework. The continuum interpretation as the free Majorana/Ising theory is established background, not a consequence of the equality $H_\phi=\one$ alone.

The Hamiltonian result extends beyond free fermions. For $k\ge2$ and $\delta=2\cos(\pi/(k+2))$, let $\C_k$ be the even part of the unitary Temperley--Lieb--Jones category and $B_k$ the fusion chain $\overline{\bigcup_n\End_{\C_k}(Y^{\otimes n})}$, $Y=X\otimes X$, so that $\DHR(B_k)\simeq\mathcal Z(\C_k)$. At $k=2$, $B_2=B$ and $-\sum_ie_i=H_1+\text{const}$ (Remark~\ref{rem:k2chain}).

Let $\omega_k$ be the limit of the Perron--Frobenius ground vectors of the finite volume Hamiltonians $-\sum_ie_i$ restricted to the vacuum blocks $\Hom(\one,X^{\otimes2N})$, the spaces of morphisms from the tensor unit (Section~\ref{sec:TLJchain}).

\begin{theorem}
\label{thm:main2}
For every $k\ge2$, the state $\omega_k$ is pure and translation invariant, and $H_{\omega_k}=\one$ in $\mathcal Z(\C_k)$. At $k=2$, $\omega_2\circ\EE=\phi_\sharp$.
\end{theorem}

The mechanism differs from that of Theorem~\ref{thm:main}. The finite volume states are expectations for the critical quantum random-cluster model at $q=\delta^2$. The charged translation coefficients become matrix coefficients of affine Temperley--Lieb modules, each loop separating the two insertions costs a factor less than one, and crossing estimates give logarithmically many such loops. Decay then excludes every nontrivial sector. At $k=2$ the coefficients are the order, disorder and fermion correlators.

Part~\ref{part:Ising} treats the Ising states. Section~\ref{sec:states} fixes the state and sign conventions. Section~\ref{sec:duality} proves general half-line duality. Section~\ref{sec:sectors} computes the algebra objects. Section~\ref{sec:conjugacy} proves conjugacy and compares the dynamics, and Section~\ref{sec:TL} gives the Temperley--Lieb formulae. Part~\ref{part:TLJ} treats the Temperley--Lieb--Jones chains. Section~\ref{sec:TLJchain} sets up the chain and the case $k=2$, Sections~\ref{sec:loops}--\ref{sec:purity} construct the state, Sections~\ref{sec:affine}--\ref{sec:decay} prove Theorem~\ref{thm:main2}, and Section~\ref{sec:scope} gives an overview of what extends and what remains open. Appendix~\ref{sec:standardness} establishes half-line standardness.

\part{Ising states}
\label{part:Ising}

\section{The even algebra and the two critical states}
\label{sec:states}

\subsection{The categorical invariant}
\label{sec:invariant}

Write $\EE:A^P\to B$ for the averaging expectation,
\[
 \EE(a)=\frac12(a+\Lambda(a)).
\]
A symmetric state on $A^P$ is precisely $\psi\circ\EE$ for a state $\psi$ on $B$. If $\psi$ is pure, the internal end of its GNS representation is a connected Q-system \cite{EvansJones2025}
\begin{equation}
 H_\psi=\underline{\End}\bigl(L^2(B,\psi)\bigr)\in\DHR(B).
 \label{eq:internalend}
\end{equation}
We use $H_\phi$ for this same object when $\phi=\psi\circ\EE$. A state is \emph{local} if $H_\phi$ is braided commutative, \emph{topological} if it is Lagrangian, and \emph{gapless} if it is not Morita equivalent to a Lagrangian algebra \cite[Section~5]{EvansJones2025}. The last two terms are categorical definitions. Assertions about a spectral gap require a specified dynamics.

In the toric code the connected commutative Q-systems are $\one$, $\one\oplus e$, and $\one\oplus m$, and the last two are Lagrangian. A connected algebra Morita equivalent to $\one$ is the internal end $X\otimes\overline X$ of a simple object $X$ \cite{Ost03}, which is $\one$ here since all simple objects are invertible. Thus the tensor unit is local and categorically gapless. Commutativity without the Lagrangian property does not make a state gapped.

\subsection{CAR and Pauli conventions}
\label{sec:CAR}

Let $A^F$ be the CAR algebra with generators $c_j,c_j^*$, and put
\[
 a_j=c_j+c_j^*,\qquad b_j=-i(c_j-c_j^*),\qquad
 m_{2j-1}=a_j,\quad m_{2j}=b_j.
\]
The $m_k$ are self-adjoint unitaries satisfying $\{m_k,m_l\}=2\delta_{kl}$. The Jordan--Wigner identification of even algebras gives
\begin{equation}
 B=(A^P)_+=(A^F)_+,\qquad
 \sigma_z^{(j)}=i a_jb_j,\quad
 \sigma_x^{(j)}\sigma_x^{(j+1)}=i b_ja_{j+1}.
 \label{eq:JW}
\end{equation}
For an interval or a ray $I$ in the site lattice, $B_I$ is the even CAR algebra generated by the Majoranas of the sites in $I$. This agrees with the fixed-point algebra of the Pauli observables supported in $I$.

For self-dual CAR notation take $\K=\ell^2(\Z)\oplus\ell^2(\Z)$,
\[
 C(f,g)=\sum_j(c_j^*f_j+c_jg_j),\qquad
 \mathfrak c(f,g)=(\bar g,\bar f).
\]
Here $\mathfrak c$ is the antiunitary involution implementing the adjoint, $C(\mathfrak ch)=C(h)^*$. We use an inner product linear in its second variable. A basis projection is an orthogonal projection $E$ with $\mathfrak c E\mathfrak c=1-E$. Its pure quasifree state is specified by
\begin{equation}
 \omega_E(C(h)^*C(k))=\langle h,(1-E)k\rangle.
 \label{eq:covariance}
\end{equation}
Thus $E$ is the occupied, negative-energy projection in the convention $C(f,g)=c^*f+cg$. Equivalently, $\omega_E(C(h)C(k)^*)=\langle k,Eh\rangle$. In particular $\omega_E(c_j^*c_j)=(E_{11})_{jj}$. A unitary $v$ commuting with $\mathfrak c$ induces $\tau(v)C(h)=C(vh)$, and
\begin{equation}
 \omega_E\circ\tau(v)=\omega_{v^*Ev}.
 \label{eq:pullback}
\end{equation}
Through the Jordan--Wigner identification \eqref{eq:JW} of the even CAR algebra with $B=(A^P)^\Lambda$, the restriction $\psi_E=\omega_E|_B$ is a state on $B$. Its symmetric Pauli extension is the state $\phi_E=\psi_E\circ\EE$ on $A^P$, with $\EE$ the averaging expectation of Section~\ref{sec:invariant}. It is the unique $\Lambda$-invariant state on $A^P$ restricting to $\psi_E$.

The restriction $\psi_E$ is pure, since a pure even state restricts to a pure state on the fixed-point algebra \cite[Lemma~4.1]{ArakiEvans83}. Purity passes from a pure even state to its restriction to the fixed-point algebra, but not from $\psi_E$ up to the extension $\phi_E$ on $A^P$. For $\beta>\beta_c$ the transfer matrix state $\phi_\beta$ has pure even restriction and is the average of two inequivalent pure states, as recalled in Section~\ref{sec:intro}.

\subsection{The critical symbols}
\label{sec:symbols}

In Fourier representation define
\begin{equation}
 V(\Theta)=\begin{pmatrix}\cos\Theta&-i\sin\Theta\\ i\sin\Theta&-\cos\Theta\end{pmatrix},
 \qquad E(\Theta)=\frac12(1-V(\Theta)).
 \label{eq:V}
\end{equation}
The basis-projection condition holds exactly when $\cos\Theta$ is even and $\sin\Theta$ is odd, or equivalently $\Theta(-\theta)\equiv-\Theta(\theta)\pmod{2\pi}$.

For the transfer matrix state, recall from Section~\ref{sec:models} that criticality is $K_1^*=K_2=:K$, and put $z=\tanh K\in(0,1)$. In the bond-symmetric convention, the critical symbol is
\begin{align}
 \Theta_\sharp(\theta)&=\frac\theta2+\frac\pi2\sgn\theta,
 &\chi_z(\theta)&=\arg(1-z^2e^{i\theta}),\label{eq:angles}\\
 \Theta_c(\theta)&=\Theta_\sharp(\theta)+\chi_z(\theta),
 &E_\sharp&=E(\Theta_\sharp),\quad E_c=E(\Theta_c),\label{eq:criticalangles}
\end{align}
for $-\pi<\theta<\pi$, with values at individual endpoints immaterial. Here $\chi_z$ is the continuous odd branch with $\chi_z(0)=0$. Proposition~\ref{prop:transferconvention} derives this symbol and selects its branch directly, since oddness alone does not choose between a projection and its complement. We also allow $z=0$ in the formulae, where $E_c=E_\sharp$, to denote the limiting state as $K\downarrow0$.

\begin{proposition}[Transfer convention and the occupied projection]
\label{prop:transferconvention}
For $K>0$, take the bond-symmetric transfer matrix \eqref{eq:transfer} at criticality, $K_1^*=K_2=K$,
with sums over finite intervals and open boundary conditions.
Their normalized Perron--Frobenius vector states converge on every local
Pauli observable to the symmetric extension of the quasifree state with
occupied projection $E_c=(1-V(\Theta_\sharp+\chi_z))/2$.
\end{proposition}
This is the critical case of the construction of Araki and Evans \cite[Section~3]{ArakiEvans83}. The limit is the even restriction of the Fock state of the negative spectral projection of the one-particle generator $H$ of \cite[(3.5)]{ArakiEvans83}, whose symbol is \cite[(3.7)--(3.9)]{ArakiEvans83}, and the passage to the limit uses \cite[Lemmas~3.1--3.3]{ArakiEvans83}. We include the derivation in our conventions because it fixes the branch $\Theta_\sharp+\chi_z$ used below.
\begin{proof}
Put $n(\theta)=\cos\theta\,\sigma_3+\sin\theta\,\sigma_2$,
$c_K=\cosh(2K)$, $s_K=\sinh(2K)$, and
$A_\theta=c_K-(c_K-1)\cos\theta$. The CAR commutators of the field
and bond terms give the conjugation matrix
\[
 U_T=e^{-Kn(\theta)}e^{2K\sigma_3}e^{-Kn(\theta)}
     =\cosh\gamma\,1+s_K(\sigma_3-A_\theta n(\theta)),
 \qquad
 \cosh\gamma=c_K^2-s_K^2\cos\theta.
\]
The effective one-particle Hamiltonian is $h_{\mathrm{eff}}=-\log U_T$.
Writing $h_{\mathrm{eff}}/\gamma=v_3\sigma_3+v_2\sigma_2$, and setting
$\zeta=e^{i\theta}$ and $a=z^2$, one finds
\[
 v_3+iv_2=\frac{s_K(A_\theta\zeta-1)}{\sinh\gamma},\qquad
 A_\theta\zeta-1=\frac{(\zeta-1)(1-a\zeta)}{1-a}.
\]
Moreover
\[
 \sinh\gamma=
 \frac{2s_K|\sin(\theta/2)|\,|1-a\zeta|}{1-a}.
\]
Therefore, almost everywhere on $(-\pi,\pi)$,
\[
 v_3+iv_2
 =\frac{\zeta-1}{2|\sin(\theta/2)|}
      \frac{1-a\zeta}{|1-a\zeta|}
 =e^{i\Theta_\sharp(\theta)}e^{i\chi_z(\theta)}.
\]
Thus $h_{\mathrm{eff}}=\gamma V(\Theta_\sharp+\chi_z)$, whose negative
spectral projection is $E_c$. We justify its identification with the
thermodynamic transfer matrix state as follows.

Let $\Lambda_N=[-N,N]$, let $P_N$ be the projection onto both self-dual
components at those sites, and let $T_N$ be the displayed spin transfer
matrix with only bonds contained in $\Lambda_N$. Write $n$ also for the
Laurent operator with symbol $n(\theta)$, and put
\[
 n_N=P_NnP_N,\qquad s_N=P_N\sigma_3P_N,\qquad
 U_N=e^{-Kn_N}e^{2Ks_N}e^{-Kn_N}.
\]
Here $U_N$ is extended by the identity on $(1-P_N)\K$. The finite CAR
calculation gives $T_N C(f)T_N^{-1}=C(U_Nf)$ for $f\in P_N\K$.
Taking the logarithm of this conjugation action gives
\[
 [-\log T_N,C(f)]=C(h_Nf),\qquad h_N=-\log U_N.
\]
Indeed, in an eigenbasis of $T_N$ with positive eigenvalues $t_i$, both
$\log(\Ad T_N)$ and the commutator with $\log T_N$ act on the matrix
unit $e_{ij}$ by multiplication by $\log t_i-\log t_j$.
Consequently $-\log T_N$ is the quadratic Hamiltonian with symbol
$h_N$, up to a scalar. Its ground vector is the dominant eigenvector
of $T_N$.

In the product $\sigma_x$ basis, the field factor has strictly positive
matrix entries and the bond factors are positive diagonal matrices.
Perron--Frobenius therefore makes this eigenvector unique. In particular,
$h_N|_{P_N\K}$ has no zero mode, and the finite volume vector state is
quasifree with occupied projection
$E_N=\mathbf1_{(-\infty,0)}(h_N)|_{P_N\K}$.
Moreover, $T_N$ commutes with the spin flip
$\prod_{j\in\Lambda_N}\sigma_z^{(j)}$, so its unique vector state is
spin-flip invariant.

The operators $n_N$ and $s_N$ converge strongly to $n$ and $\sigma_3$
and have norm at most one. Hence $U_N\to U_T$ strongly, with common
spectral bounds
\[
 e^{-4K}1\le U_N\le e^{4K}1.
\]
Continuous functional calculus gives $h_N\to h_{\mathrm{eff}}$ strongly.
The multiplier $h_{\mathrm{eff}}$ has no zero eigenvector in $\K$. Indeed
$\gamma(\theta)>0$ except at $\theta=0$, a set of measure zero.
The theorem on strong convergence of spectral projections therefore gives
\[
 \mathbf1_{(-\infty,0)}(h_N)\longrightarrow E_c
 \quad\hbox{strongly on }\K.
\]
Thus all two-point functions with fixed finite support converge in the
convention \eqref{eq:covariance}, and Wick's rule gives convergence on every
local even CAR algebra, hence on every $B_I$. Spin-flip invariance makes
all local odd Pauli expectations zero. The full local limit is therefore
$\psi_{E_c}\circ\EE=\phi_c$, as claimed.
\end{proof}

Write $\omega_c=\omega_{E_c}$ and $\omega_\sharp=\omega_{E_\sharp}$. Set $\psi_c=\psi_{E_c}$, $\psi_\sharp=\psi_{E_\sharp}$, and let $\phi_c,\phi_\sharp$ denote their symmetric Pauli extensions. Purity on $B$ is \cite[Lemma~4.1]{ArakiEvans83}. Purity of the critical spin extension also follows from the sector calculation in Section~\ref{sec:sectors}.

\subsection{Ground-state selection and Kramers--Wannier duality}

\begin{proposition}[The critical quantum chain state]
\label{prop:groundstate}
The state $\phi_\sharp$ is the unique ground state of \eqref{eq:H1}. In particular it is pure on $A^P$.
\end{proposition}
This is the case $J=1/4$, $\gamma=1$, $\lambda=1$ of \cite[Theorem~1 and Proposition~3.3]{ArakiMatsui85}. The proof records the identification of the ground state with $\psi_\sharp\circ\EE$ in our conventions.
\begin{proof}
In the conventions \eqref{eq:JW}--\eqref{eq:covariance}, the one-particle symbol of $H_\lambda$ in \eqref{eq:Hlambda} is
\[
 h_\lambda(\theta)=(\cos\theta-\lambda)\sigma_3+\sin\theta\,\sigma_2,
\]
with eigenvalues $\pm\varepsilon_\lambda(\theta)$, $\varepsilon_\lambda(\theta)=(1-2\lambda\cos\theta+\lambda^2)^{1/2}$.
At $\lambda=1$ its eigenvalues are $\pm\varepsilon_\sharp(\theta)$, where $\varepsilon_\sharp(\theta)=2|\sin(\theta/2)|$, and
\[
 \frac{h_1(\theta)}{\varepsilon_\sharp(\theta)}
 =-|\sin(\theta/2)|\sigma_3
   +\sgn(\theta)\cos(\theta/2)\sigma_2
 =V(\Theta_\sharp(\theta))
\]
almost everywhere. Its negative-energy projection is therefore $E_\sharp$. Since $\varepsilon_\sharp>0$ almost everywhere, the one-particle generator has absolutely continuous spectrum and $0$ is not an eigenvalue. By \cite[Theorem~3(1)]{ArakiMatsui85}, the quasifree dynamics of $A^F$ then has a unique ground state, the Fock state in which the fields of one-particle vectors of negative energy annihilate the vacuum. In the convention \eqref{eq:covariance} this is $\omega_\sharp$. The restriction to $B$ of any ground state of the Pauli dynamics is a ground state of the common even dynamics \cite[Section~3]{ArakiMatsui85}. For absolutely continuous spectrum the even dynamics has a unique ground state, namely the restriction $\omega_\sharp|_B=\psi_\sharp$ \cite[Theorem~4(1)($\alpha$)]{ArakiMatsui85}. A $\Lambda$-invariant state is determined by its restriction to $B$, so the unique $\Lambda$-invariant ground state is $\psi_\sharp\circ\EE=\phi_\sharp$ \cite[Proposition~3.2]{ArakiMatsui85}. For $|\lambda|\ge1$ this state is pure and is the unique ground state of the Pauli dynamics \cite[Proposition~3.3(1)]{ArakiMatsui85}, see also \cite{Araki84}.
\end{proof}

Kramers--Wannier duality is the bounded-spread automorphism of $B$ defined by
\begin{equation}
 \kappa(\sigma_z^{(j)})=\sigma_x^{(j)}\sigma_x^{(j+1)},\qquad
 \kappa(\sigma_x^{(j)}\sigma_x^{(j+1)})=\sigma_z^{(j+1)}.
 \label{eq:KW}
\end{equation}
Its square is the site shift. In the CAR presentation it is the restriction of $\tau(W)$ \cite[Section~3]{EvansLewis86}, with
\[
 W(\theta)=\frac i2
 \begin{pmatrix}1-e^{i\theta}&1+e^{i\theta}\\-1-e^{i\theta}&e^{i\theta}-1\end{pmatrix},
 \qquad W^2=e^{i\theta}1.
\]
Direct matrix multiplication gives
\begin{equation}
 W^*V(\Theta)W=V(\pi-\Theta+\theta).
 \label{eq:KWangle}
\end{equation}
Consequently $\psi_\sharp\circ\kappa=\psi_\sharp$, while KW sends the angle $\Theta_c=\Theta_\sharp+\chi_z$ to $\Theta_\sharp-\chi_z$. For $z>0$ these define distinct states. In fact their Fock representations are inequivalent. The difference of their basis projections is a nonzero multiplication operator, hence is not Hilbert--Schmidt \cite{Araki70}. Equivalence of their even vacuum representations would imply equivalence of the induced full CAR representations, so $\psi_c$ is not $\kappa$-covariant either.

\begin{remark}[Invariance does not select a ground state]
Both $E_\sharp$ and $1-E_\sharp$ are basis projections fixed by $W$, and $1-E_\sharp=E(\Theta)$ for the odd angle $\Theta=\theta/2-(\pi/2)\sgn\theta$. Thus oddness and KW invariance do not characterize $\phi_\sharp$. Proposition~\ref{prop:groundstate} selects it by the sign of the one-particle energy. KW invariance can also be checked from the full covariance in Section~\ref{sec:TL}. Equal expectations of the generators alone would not suffice.
\end{remark}

\section{Half-line duality for pure quasifree states}
\label{sec:duality}

For a pure state $\psi$ on $B$, bounded-spread cone Haag duality means that some $R\ge0$ satisfies
\[
 \pi_\psi(B_{L^c})'\subseteq\pi_\psi(B_{L^{(R)}})''
 \quad\text{in }\mathcal B(L^2(B,\psi))
\]
for every ray $L$. Here $L^{(R)}$ is the $R$-neighbourhood of $L$ in the site lattice. Since $B$ is simple, this condition implies that $H_\psi$ is local \cite[Theorem~5.6]{EvansJones2025}. For pure quasifree states the stronger equality with $R=0$ follows directly from twisted duality.

\begin{theorem}[Observable duality]
\label{thm:duality}
Let $\omega_E$ be any pure quasifree state of the CAR chain and $\psi_E=\omega_E|_B$. Then, for every ray of physical sites $L$,
\[
 \pi_{\psi_E}(B_{L^c})'=\pi_{\psi_E}(B_L)''
 \quad\text{on }L^2(B,\psi_E).
\]
\end{theorem}
\begin{proof}
Let $(\pi_E,\HH_E,\Omega_E)$ be the Fock representation of $\omega_E$, and $P$ the self-adjoint unitary implementing the grading with $P\Omega_E=\Omega_E$. The vacuum is cyclic for $\pi_E(B)$ in the even subspace $\HH_E^+=\ker(P-1)$, so $\HH_E^+=L^2(B,\psi_E)$. For a set $S$ of Majorana indices let $A^F(S)$ be the CAR algebra on those indices. In the full Fock representation put
\[
 M_S=\pi_E(A^F(S))'',\qquad M_{S,+}=M_S\cap\{P\}',\qquad p_+=\frac{1+P}{2}.
\]
Parity averaging and bounded strong approximation give $M_{S,+}=\pi_E(A^F(S)_+)''$, linking these fixed-point von Neumann algebras to the local observable algebras.

Purity provides a complex one-particle Hilbert space whose underlying real Hilbert space is the real Majorana space, with its original real inner product. The closed real field subspaces associated to $S$ and $S^c$ are therefore real orthogonal complements. Foit's abstract twisted duality theorem \cite[Theorem~2.4(v)]{Foit83} gives
\begin{equation}
 M_S'=Z M_{S^c} Z^*,
 \label{eq:twistedduality}
\end{equation}
where $Z$ is the Klein twist, a unitary function of $P$. No standardness assumption is required.

We explicitly pass to the even vacuum sector. The commutant identity for an intersection of von Neumann algebras gives
\[
 M_{S,+}'=(M_S\cap\{P\}')'
 =M_S'\vee\{P\}''=M_{S^c}\vee\{P\}''.
\]
Because $p_+$ commutes with $M_{S,+}$, its reduced commutant is
\begin{equation}
 \bigl(M_{S,+}|_{\HH_E^+}\bigr)'
 =p_+\bigl(M_{S^c}\vee\{P\}''\bigr)p_+\big|_{\HH_E^+}.
 \label{eq:compression}
\end{equation}
The algebra $M_{S^c}$ is normalized by $P$. Hence the linear span of $x+yP$, with $x,y\in M_{S^c}$, is strongly dense in $M_{S^c}\vee\{P\}''$. Moreover,
\[
 p_+xp_+=p_+\frac{x+PxP}{2}p_+.
\]
Thus the right-hand side of \eqref{eq:compression} equals $M_{S^c,+}|_{\HH_E^+}$, which is weakly closed. Taking $S$ to contain both Majoranas of each site of a ray proves the result, and interchanging the ray and its complement gives the displayed orientation.
\end{proof}

\begin{corollary}
\label{cor:criticalduality}
Both $\psi_c$ and $\psi_\sharp$ satisfy cone Haag duality with $R=0$. Their algebra objects are local.
\end{corollary}

\begin{remark}[Standardness and region conventions]
The proof applies to any set $S$ if its observable algebra is defined as the even part of the CAR algebra on the whole set. For a disconnected set, this includes products of odd fields in different components and can exceed the algebra generated by the even algebras of the separate components. We use only physical rays, so no disconnected-region identification is needed.

Generic position of the one-particle projections is not an input to Theorem~\ref{thm:duality}. It has a different consequence. The half-line field subspace is standard, and the Fock vacuum is cyclic and separating for its field algebra. We prove this additional property for the critical states in Appendix~\ref{sec:standardness}. Neither standardness nor global purity determines the type of a half-line von Neumann algebra. Factoriality of the critical half-line algebras is proved separately in Proposition~\ref{prop:factor}.
\end{remark}

\section{A direct calculation of the algebra object}
\label{sec:sectors}

For a basis projection $E$, write $\pi_E^+$ for the representation of $B$ on the even Fock space $\HH_E^+=L^2(B,\psi_E)$ of the proof of Theorem~\ref{thm:duality}. In this section, implementation of a sector always means implementation in $\pi_E^+$, and implementation of a full CAR automorphism means implementation on $\HH_E$.

\subsection{Sector representatives and parity}

Let $Q_+$ and $Q_-$ be the projections on the sites $j\ge1$ and $j\le0$ in $\K$, including both self-dual components. Define
\[
 \theta_-=Q_+-Q_-,\qquad \Lambda_-=\tau(\theta_-).
\]
Put $q=c_1+c_1^*$, an odd self-adjoint unitary at the right endpoint of this cut.
The four sectors are represented by
\begin{equation}
 \rho_1=\mathrm{id},\qquad
 \rho_e=\Lambda_-|_B,\qquad
 \rho_f=\operatorname{Ad}(q)|_B,\qquad
 \rho_m=\rho_f\rho_e,\qquad f=em.
 \label{eq:sector-representatives}
\end{equation}
Here $e$ is the disorder sector and $m$ is the spin-charge sector.
With these labels the high-temperature endpoint state of
Corollary~\ref{cor:offcritical-sector-calculation} has algebra object
$\mathbf1\oplus e$. Thus our labels interchange $e$ and $m$ relative to
\cite[Example~5.28]{EvansJones2025}.

For completeness, the automorphisms in
\eqref{eq:sector-representatives} are localized and transportable for
the interval net $B_I$.  Conjugation by $q$ is the identity on every
even local algebra away from site $1$.  The half-line parity is the
identity on each even local algebra lying wholly on either side of
the cut.  Moving the Majorana $q=a_1$ to another Majorana $m_k$ is
implemented by the even unitary $m_kq\in B$, since
$\Ad(m_k)=\Ad(m_kq)\circ\Ad(q)$. Moving the cut through finitely many
complete sites is implemented by the product of their on-site parity
operators $\sigma_z^{(j)}=i\,a_jb_j$.  For an invertible sector $\rho$ the self-braiding $\varepsilon(\rho,\rho)\in\End(\rho^2)$ is a scalar, independent of the choice of localized transporters. With the transporters above one finds $\varepsilon(\rho_f,\rho_f)=-1$ and $\varepsilon(\rho_e,\rho_e)=+1$.
The half-line parity is outer. It changes the sign of
$i m_jm_k$ with $j$ arbitrarily far to the left and $k$ arbitrarily
far to the right, while any fixed local even operator commutes with
these bilinears.  Approximation by local operators therefore rules
out an inner implementer in $B$.
Consequently these are the four distinct invertible sectors of
$\DHR(B)$. The remaining nontrivial sector $m=ef$ is the other boson.
In the crossed-product Jordan--Wigner realization of \cite[Section~2]{ArakiEvans83}, following \cite{Araki84}, see also \cite[Section~2]{EvansLewis86}, an auxiliary self-adjoint unitary $T$ implements $\Lambda_-$ and $\sigma_x^{(1)}=Tq$. Hence
\begin{equation}
 \operatorname{Ad}(\sigma_x^{(1)})|_B=\rho_m,
 \qquad A^P=B\rtimes_{\rho_m}\mathbb Z_2.
 \label{eq:spin-crossed-product}
\end{equation}

\begin{proposition}[The multiplicities]
\label{prop:sector-multiplicities}
Let $E$ be a basis projection and $X_E=L^2(B,\psi_E)$. Then
\[
 H_{\psi_E}=\underline{\operatorname{End}}(X_E)
 \cong\bigoplus_{a\in\{1,e,m,f\}}N_a\,a,
 \qquad N_a=\dim\operatorname{Hom}_B(\pi_E^+\circ\rho_a,\pi_E^+)\in\{0,1\},
\]
and $N_a=1$ exactly when $\psi_E\circ\rho_a$ and $\psi_E$ are equivalent. Write $\sigma(E,F)=(-1)^{\dim(E\wedge(1-F))}$ for the $\Z_2$-index of basis projections with $E-F$ Hilbert--Schmidt \cite[(5.3)]{ArakiEvans83}. Then
\begin{enumerate}[leftmargin=*]
\item $N_f=0$.
\item $N_e=1$ if and only if $E-\theta_-E\theta_-$ is Hilbert--Schmidt and $\sigma(E,\theta_-E\theta_-)=1$.
\item $N_m=1$ if and only if $E-\theta_-E\theta_-$ is Hilbert--Schmidt and $\sigma(E,\theta_-E\theta_-)=-1$.
\end{enumerate}
\end{proposition}
\begin{proof}
The multiplicity formula is the definition of the internal end, with sector action realized by composition with $\rho_a$. All four representatives are involutive, so the possible inverse in the bimodule convention makes no difference. Since $\rho_a$ is an automorphism and $\pi_E^+$ is irreducible \cite[Lemma~4.1]{ArakiEvans83}, $\pi_E^+\circ\rho_a$ is the irreducible GNS representation of $\psi_E\circ\rho_a$, and Schur's lemma gives the remaining assertions. By \eqref{eq:pullback}, $\psi_E\circ\rho_e$ is the even restriction of $\omega_{\theta_-E\theta_-}$. In the notation of \cite[Section~4]{ArakiEvans83} take $\alpha$ the grading, $\beta=\Lambda_-$ and $U=q$, so that $\psi_E\circ\rho_f=\psi_E\circ\Ad U$ and, since $\Lambda_-(q)=q$, $\psi_E\circ\rho_m=\psi_E\circ\beta\circ\Ad U$. Then (1) is \cite[Lemma~4.2]{ArakiEvans83}, (2) is \cite[Theorem~4]{ArakiEvans83}, and (3) is \cite[Lemma~4.3]{ArakiEvans83} combined with \cite[Theorem~4]{ArakiEvans83} and the Hilbert--Schmidt criterion for equivalence of Fock states \cite{Araki70}.
\end{proof}

By \cite{Araki70}, $\Lambda_-$ is implemented on $\HH_E$ exactly when $[\theta_-,E]\in\HS$, and then by a unitary of parity $\sigma(E,\theta_-E\theta_-)$, unique up to a scalar \cite[Section~7]{ArakiEvans83}, see also \cite[(5.16)--(5.17)]{EvansLewis86}.

\subsection{The critical states}

\begin{lemma}[The critical half-line commutator]
\label{lem:projection-jump}
For $E=E_c$ and $E=E_\sharp$, the operator $E-\theta_-E\theta_-$, equivalently $[\theta_-,E]$, is not Hilbert--Schmidt.
\end{lemma}
\begin{proof}
For $E_c$ this is \cite[Lemma~5.1]{ArakiEvans83}. It states that $\omega_{E_c}$ and $\omega_{\theta_-E_c\theta_-}$ are inequivalent, which is proved in \cite[Section~6(2)]{ArakiEvans83} from the jump of the symbol at $\theta=0$. For $E_\sharp$ it is the case $|\lambda|=1$, $\gamma\ne0$ of \cite[Lemma~4.5(1)]{ArakiMatsui85}, computed there in Section~7.
\end{proof}

\begin{theorem}[The critical algebra objects]
\label{thm:critical-algebra-direct}
The two critical states satisfy
\[
 H_{\phi_\sharp}\cong\mathbf1,
 \qquad H_{\phi_c}\cong\mathbf1
\]
as algebra objects in $\DHR(B)$. In both cases every nontrivial sector has multiplicity zero.
\end{theorem}
\begin{proof}
By Proposition~\ref{prop:sector-multiplicities} and Lemma~\ref{lem:projection-jump}, $N_e=N_m=N_f=0$. The only surviving summand is the unit with multiplicity one, and its unital algebra structure is that of the unit algebra.
\end{proof}

In terms of the multiplicities of Proposition~\ref{prop:sector-multiplicities}, $\dim\End_{A^P}L^2(A^P,\phi_E)=1+N_m$, by induction along $A^P\cong\mathbf1\oplus m$ as a $B$--$B$ bimodule \eqref{eq:spin-crossed-product} and Frobenius reciprocity, see \cite[Theorem~5.12]{EvansJones2025}. So $N_m=0$ is purity of the symmetric Pauli state. For $\phi_c$ this is the critical case of \cite[Theorem~2 and Lemma~5.1]{ArakiEvans83}, and for $\phi_\sharp$ it is \cite[Proposition~3.3(1)]{ArakiMatsui85}.

For $\psi_\sharp$ there is a second categorical route. Corollary~\ref{cor:criticalduality} gives locality and KW invariance permutes $e$ and $m$. Among the three connected commutative algebras $\one$, $\one\oplus e$, and $\one\oplus m$, only $\one$ is fixed. This is the mechanism of \cite[Example~5.28]{EvansJones2025}. The calculation above also handles $\psi_c$, which is not KW invariant, and uses neither transport of DHR sectors along an automorphism with decaying tails nor a classification of Morita classes.

\subsection{Half-line factors}

The critical half-line commutator also determines the centres of the half-line algebras. For a physical ray $L$, let $S_L$ be the set of Majorana indices of its sites and, as in the proof of Theorem~\ref{thm:duality},
\[
 M_L=\pi_E(A^F(S_L))'',\qquad M_{L,+}=M_L\cap\{P\}'.
\]
Write $\theta_L$ for the self-adjoint unitary on $\K$ equal to $1$ on the sites of $L$ and $-1$ on the sites of $L^c$. Then $\theta_L=\pm\theta_-$ up to a site translation, so Lemma~\ref{lem:projection-jump} and translation invariance give $[\theta_L,E]\notin\HS$ for $E\in\{E_\sharp,E_c\}$.

\begin{proposition}[Half-line factors]
\label{prop:factor}
Let $E\in\{E_\sharp,E_c\}$ and let $L$ be a physical ray. Then $M_L$, $M_{L,+}$ and $\pi_E^+(B_L)''$ are factors.
\end{proposition}
\begin{proof}
We use graded locality. Even elements of $M_L$ commute with $M_{L^c}$, and odd elements of $M_L$ anticommute with odd elements of $M_{L^c}$. Also $M_L\vee M_{L^c}=\pi_E(A^F)''=\mathcal B(\HH_E)$. If a unitary $u$ on $\HH_E$ commutes with $M_L$ and satisfies $uxu^*=PxP$ for $x\in M_{L^c}$, then $\Ad u$ implements $\tau(\theta_L)$. This is impossible, since $[\theta_L,E]\notin\HS$ \cite{Araki70}. The same holds with the roles of $L$ and $L^c$ exchanged, since $\theta_{L^c}=-\theta_L$.

\emph{$M_L$.} Since $\Ad P$ preserves $M_L$, it preserves its centre, and every central element is a sum of an even and an odd central element. An even central element commutes with $M_L$ and, by graded locality, with $M_{L^c}$, so it is a scalar. If $\gamma$ is a nonzero odd central element, then $\gamma^*\gamma$ is even and central, hence a positive scalar, and we may take $\gamma$ unitary. It commutes with $M_L$, commutes with the even part of $M_{L^c}$ and anticommutes with its odd part. So it acts trivially on the fields of $L$ and by the grading on those of $L^c$, and would implement $\tau(\theta_L)$, which is excluded by the first paragraph. Hence $M_L$ is a factor.

\emph{$M_{L,+}$.} Choose an odd self-adjoint unitary $a\in\pi_E(A^F(S_L))$, for instance a Majorana operator. Then $M_L=M_{L,+}+M_{L,+}a$, and $\sigma=\Ad a$ restricts to an involutive automorphism of the centre $\mathcal Z$ of $M_{L,+}$. An element of $\mathcal Z$ fixed by $\sigma$ commutes with $M_{L,+}$ and $a$, so it lies in the centre of $M_L$ and is a scalar. If $y\in\mathcal Z$ satisfies $\sigma(y)=-y$, then $y^*y$ is $\sigma$-fixed, hence a scalar. For two such elements $y_1,y_2$ also $y_1y_2^*$ is a scalar. Thus $\mathcal Z\subseteq\mathbb C1+\mathbb Cu$ for a unitary $u$ with $aua=-u$, which we may take self-adjoint since $u^2$ is a scalar. Such a $u$ is even, central in $M_{L,+}$ and anticommutes with $a$, hence with every odd element of $M_L$. By graded locality it commutes with $M_{L^c}$. Thus $\Ad u$ is the identity on $M_{L^c}$ and the grading on $M_L$, which is excluded by the first paragraph with $L$ and $L^c$ exchanged. Hence $\mathcal Z=\mathbb C1$.

\emph{$\pi_E^+(B_L)''$.} By the proof of Theorem~\ref{thm:duality}, $M_{L,+}=\pi_E(B_L)''$, and $p_+=(1+P)/2$ lies in its commutant. Therefore $\pi_E^+(B_L)''=M_{L,+}p_+$, and $x\mapsto xp_+$ is a normal $*$-homomorphism from $M_{L,+}$ onto this algebra. Its kernel is a weakly closed two-sided ideal, hence of the form $M_{L,+}(1-z)$ with $z$ a central projection. Since $M_{L,+}$ is a factor and $p_+\ne0$, the kernel is zero. So $\pi_E^+(B_L)''\cong M_{L,+}$ is a factor.
\end{proof}

Away from criticality the half-line twist is implementable. By the proof of Corollary~\ref{cor:offcritical-sector-calculation} and translation invariance, $\tau(\theta_L)$ has an implementer $U$ on $\HH_E$, unique up to a scalar, which is even for $\beta<\beta_c$ and odd for $\beta>\beta_c$. Factoriality then fails in exactly one of the two half-line algebras, according to the parity of $U$.

If $\beta<\beta_c$, then $N_e=1$, and $U$ restricts to an implementer $U_+$ of $\tau(\theta_L)|_B$ on $\HH_E^+$. It commutes with $\pi_E^+(B_L)$ and $\pi_E^+(B_{L^c})$, because $\tau(\theta_L)$ is trivial on the even algebra of either half-line. By Theorem~\ref{thm:duality}, $U_+$ lies in the centre of $\pi_E^+(B_L)''$. It is not a scalar, since $\rho_e$ is not the identity and $\pi_E^+$ is faithful on the simple algebra $B$. Hence $\pi_E^+(B_L)''$, and with it $M_{L,+}$, is not a factor. On the other hand $M_L$ is a factor. In the proof of Proposition~\ref{prop:factor}, a nonscalar central element of $M_L$ would give an odd implementer of $\tau(\theta_L)$, and there is none since $U$ is even.

If $\beta>\beta_c$, then $N_m=1$ and $U$ is odd. It commutes with $M_L$, and $PU$ commutes with $M_{L^c}$. By twisted duality \eqref{eq:twistedduality}, $PU\in M_{L^c}'=ZM_LZ^*$. For odd $x$, $Z^*xZ$ is a scalar multiple of $Px$, since $Z$ is the Klein twist. Hence $U\in M_L$. Thus $U$ is a nonscalar odd element of the centre of $M_L$, and $M_L$ is not a factor. By the argument for $M_L$ in the proof of Proposition~\ref{prop:factor}, the centre of $M_L$ is $\mathbb C1+\mathbb CU$, whose even part is $\mathbb C1$. In the argument for $M_{L,+}$, a nonscalar central element $u$ would be an even implementer of $\tau(\theta_{L^c})=\Theta\circ\tau(\theta_L)$, where $\Theta=\tau(-1)$ is the grading, whose implementers are the multiples of the odd unitary $PU$. Hence $M_{L,+}$ is a factor, and so is $\pi_E^+(B_L)''\cong M_{L,+}$.

So for $\beta<\beta_c$ the half-line field algebra is a factor and the half-line observable algebra is not. For $\beta>\beta_c$ the reverse holds. At $\beta_c$ both are factors (Proposition~\ref{prop:factor}). The behaviour of the field algebra is consistent with the result of Lewis and Winnink that the restriction of the state to the half-line CAR algebra is primary for $\beta<\beta_c$ and not primary for $\beta>\beta_c$ \cite{LW81}, as recalled in \cite[Section~1]{ArakiEvans83}.

\subsection{Off-critical calibration}

Fix ferromagnetic couplings $J_1,J_2>0$. For $0<\beta<\infty$ let $\phi_\beta$ be the Araki--Evans state on the Pauli algebra $A^P$ of Section~\ref{sec:models}. This is the limit of the Perron--Frobenius vector states of the row transfer matrices \eqref{eq:transfer}, with $K_2=\beta J_2$ and $\tanh K_1^*=e^{-2\beta J_1}$. Put $a=K_2$ and $b=K_1^*$. By \cite[Lemma~3.2]{ArakiEvans83}, $\phi_\beta=\psi_{E_\beta}\circ\EE$ for the basis projection $E_\beta$ on $\K$ with symbol
\begin{equation}
 U_{a,b}(\theta)=e^{-a n(\theta)}e^{2b\sigma_3}e^{-a n(\theta)},
 \qquad E_\beta(\theta)=E_{a,b}(\theta)=1_{(-\infty,0)}(-\log U_{a,b}(\theta)),
 \label{eq:offcritical-symbol}
\end{equation}
where $n(\theta)=\cos\theta\,\sigma_3+\sin\theta\,\sigma_2$. Here $U_{a,b}$ is the symbol of the one-particle transfer operator $e^{2H}$ of \cite[(3.5)]{ArakiEvans83}, in the convention of Proposition~\ref{prop:transferconvention}, and $E_\beta=1_{(0,\infty)}(H)$ is the projection onto the vectors $h$ for which $C(h)$ annihilates the vacuum, as in \eqref{eq:covariance}. It is exactly the projection $E_\beta=(1-V_\beta)/2$ of \cite[(2.18)--(2.20)]{EvansLewis86}, equivalently \cite[(3.7)--(3.9)]{ArakiEvans83}. Araki and Evans label the same Fock state by the complementary projection $E_-=1-E_\beta$, the spectral projection of $H$ for $(-\infty,0)$. In both descriptions the annihilation operators are the $C(h)$ with $h$ in the range of $E_\beta$ \cite[(3.4) and Lemma~3.2]{ArakiEvans83}. The states $\phi_0$ and $\phi_\infty$ are the limits $\beta\to0$ and $\beta\to\infty$, with $E_0=(1+\sigma_3)/2$ and $E_\infty=(1-n(\theta))/2$. The unique critical inverse temperature satisfies $a=b$, and $\phi_{\beta_c}=\phi_c$.
Here $\beta$ is the inverse temperature of the two-dimensional
classical model, and $\phi_\beta$ is a transfer matrix state, not a
finite-temperature Gibbs state of $H_1$.

\begin{corollary}[The two gapped regimes]
\label{cor:offcritical-sector-calculation}
For the Araki--Evans states $\phi_\beta$, $0\le\beta\le\infty$,
\[
 H_{\phi_\beta}\cong
 \begin{cases}
  \mathbf1\oplus e,&0\leq\beta<\beta_c,\\
  \mathbf1,&\beta=\beta_c,\\
  \mathbf1\oplus m,&\beta_c<\beta\leq\infty.
 \end{cases}
\]
\end{corollary}
\begin{proof}
By \cite[Lemma~5.1]{ArakiEvans83}, $E_\beta-\theta_-E_\beta\theta_-$ is Hilbert--Schmidt for $\beta\ne\beta_c$, and by \cite[Lemma~5.2]{ArakiEvans83}, $\sigma(E_\beta,\theta_-E_\beta\theta_-)$ is $1$ for $\beta<\beta_c$ and $-1$ for $\beta>\beta_c$. The endpoints $\beta=0,\infty$ are the cases $K_2=0$ and $K_1^*=0$ of \cite[Section~6(3)]{ArakiEvans83}. Araki and Evans state these results for the complementary projection $1-E_\beta$. Since $\mathfrak c\{E\wedge(1-F)\}\mathfrak c=(1-E)\wedge F$, the index is the same \cite[Section~5]{ArakiEvans83}. Proposition~\ref{prop:sector-multiplicities} now gives $N_e=1$, $N_m=0$ for $\beta<\beta_c$, and $N_e=0$, $N_m=1$ for $\beta>\beta_c$, with $N_f=0$ throughout. The critical case is Theorem~\ref{thm:critical-algebra-direct}. On $\mathbf1\oplus e$ and $\mathbf1\oplus m$ the connected $\mathbb Z_2$ algebra structure is unique up to isomorphism, by rescaling its generating intertwiner to have square $1$.
\end{proof}

\begin{center}
\begin{tabular}{lll}
\toprule
Regime & Algebra object & Symmetric spin state\\
\midrule
$\beta<\beta_c$ & $\one\oplus e$ & pure, disordered\\
$\beta=\beta_c$ & $\one$ & pure, critical\\
$\beta>\beta_c$ & $\one\oplus m$ & mixture of two ordered pure states\\
\bottomrule
\end{tabular}
\end{center}
With the labels \eqref{eq:sector-representatives}, the first and third rows recover the noncritical algebra objects of \cite[Example~5.28]{EvansJones2025}, where $e$ and $m$ are interchanged. They are topological in the categorical sense, and the critical row is categorically gapless. The purity column is \cite[Theorems~1 and~2]{ArakiEvans83}. No completeness theorem for automorphic equivalence of states is assumed.

\section{Conjugacy and the two dynamics}
\label{sec:conjugacy}

Fix $z\in[0,1)$ and use the critical angles of \eqref{eq:angles}--\eqref{eq:criticalangles}. The difference $\chi_z$ is analytic. Removing it is a rotation of the covariance, but a locality statement on the Pauli algebra requires a spin extension. Its existence is a theorem of Araki. A quasifree dynamics of the CAR algebra whose symbol is holomorphic in $e^{i\theta}$ near the unit circle extends uniquely to a strongly continuous one-parameter group of automorphisms of the Pauli algebra \cite[Proposition~5.1]{Araki86}, by the Evans--Lewis criterion \cite{EvansLewis86}. What we add is a uniform exponential locality bound, for which we realize the extension as the evolution of an even interaction.

For $z=0$, $w_t=1$ and $\alpha_t=\id$. Hence assume $0<z<1$ and abbreviate $\chi_z$ to $\chi$.

\begin{proposition}[A spin extension with uniform exponential locality]\label{prop:spinlocality}
Fix $0<z<1$, write $a=z^2$, and set
\[
 \chi(\theta)=\operatorname{Im}\log(1-ae^{i\theta})
      =-\sum_{n\ge1}\frac{a^n}{n}\sin(n\theta),
 \quad d(\theta)=\frac12\chi(\theta)\sigma_1,
 \quad w_t=e^{itd}.
\]
There is a unique strongly continuous group $\alpha_t$ of automorphisms of the Pauli algebra whose restriction
to $B$ is $\tau(w_t)|_B$ \cite[Proposition~5.1]{Araki86}. It is graded and
translation covariant, and for every $0<\mu<-\log a$ and every $T<\infty$
there is $C_{\mu,T}<\infty$ such that, for every finite interval
$I\subset\mathbb Z$, every $r\in\mathbb N_0$, $|t|\le T$, and
$A\in A^P_I$,
\[
 \operatorname{dist}\bigl(\alpha_t(A),A^P_{I^{(r)}}\bigr)
 \le C_{\mu,T}\|A\|e^{-\mu r}.
\]
The constant is independent of the position and length of $I$. This
estimate applies equally to $\alpha_t^{-1}=\alpha_{-t}$.
\end{proposition}
\begin{proof}
The symbol $\chi(\theta)=\frac1{2i}\bigl(\log(1-az)-\log(1-a/z)\bigr)$, $z=e^{i\theta}$, is holomorphic in the annulus $a<|z|<1/a$, so \cite[Proposition~5.1]{Araki86} gives existence and uniqueness of $\alpha_t$. We now construct $\alpha_t$ as the dynamics of an even interaction. By uniqueness this is the same group, and the construction gives the grading, the covariance and the locality bound. The self-adjoint matrix $d$ satisfies $\mathfrak c d\mathfrak c=-d$, so $w_t$
commutes with $\mathfrak c$. Its spatial matrix coefficients obey
$\|d_{j,k}\|\le C a^{|j-k|}/|j-k|$ for $j\ne k$, and $d_{j,j}=0$.
In a finite interval $\Lambda$, put
\[
 Q_\Lambda=\frac12\sum_{u,v\in\Lambda\times\{1,2\}}
        (d_\Lambda)_{uv} C(e_u)C(e_v)^*,
 \qquad d_\Lambda=P_\Lambda dP_\Lambda.
\]
The CAR imply $[Q_\Lambda,C(f)]=C(d_\Lambda f)$, so its Heisenberg
evolution is $\tau(e^{itd_\Lambda})$. Group the quadratic terms on two sites by their endpoint interval,
and apply the Jordan--Wigner transformation to these even terms. This
gives a translation invariant even spin interaction $\Phi$
supported on intervals, with
\[
 \|\Phi([j,j+n])\|\le C\frac{a^n}{n},\qquad n\ge1.
\]
Indeed a quadratic fermion term with endpoints $j,j+n$ becomes a Pauli
string contained in $[j,j+n]$, with unchanged norm. For any
$0<\nu<-\log a$,
\[
 \sup_x\sum_{X\ni x}|X|e^{\nu\operatorname{diam}X}\|\Phi(X)\|
 \le C\sum_{n\ge1}(n+1)^2e^{\nu n}\frac{a^n}{n}<\infty.
\]
The thermodynamic-limit spin dynamics therefore exists and satisfies
an exponential Lieb--Robinson bound, see \cite[Theorem 3.1 and
Section 3.3]{NSY19}. The finite volume spin and fermion evolutions agree
on the common even algebra. Passing to the limit shows
$\alpha_t|_B=\tau(w_t)|_B$. Evenness and translation covariance pass
to the limit as well.

To make the stated uniformity explicit, choose $\mu<\nu<-\log a$ and
apply the Lieb--Robinson estimate with an exponentially decreasing
$F$-function. For $A\in A^P_I$ and local
$D\in A^P_{(I^{(r)})^c}$ it gives, uniformly for $|t|\le T$,
\[
 \|[\alpha_t(A),D]\|
 \le C\|A\|\|D\|
      \sum_{x\in I}\sum_{y\notin I^{(r)}}e^{-\mu|x-y|}
 \le \frac{2Ce^{-\mu(r+1)}}{(1-e^{-\mu})^2}\|A\|\|D\|.
\]
The last estimate is a geometric sum over the two endpoints of $I$,
and has no dependence on $|I|$. Applying the product-trace conditional
expectation onto $A^P_{I^{(r)}}$, as in
\cite[Corollary 4.4]{NSY19}, proves the claimed distance bound.
The construction works for positive and negative times, proving the
same estimate for the inverse.
\end{proof}

Without the locality bounds, the extension $\alpha_1$, and hence the equivalence of Corollary~\ref{cor:conjugacy}, is therefore already given by Araki's construction \cite[Section~3 and Proposition~5.1]{Araki86}.

\begin{corollary}[Conjugacy and the transported duality]\label{cor:conjugacy}
Set $w=w_1$ and $\alpha=\alpha_1$. Then
\[
 \phi_\sharp=\phi_c\circ\alpha,
 \qquad
 \kappa_c=\alpha|_B\circ\kappa\circ(\alpha|_B)^{-1}
 \quad\hbox{satisfies}\quad \psi_c\circ\kappa_c=\psi_c.
\]
In particular,
\[
 \kappa_c=\tau(w)^2\circ\kappa
          =\kappa\circ\tau(w)^{-2}\quad\hbox{on }B.
\]
\end{corollary}
\begin{proof}
Writing $V(\Theta)=R(\Theta)^*\sigma_3R(\Theta)$ with
$R(s)=e^{-is\sigma_1/2}$ and $\Theta_c=\Theta_\sharp+\chi$, one obtains
\[
 w^*V_cw=V_\sharp,
 \qquad \omega_c\circ\tau(w)=\omega_\sharp.
\]
The symmetric spin extension of a state on $B$ is its composition
with the grading expectation. Since $\alpha$ is graded, equality on
$B$ gives the asserted equality on the Pauli algebra. The transported
duality fixes $\psi_c$ because $\psi_\sharp$ is $\kappa$-invariant.
The matrix $W$ anticommutes with $\sigma_1$, hence
$WwW^{-1}=w^{-1}$ and
$\kappa\tau(w)\kappa^{-1}=\tau(w)^{-1}$, yielding the final formula.
\end{proof}

\begin{remark}
Exponential locality of a CAR multiplier alone does not establish
exponential locality of an arbitrary spin extension. The interaction
construction above supplies that missing step, including uniform control
of the inverse. It does not establish bounded spread or justify
transporting the categorical invariant using a theorem that assumes
bounded spread. The invariant should be calculated independently for
the two states.
\end{remark}

\begin{remark}[Anisotropic limit and continuum limit]
\label{rem:anisotropic}
The limit $z\to0$ removes the analytic factor $\chi$ and gives
$E_c\to E_\sharp$ in operator norm. It is the anisotropic Hamiltonian
limit of the transfer matrix family. A continuum scaling limit also
rescales space and time and is a distinct construction. The covariance
calculation and the automorphism $\alpha$ do not by themselves construct
that scaling limit or prove preservation of central charge. It is
therefore preferable to describe the conformal interpretation as the
established Ising/free Majorana interpretation of the critical chain,
and not as a new theorem proved by this conjugacy.
\end{remark}

\subsection{Comparison of the dispersions}

The transfer matrix dispersion is given by $\cosh\gamma=\cosh2K_1^*\cosh2K_2-\sinh2K_1^*\sinh2K_2\cos\theta$ \cite[(3.8)]{ArakiEvans83}, \cite[(2.18)]{EvansLewis86}, going back to Onsager and Kaufman. At criticality it obeys
\[
 \cosh\gamma(\theta)=\cosh^2(2K)-\sinh^2(2K)\cos\theta,
 \qquad z=\tanh K.
\]
Consequently,
\begin{equation}
 \sinh\bigl(\gamma(\theta)/2\bigr)
   =\sinh(2K)|\sin(\theta/2)|.
 \label{eq:dispersion}
\end{equation}
The normalized one-particle symbols are $V_c$ and $V_\sharp$. Conjugation by $w=e^{i\chi\sigma_1/2}$ gives
\[
 w^*\bigl(\gamma V_c\bigr)w=\gamma V_\sharp,
\]
whereas the critical quantum chain Hamiltonian has symbol $\varepsilon_\sharp V_\sharp$ with $\varepsilon_\sharp=2|\sin(\theta/2)|$. Thus the conjugacy matches ground state projections, but does not intertwine these two specified dynamics. Their velocities at low momentum are $\sinh(2K)$ and $1$, respectively.

As $K\to0$, $\gamma(\theta)=4K|\sin(\theta/2)|+O(K^3)$ uniformly in $\theta$, so rescaling the transfer matrix generator by $1/(2K)$ recovers the quantum chain dispersion. This agrees with $E_c\to E_\sharp$ as $z\to0$. It is a limit within the family of lattice models, not a construction of a continuum theory.

\section{The Temperley--Lieb state and the conformal interpretation}
\label{sec:TL}

\subsection{The complete lattice covariance}

The common even algebra is generated by
\begin{equation}
 e_{2j-1}=\frac12(1+\sigma_z^{(j)}),\qquad
 e_{2j}=\frac12(1+\sigma_x^{(j)}\sigma_x^{(j+1)}).
 \label{eq:TLgenerators}
\end{equation}
They satisfy the projection-normalized Temperley--Lieb relations
\[
 e_i^2=e_i=e_i^*,\qquad e_ie_{i\pm1}e_i=\frac12e_i,
 \qquad [e_i,e_j]=0\quad(|i-j|\ge2).
\]
In Majorana notation these are uniformly $e_k=(1+i m_km_{k+1})/2$. Thus $H_1=-\sum_k e_k$ up to its formal additive constant, and KW is the shift $e_k\mapsto e_{k+1}$. The Majorana translation realizes this shift on the even algebra \cite[Remark~3.2]{EvansLewis86}.

For a translation invariant Majorana state define its real antisymmetric covariance by
\begin{equation}
 \Omega(n)=\frac i2\omega([m_k,m_{k+n}]).
 \label{eq:OmegaDef}
\end{equation}
For $n\ne0$ this equals $\omega(i m_km_{k+n})$, whereas $\Omega(0)=0$. The uncentered expression $\omega(i m_k^2)$ is $i$ and is not the diagonal covariance.

\begin{proposition}
\label{prop:covariance}
For $\omega_\sharp$, the covariance is independent of the starting index and is
\begin{equation}
 \Omega(n)=\begin{cases}2/(\pi n),&n\text{ odd},\\0,&n\text{ even}.
 \end{cases}
 \label{eq:Omega}
\end{equation}
In particular $\omega_\sharp$ is invariant under Majorana translation and $\psi_\sharp$ is KW invariant.
\end{proposition}
\begin{proof}
Use the Fourier convention
\[
 G_n=\frac1{2\pi}\int_{-\pi}^{\pi}e^{i\Theta_\sharp(\theta)}e^{-in\theta}\,d\theta.
\]
Since $e^{i\Theta_\sharp}=i\sgn(\theta)e^{i\theta/2}$, writing $u=1/2-n$ gives
\[
 G_n=\frac{i}{2\pi}\left(\int_0^\pi e^{iu\theta}\,d\theta-
 \int_{-\pi}^0 e^{iu\theta}\,d\theta\right)
 =-\frac1{\pi u}=\frac2{\pi(2n-1)}.
\]
The CAR covariance \eqref{eq:covariance} yields
\[
 \omega_\sharp(i a_jb_{j+r})=-G_{-r},\qquad
 \omega_\sharp(i b_ja_{j+r+1})=G_{r+1}.
\]
Both are $2/(\pi(2r+1))$. Covariances within a sublattice vanish at distinct sites, giving the statement for even lags. All two-point functions, including the diagonal CAR term, are therefore invariant under Majorana translation. Quasifreeness then gives invariance of the entire state, not only of individual generator expectations.
\end{proof}

\begin{corollary}
For every $k$,
\[
 \psi_\sharp(e_k)=\frac12+\frac1\pi,\qquad
 \phi_\sharp(\sigma_z^{(j)})=
 \phi_\sharp(\sigma_x^{(j)}\sigma_x^{(j+1)})=\frac2\pi.
\]
\end{corollary}
\begin{proof}
Use $\psi_\sharp(e_k)=(1+\Omega(1))/2$ and \eqref{eq:TLgenerators}.
\end{proof}

For distinct indices $j_1,\ldots,j_{2s}$, Wick's rule is
\begin{equation}
 \omega_\sharp(m_{j_1}\cdots m_{j_{2s}})=\Pf K,\qquad
 K_{ab}=-i\Omega(j_b-j_a)\quad(a<b),\quad K_{ba}=-K_{ab}.
 \label{eq:Pfaffian}
\end{equation}
Odd monomials have zero expectation. To evaluate a TL word, expand each $e_k$, reduce repeated Majoranas using the CAR, and apply \eqref{eq:Pfaffian}. The factors of $i$ in this formula are essential. This determines the state on every local element of $B$.

\subsection{What the categorical invariant records}

The Hamiltonian and covariance above are those of the standard critical Ising chain. Its free Majorana continuum description has central charge $c=1/2$, see for example the explicit discussion of the continuum limit in \cite{SS24}. This conformal identification is background about this particular critical dynamics. We do not infer central charge from $H_\phi$, or prove its invariance under an arbitrary automorphism of the quasi-local algebra.

The equality $H_\phi=\one$ says that none of the three nontrivial toric code sectors is implementable in the even vacuum representation. It does not distinguish all critical theories with the same symmetry, determine their stress tensors, or classify their continuum limits. The transfer matrix and quantum chain states have the same algebra object by independent calculations in Section~\ref{sec:sectors}. Their exponentially local conjugacy is an additional lattice statement. The Temperley--Lieb--Jones chains of Part~\ref{part:TLJ} give critical states with trivial algebra object and different central charges (Section~\ref{sec:scope}).

The distinction between categorical and spectral gaplessness is also concrete here. For $H_1$, $\varepsilon_\sharp(\theta)=2|\sin(\theta/2)|$ has arbitrarily small positive values. The transfer matrix dispersion in \eqref{eq:dispersion} has the same property. One- and two-particle wave packets with momenta near zero give energies tending to zero in the full Fock representation and in the even vacuum sector, respectively. Thus these specified critical dynamics have no positive spectral gap. The calculation of the algebra object does not require selecting these Hamiltonians.

\part{Temperley--Lieb--Jones chains below index four}
\label{part:TLJ}

\section{The Temperley--Lieb--Jones chain}
\label{sec:TLJchain}

\begin{figure}[ht]
\centering
\begin{tikzpicture}[scale=0.8]
\begin{scope}
 \draw[very thick] (0,0) -- (0,3.2);
 \draw[very thick] (3,0) -- (3,3.2);
 \draw[very thick] (1,0) -- (1,0.6) to[out=90,in=90] (2,0.6) -- (2,0);
 \draw[very thick] (1,3.2) -- (1,2.6) to[out=-90,in=-90] (2,2.6) -- (2,3.2);
 \node[below] at (0,0) {\scriptsize $i-1$};
 \node[below] at (1,0) {\scriptsize $i$};
 \node[below] at (2,0) {\scriptsize $i+1$};
 \node[below] at (3,0) {\scriptsize $i+2$};
 \node at (1.5,1.25) {\scriptsize $R_X^*$};
 \node at (1.5,1.95) {\scriptsize $R_X$};
 \node at (1.5,-1.0) {(a) $E_i=\delta e_i$};
\end{scope}
\begin{scope}[xshift=6cm]
 \foreach \x in {0,1.6,3.2} {\draw[very thick] (\x,2.0) -- (\x,1.2) to[out=-90,in=-90] (\x+0.8,1.2) -- (\x+0.8,2.0);}
 \draw[decorate,decoration={brace,mirror,amplitude=3pt}] (-0.05,0.55) -- (0.85,0.55);
 \node[below] at (0.4,0.45) {\scriptsize site $1$};
 \node[below] at (3.6,0.45) {\scriptsize site $N$};
 \node at (2.0,1.6) {$\cdots$};
 \node at (2.0,-1.0) {(b) $\mathrm{cup}_N=\delta^{N/2}\xi_N$};
\end{scope}
\end{tikzpicture}
\caption{(a) The cup--cap $E_i$. It is the cap $R_X^*$ on the strands $i,i+1$ followed by the cup $R_X$, with identities on the other strands. Diagrams are read from bottom to top. (b) The cup vector in $V_N=\Hom(\one,X^{\otimes2N})$, one cup on each $Y$-site.}
\label{fig:TLgen}
\end{figure}
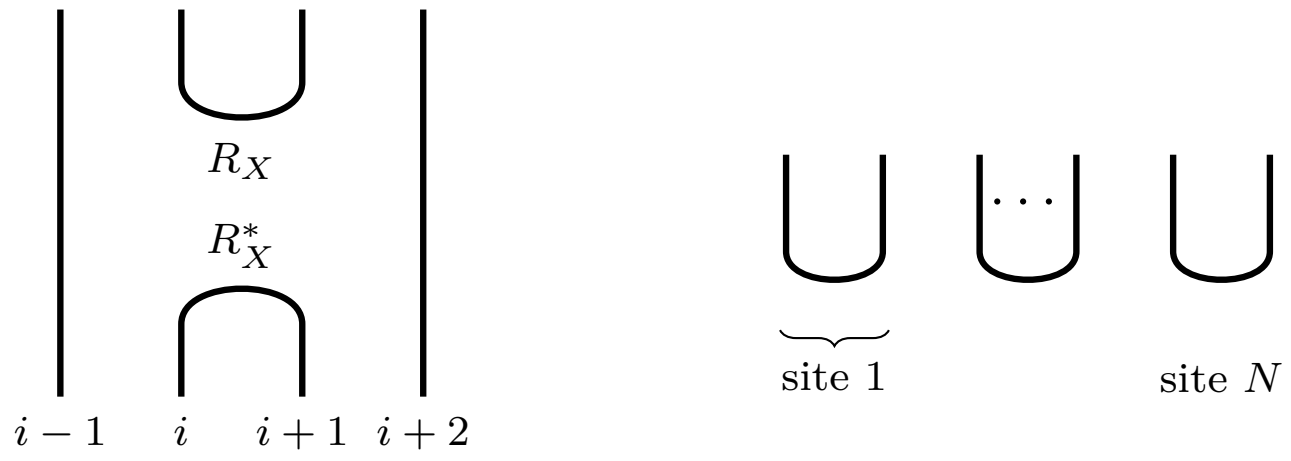

Fix an integer $k\ge2$ and put $\delta=2\cos(\pi/(k+2))$, so that $\delta^2\in[2,4)$. Let $\mathcal T_k$ be the unitary Temperley--Lieb--Jones category, in which the Jones--Wenzl idempotent $p_{k+1}$ is zero. Write $X$ for the fundamental object, $Y=X\otimes X$, and $\C_k=\langle Y\rangle$ for the even part. We use unnormalized spherical traces, $\mathrm{Tr}_A(\id_A)=d(A)$, a fixed unitary strict model, and standard solutions of the conjugate equations. Put
\[
 B_k(I)=\End_{\C_k}\bigl(Y^{\otimes|I|}\bigr),\qquad
 B_k=\overline{\textstyle\bigcup_IB_k(I)},
\]
with inclusions given by adding identity strands. Inner products are linear in the second variable, as in Section~\ref{sec:CAR}. The $X$-strands of the chain are numbered by $i\in\Z$, and the strands $2j-1,2j$ form the $Y$-site $j$. Let $R_X:\one\to X\otimes X$ be a standard solution of the conjugate equations for the self-dual object $X$, so that $R_X^*R_X=\delta$. We draw $R_X$ as a cup and $R_X^*$ as a cap. For $i\in\Z$ let
\[
 E_i=\id\otimes R_XR_X^*\otimes\id,\qquad e_i=\delta^{-1}E_i,
\]
where $R_XR_X^*$ acts on the strands $i$ and $i+1$ and the identities act on the other strands. So $E_i$ is the cup--cap of Figure~\ref{fig:TLgen}(a). For every interval $I$ of $Y$-sites containing the sites of the strands $i$ and $i+1$, this defines elements $E_i,e_i\in B_k(I)=\End_{\C_k}(X^{\otimes2|I|})$, compatible with the inclusions. The $e_i$ are the Temperley--Lieb generators. They satisfy $E_i^2=\delta E_i$ and
\[
 e_i=e_i^*=e_i^2,\qquad e_ie_{i\pm1}e_i=\delta^{-2}e_i,\qquad[e_i,e_j]=0\ (|i-j|\ge2).
\] Since $\End(X^{\otimes2n})$ is the quotient of the Temperley--Lieb algebra by the ideal generated by $p_{k+1}$, each $B_k(I)$ is generated by the $e_i$ with both strands in $I$. The relations, and the ideals generated by the Jones--Wenzl idempotents, are invariant under $e_i\mapsto e_{i+1}$, and the local algebras are finite dimensional with unique $C^*$-norms. Hence there is an automorphism $\varsigma$ of $B_k$ with
\[
 \varsigma(e_i)=e_{i+1},
\]
and $\varsigma^2$ is translation by one $Y$-site. Since $Y\supseteq\one$ tensor generates $\C_k$, some power of $Y$ contains every simple object of $\C_k$, and
\[
 \DHR(B_k)\simeq\mathcal Z(\C_k)
\]
\cite{Jon24}, see also \cite[Example~2.5 and Section~6]{JNPW}. This is the centre of the even part, not of the full category containing $X$. It is the quantum double of the Jones subfactor of type $A_{k+1}$, whose even part is $\C_k$, and it was computed by Evans and Kawahigashi through Ocneanu's asymptotic inclusion \cite[Section~4]{EK98b}. The categorical invariant of a pure state $\omega$ of $B_k$ is $H_\omega=\underline{\End}(L^2(B_k,\omega))\in\mathcal Z(\C_k)$, as in Section~\ref{sec:invariant}, and the terms local, topological and gapless have the same meaning.

Let $I_N=\{x\in\Z:-\lfloor N/2\rfloor<x\le\lceil N/2\rceil\}$ be the interval of the $N$ $Y$-sites around the origin, carrying the strands $2x-1,2x$ for $x\in I_N$, and let $S_N$ be the sum of the $2N-1$ generators $e_i$ with both strands in $I_N$. The finite volume Hamiltonian and its vacuum block are
\[
 h_N=-S_N,\qquad V_N=\Hom(\one,X^{\otimes2N}).
\]
Each local $a\in B_k(I)$ with $I\subseteq I_N$ acts on $V_N$ by composition, and $V_N$ has the inner product $\langle v,w\rangle=v^*w\in\End(\one)=\mathbb C$. The cup vector is
\[
 \mathrm{cup}_N=R_X^{\otimes N}\in V_N,
\]
the $N$ adjacent cups, one on each $Y$-site of $I_N$ (Figure~\ref{fig:TLgen}(b)). Its norm is $\delta^{N/2}$, and $\xi_N=\delta^{-N/2}\mathrm{cup}_N$ is the normalized cup vector.

\begin{remark}[The case $k=2$]
\label{rem:k2chain}
At $k=2$, $\delta=\sqrt2$ and $\C_2\simeq\Vect(\Z_2)$. The elements $\frac12(1+i\,m_im_{i+1})$ of \eqref{eq:TLgenerators} satisfy the relations above with $\delta^{-2}=\frac12$ (Section~\ref{sec:TL}), and with $E_i=\sqrt2e_i$ the Jones--Wenzl idempotent $p_3=1-\sqrt2(E_1+E_2)+(E_1E_2+E_2E_1)$ vanishes in the Clifford algebra. Since both local algebras have dimension $2^{2n-1}$, the assignment $e_i\mapsto\frac12(1+i\,m_im_{i+1})$ extends to an isomorphism of $B_2$ onto the even Pauli algebra $B$ of Part~\ref{part:Ising}, with $\varsigma=\kappa$ by \eqref{eq:KW}. Under it $\DHR(B_2)$ is the toric code of Part~\ref{part:Ising} and $-\sum_ie_i=H_1+\text{const}$. The cup vector $\xi_N$ is fixed by every $e_{2j-1}$, which singles out $V_N$ among the two irreducible representations of $B_{I_N}$, and corresponds to the vector with $\sigma_z^{(j)}=+1$ in the even subspace of the open chain on $I_N$. The state $\omega_2$ constructed below is a ground state of the even dynamics of $H_1$ (Proposition~\ref{prop:ground}), so $\omega_2=\psi_\sharp$ by \cite[Theorem~4(1)($\alpha$)]{ArakiMatsui85}, as in the proof of Proposition~\ref{prop:groundstate}.
\end{remark}

\section{Loop representation and the infinite volume state}
\label{sec:loops}

This section constructs the state $\omega_k$ and the probabilistic tools used in Sections~\ref{sec:purity}--\ref{sec:decay}. The finite volume ground state is obtained from the Perron--Frobenius theorem in the path representation (Section~\ref{sec:loops}.1). The time-ordered expansion of the exponential, in the form used by Aizenman and Nachtergaele, turns its imaginary time correlations into expectations for a loop model. This loop model is the image of the quantum random-cluster model with cluster weight $q=\delta^2$ at its self-dual point (Section~\ref{sec:loopsub}). The cups of $\xi_N$ give the free boundary condition, so the state is built from the free measures alone, and its limit is the free infinite volume measure $\phi^0$. The wired measures enter only as a comparison. Conditioning on the configuration outside a region induces on it a boundary condition that lies between free and wired. Primal crossings are least likely under the free boundary condition, and dual crossings under the wired one. For this model Duminil-Copin, Li and Manolescu proved the Russo--Seymour--Welsh crossing bounds and that the free and wired measures coincide, $\phi^0=\phi^1$ (Lemma~\ref{lem:critical-loops}). The crossing bounds hold at every scale and whatever the outside configuration. With positive association and the domain Markov property they give separating loops in annuli at every scale (Lemma~\ref{lem:separating}) and polynomial bounds on loop sizes (Lemma~\ref{lem:looptail}). These separating loops drive the thermodynamic limit, purity (Section~\ref{sec:purity}) and the decay of charged coefficients (Section~\ref{sec:decay}). The equality $\phi^0=\phi^1$ is used only for mixing. Any boundary condition induced by a far configuration lies between $\phi^0$ and $\phi^1$, and is forgotten because they coincide (Lemma~\ref{lem:looptail}, Remark~\ref{rem:mixing}). These yield the thermodynamic limit $\omega_k$, a translation invariant ground state (Proposition~\ref{prop:ground}).

\subsection{Positivity and the finite volume ground state}

We use the path representation on the graph $A_{k+1}$ \cite{Pasquier,GHJ,Jon00,EK98}, with the Perron--Frobenius weights
\[
 \mu_a=\sin\frac{a\pi}{k+2},\qquad\sum_{b\sim a}\mu_b=\delta\mu_a .
\]
The space $V_N$ has a basis of paths of length $2N$ on $A_{k+1}$ from the vertex $1$ to itself. The generator $e_i$ acts on the $i$th step. On a path passing through $a,b,a$ at positions $i-1,i,i+1$ it acts by the rank-one projection onto
\[
 \textstyle\sum_{b\sim a}\sqrt{\mu_b/(\delta\mu_a)}\,|b\rangle,
\]
and it is zero unless the two outer vertices agree. This is the embedding of the Temperley--Lieb planar algebra in the planar algebra of the graph $A_{k+1}$ \cite{Jon00}, see also \cite[Section~2.2]{GJSZ}. All matrix entries of the $e_i$ in the path basis are therefore nonnegative. Every path from $1$ to $1$ reduces to the alternating path $1,2,1,\dots,1$ by lowering peaks, and each lowering and its reverse has a positive entry in some $e_i$. So the matrix of $S_N$ on $V_N$ is nonnegative and irreducible, and the Perron--Frobenius theorem applies to it. Its largest eigenvalue is simple, with an eigenvector $\Omega_N$ of positive coordinates, which we normalize. Equivalently, the lowest eigenvalue
\[
 \epsilon_N=\min\mathrm{Spec}(h_N|_{V_N})
\]
of $h_N=-S_N$ is simple with eigenvector $\Omega_N$, and the rest of the spectrum of $h_N|_{V_N}$ lies in $[\epsilon_N+g_N,\infty)$ for some $g_N>0$. Consequently $e^{-T(h_N-\epsilon_N)}\to|\Omega_N\rangle\langle\Omega_N|$ as $T\to\infty$, so that
\begin{equation}
 e^{-T(h_N-\epsilon_N)}\psi_N\longrightarrow\langle\Omega_N,\psi_N\rangle\,\Omega_N\qquad(T\to\infty)
 \label{eq:PF}
\end{equation}
for every $\psi_N\in V_N$, with a nonzero limit whenever $\langle\Omega_N,\psi_N\rangle\ne0$, in particular for every nonzero $\psi_N$ with nonnegative coordinates. We use the cup vector $\xi_N$ of Section~\ref{sec:TLJchain}, which in the path basis is a positive multiple of the basis vector of the alternating path $1,2,1,\dots,1$, since a cup starting at the vertex $1$ passes through its only neighbour $2$. So $\langle\Omega_N,\xi_N\rangle>0$. It is chosen not for \eqref{eq:PF} but because it makes $\omega_{N,T}$ a loop expectation with free boundary conditions in time (Section~\ref{sec:loopsub}).
Write $\omega_N(a)=\langle\Omega_N,a\Omega_N\rangle$ for the finite volume ground state.

\paragraph{Schwinger functions.} Fix an imaginary time cutoff $T>0$. For local $a_1,\dots,a_r\in B_k$ supported in $I_N$ (Section~\ref{sec:TLJchain}) and times $-T\le t_1\le\dots\le t_r\le T$, the finite volume Schwinger functions are the imaginary time correlation functions
\begin{equation}
 \mathcal S_{N,T}(a_1,t_1;\dots;a_r,t_r)=\frac{\langle\xi_N,e^{-(T-t_r)h_N}a_re^{-(t_r-t_{r-1})h_N}\cdots e^{-(t_2-t_1)h_N}a_1e^{-(T+t_1)h_N}\xi_N\rangle}{\langle\xi_N,e^{-2Th_N}\xi_N\rangle}.
 \label{eq:schwinger}
\end{equation}
In the loop diagrams of Section~\ref{sec:loopsub}, as in Figure~\ref{fig:loopcluster}(a), imaginary time is the vertical coordinate $t\in[-T,T]$. The vector $\xi_N$ sits at the bottom edge $t=-T$ and its adjoint at the top edge $t=T$, the factor $e^{-(t_{j+1}-t_j)h_N}$ acts on the time interval between consecutive insertions, and $a_j$ is inserted at height $t_j$. A single insertion at time $0$ gives the state
\begin{equation}
 \omega_{N,T}(a)=\mathcal S_{N,T}(a,0)=\frac{\langle\xi_N,e^{-Th_N}ae^{-Th_N}\xi_N\rangle}{\langle\xi_N,e^{-2Th_N}\xi_N\rangle}.
 \label{eq:finitestate}
\end{equation}
By \eqref{eq:PF}, as $T\to\infty$ the Schwinger functions converge to those of the finite volume ground state,
\[
 \mathcal S_{N,T}(a_1,t_1;\dots;a_r,t_r)\longrightarrow\langle\Omega_N,a_re^{-(t_r-t_{r-1})(h_N-\epsilon_N)}\cdots e^{-(t_2-t_1)(h_N-\epsilon_N)}a_1\Omega_N\rangle,
\]
and in particular $\omega_{N,T}\to\omega_N$ on local elements. A word of length $r$ is a product $W=e_{i_1}\cdots e_{i_r}$ of generators. Each $B_k(I)=\End_{\C_k}(X^{\otimes2|I|})$ is spanned by words \cite{GHJ}, so it suffices to treat insertions that are words.

\subsection{The loop law}
\label{sec:loopsub}

The cup--caps $E_i=\delta e_i$ of Section~\ref{sec:TLJchain} satisfy $E_i^2=\delta E_i$, and $S_N=\delta^{-1}\sum_iE_i$, the sum over the $2N-1$ bonds of $I_N$, the bond $i$ joining the strands $i$ and $i+1$. The formula \eqref{eq:poisson} below writes the power series of $e^{tS_N}$ as a sum of time-ordered products of cup--caps. Every term is a Temperley--Lieb diagram with a positive weight, so that the weights, once normalised, form a probability measure on diagrams. Aizenman and Nachtergaele write the same expansion as an expectation over Poisson processes \cite[Section~2.2, (2.11)--(2.12)]{AN94}.

\paragraph{Configurations.} An event is a pair $(i,\tau)$ of a bond $i$ and a time $\tau$. A configuration $\Gamma$ on a time interval $[s,t]$ is a finite set of events with times in $[s,t]$, and $n(\Gamma)$ is its number of events. The configurations with $n$ events and bond labels $i_1,\dots,i_n$, listed in increasing time, are parametrised by the simplex $s<\tau_1<\dots<\tau_n<t$. We write $d\Gamma$ for the measure on configurations that is Lebesgue measure on each of these simplices, summed over $n$ and over the bond labels.

\paragraph{The expansion.} For $t\ge0$,
\begin{equation}
 e^{tS_N}=\int\delta^{-n(\Gamma)}\,E(\Gamma)\,d\Gamma
 =\sum_{n\ge0}\ \sum_{i_1,\dots,i_n}\ \int_{0<\tau_1<\dots<\tau_n<t}\delta^{-n}\,E_{i_n}\cdots E_{i_1}\,d\tau_1\cdots d\tau_n ,
 \label{eq:poisson}
\end{equation}
where the integral is over the configurations on $[0,t]$ and $E(\Gamma)$ is the product of the $E_i$ over the events of $\Gamma$, ordered with later times to the left. Indeed the simplex has volume $t^n/n!$, so the term with $n$ events is $\frac{t^n}{n!}\bigl(\delta^{-1}\sum_iE_i\bigr)^n=\frac{t^n}{n!}S_N^n$. Every event carries the same weight $\delta^{-1}$, on odd and on even bonds alike. The measure $e^{-(2N-1)t/\delta}\delta^{-n(\Gamma)}d\Gamma$ is the law of independent Poisson processes of rate $1/\delta$ on the bonds, and in that form \eqref{eq:poisson} is \cite[(2.10)--(2.12)]{AN94}, with $h_b=E_i$, $J_b=1/\delta$ and $\beta=t$. The same parametrisation of the Dyson series by time-ordered labelled events, with counting measure on the labels and Lebesgue measure on the ordered times, goes back to Davies \cite{Davies69} and to Evans and Lewis \cite[(2.11)--(2.16)]{EvansLewis76}, for semigroups of completely positive maps. The times record the order of the factors and place each cup--cap at a height in a space-time picture. A configuration $\Gamma$ on $[-T,T]$ is drawn as a diagram with $2N$ vertical strands and a cup--cap between the strands $i$ and $i+1$ at height $\tau$ for each event $(i,\tau)$ of $\Gamma$. Closing it at the top and bottom with the $N$ cups of $\xi_N$ gives a planar diagram of closed loops (Figure~\ref{fig:loopcluster}(a)). We use the letter $\Gamma$ also for this diagram and its loops, and write $\ell(\Gamma)$ for the number of loops. This is the loop decomposition of \cite[Section~2.1 and Figure~1]{AN94}, with the boundary vectors $\xi_N$ in place of the trace, as in their treatment of open paths \cite[Section~2.4.b]{AN94}.

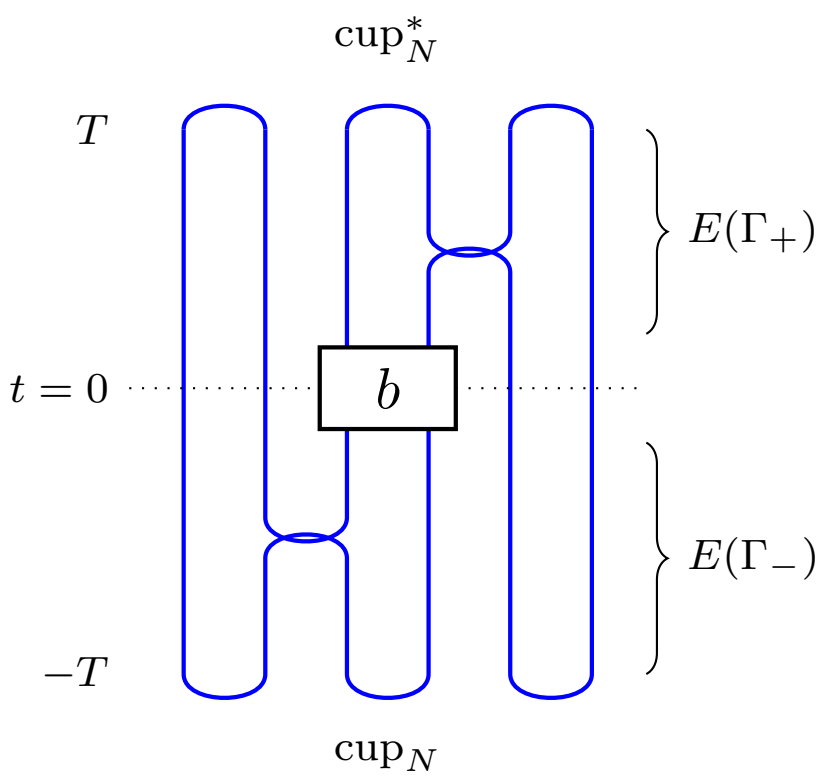
\begin{figure}[ht]
\centering
\begin{tikzpicture}[scale=0.9]
 \draw[thick,blue] (0.6,0) -- (0.6,4);
 \draw[thick,blue] (3.6,0) -- (3.6,4);
 \draw[thick,blue] (1.2,0) -- (1.2,0.85) to[out=90,in=90] (1.8,0.85) -- (1.8,0);
 \draw[thick,blue] (1.2,4) -- (1.2,1.15) to[out=-90,in=-90] (1.8,1.15) -- (1.8,1.8);
 \draw[thick,blue] (1.8,2.4) -- (1.8,4);
 \draw[thick,blue] (2.4,0) -- (2.4,1.8);
 \draw[thick,blue] (2.4,2.4) -- (2.4,2.95) to[out=90,in=90] (3.0,2.95) -- (3.0,0);
 \draw[thick,blue] (2.4,4) -- (2.4,3.25) to[out=-90,in=-90] (3.0,3.25) -- (3.0,4);
 \foreach \x/\y in {0.6/1.2,1.8/2.4,3.0/3.6} {\draw[thick,blue] (\x,0) to[out=-90,in=-90] (\y,0); \draw[thick,blue] (\x,4) to[out=90,in=90] (\y,4);}
 \draw[thick,fill=white] (1.6,1.8) rectangle (2.6,2.4);
 \node at (2.1,2.1) {$b$};
 \draw[dotted] (0.2,2.1) -- (1.6,2.1) (2.6,2.1) -- (4.0,2.1);
 \node[left] at (0.2,2.1) {\scriptsize $t=0$};
 \node[left] at (0.2,0) {\scriptsize $-T$};
 \node[left] at (0.2,4) {\scriptsize $T$};
 \draw[decorate,decoration={brace,amplitude=4pt}] (4.0,4.0) -- (4.0,2.5);
 \node[right] at (4.15,3.25) {\scriptsize $E(\Gamma_+)$};
 \draw[decorate,decoration={brace,amplitude=4pt}] (4.0,1.7) -- (4.0,0.0);
 \node[right] at (4.15,0.85) {\scriptsize $E(\Gamma_-)$};
 \node[below] at (2.1,-0.35) {\scriptsize $\mathrm{cup}_N$};
 \node[above] at (2.1,4.35) {\scriptsize $\mathrm{cup}_N^*$};
\end{tikzpicture}
\caption{The evaluation \eqref{eq:evaluation} for $N=3$. It shows the cups of $\mathrm{cup}_N$ at the bottom, the cup--caps of the events of $\Gamma$ at negative times, the element $b$ at time $0$, the cup--caps at positive times, and the caps of $\mathrm{cup}_N^*$ at the top. For $b=1$ the closed diagram has $\ell(\Gamma)$ loops and value $\delta^{\ell(\Gamma)}$.}
\label{fig:evaluation}
\end{figure}

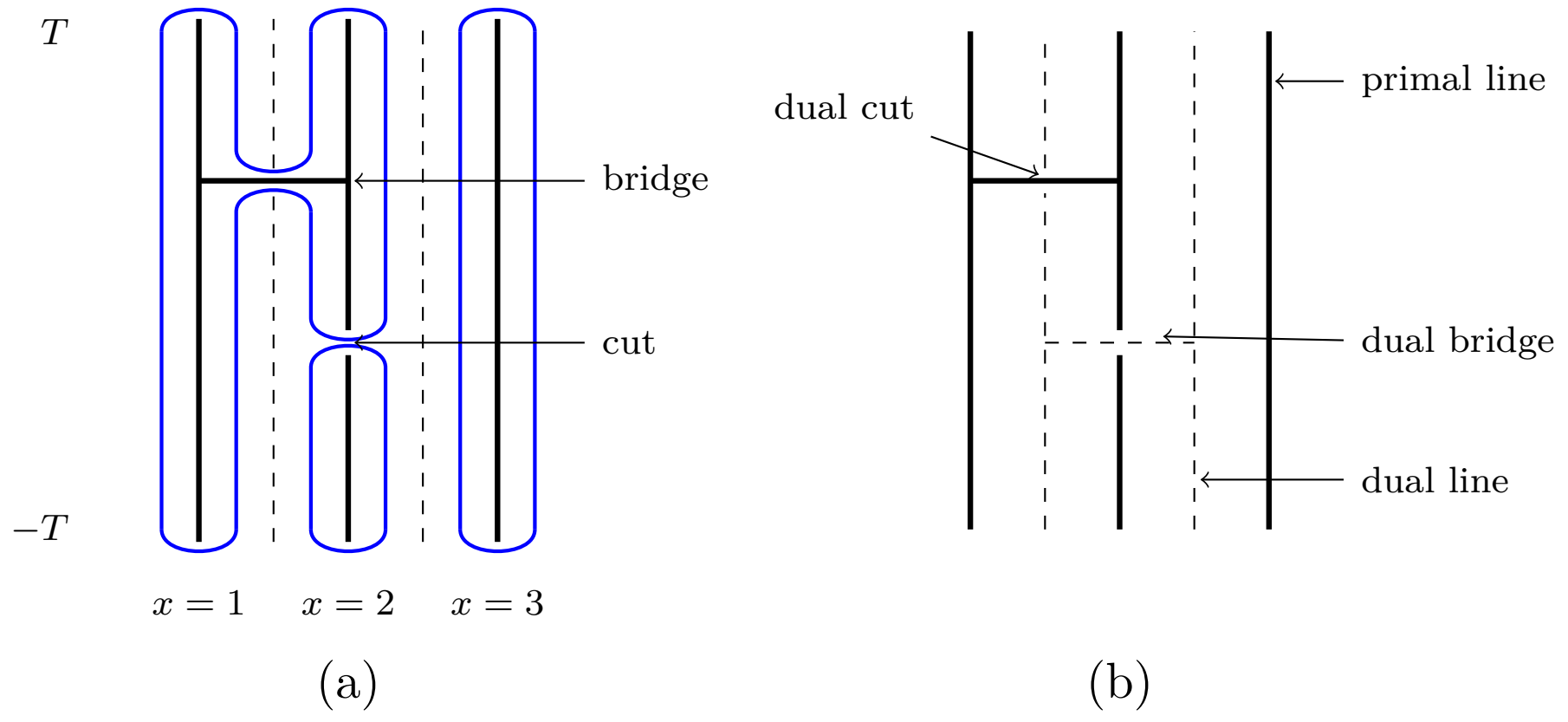
\begin{figure}[ht]
\centering
\begin{tikzpicture}[scale=0.95]
\begin{scope}
 \foreach \x in {0.6,3.0,3.6} {\draw[thick,blue] (\x,0) -- (\x,4);}
 \draw[thick,blue] (1.8,0) -- (1.8,1.3) to[out=90,in=90] (2.4,1.3) -- (2.4,0);
 \draw[thick,blue] (1.8,2.55) -- (1.8,1.7) to[out=-90,in=-90] (2.4,1.7) -- (2.4,4);
 \draw[thick,blue] (1.2,0) -- (1.2,2.55) to[out=90,in=90] (1.8,2.55);
 \draw[thick,blue] (1.2,4) -- (1.2,3.05) to[out=-90,in=-90] (1.8,3.05) -- (1.8,4);
 \foreach \x/\y in {0.6/1.2,1.8/2.4,3.0/3.6} {\draw[thick,blue] (\x,0) to[out=-90,in=-90] (\y,0); \draw[thick,blue] (\x,4) to[out=90,in=90] (\y,4);}
 \draw[very thick] (0.9,-0.1) -- (0.9,4.1);
 \draw[very thick] (2.1,-0.1) -- (2.1,1.4) (2.1,1.6) -- (2.1,4.1);
 \draw[very thick] (3.3,-0.1) -- (3.3,4.1);
 \draw[very thick] (0.9,2.8) -- (2.1,2.8);
 \draw[dashed] (1.5,-0.1) -- (1.5,2.7) (1.5,2.9) -- (1.5,4.1);
 \draw[dashed] (2.7,-0.1) -- (2.7,4.1);
 \node[right] (c) at (4.0,1.5) {\scriptsize cut};
 \draw[->,thin] (c) -- (2.15,1.5);
 \node[right] (b) at (4.0,2.8) {\scriptsize bridge};
 \draw[->,thin] (b) -- (2.15,2.8);
 \node[below] at (0.9,-0.35) {\scriptsize $x=1$};
 \node[below] at (2.1,-0.35) {\scriptsize $x=2$};
 \node[below] at (3.3,-0.35) {\scriptsize $x=3$};
 \node[left] at (0,0) {\scriptsize $-T$};
 \node[left] at (0,4) {\scriptsize $T$};
 \node at (2.1,-1.25) {(a)};
\end{scope}
\begin{scope}[xshift=6.2cm]
 \draw[very thick] (0.9,0) -- (0.9,4);
 \draw[very thick] (2.1,0) -- (2.1,1.4) (2.1,1.6) -- (2.1,4);
 \draw[very thick] (3.3,0) -- (3.3,4);
 \draw[very thick] (0.9,2.8) -- (2.1,2.8);
 \draw[dashed] (1.5,0) -- (1.5,2.7) (1.5,2.9) -- (1.5,4);
 \draw[dashed] (2.7,0) -- (2.7,4);
 \draw[dashed] (1.5,1.5) -- (2.7,1.5);
 \node[right] (pl) at (3.9,3.6) {\scriptsize primal line};
 \draw[->,thin] (pl) -- (3.35,3.6);
 \node[right] (dl) at (3.9,0.4) {\scriptsize dual line};
 \draw[->,thin] (dl) -- (2.75,0.4);
 \node[right] (db) at (3.9,1.5) {\scriptsize dual bridge};
 \draw[->,thin] (db) -- (2.45,1.55);
 \node[left] (dc) at (0.6,3.4) {\scriptsize dual cut};
 \draw[->,thin] (dc) -- (1.45,2.85);
 \node at (2.1,-1.25) {(b)};
\end{scope}
\end{tikzpicture}
\caption{(a) A configuration $\Gamma$ on three $Y$-sites, with an event on the odd bond inside the site $2$ (strands $3,4$) and an event on the even bond between the sites $1$ and $2$ (strands $2,3$), closed by the cups of $\xi_N$ at the top and bottom. Primal lines are solid, dual lines dashed. The first event cuts the primal line $x=2$, and the second is a bridge between the primal lines $x=1$ and $x=2$. (b) The primal configuration $\omega$ (solid) and the dual configuration $\omega^*$ (dashed). The bridge cuts the dual line $x=\frac32$, and the cut carries a dual bridge. The loops of (a) run between primal and dual clusters. They are the boundary curves of thin neighbourhoods of the three primal clusters, here $N=3$, $c=b=1$, $c_0=3$, $c_0^*=1$ and $\ell=c_0+c_0^*-1=2c_0-(N+c)+b=3$.}
\label{fig:loopcluster}
\end{figure}

\paragraph{Evaluations.} Let $\mathrm{cup}_N=\delta^{N/2}\xi_N$ be the unnormalized cup vector of Section~\ref{sec:TLJchain}, and let $E(\Gamma_+)$ and $E(\Gamma_-)$ be the time-ordered products of the $E_i$ over the events $(i,\tau)\in\Gamma$ with $\tau>0$ and with $\tau<0$. For $b\in B_k$ supported in $I_N$ the evaluation of $\Gamma$ with $b$ inserted at time $0$ is
\begin{equation}
 \mathrm{ev}_\Gamma(b)=\langle\mathrm{cup}_N,E(\Gamma_+)\,b\,E(\Gamma_-)\,\mathrm{cup}_N\rangle .
 \label{eq:evaluation}
\end{equation}
It is the value of the closed diagram obtained by inserting the diagram of $b$ into that of $\Gamma$ at time $0$ (Figure~\ref{fig:evaluation}).

 The uncharged evaluation is $\mathrm{ev}_\Gamma(1)=\delta^{\ell(\Gamma)}$. For a word $W=e_{i_1}\cdots e_{i_r}$, $\mathrm{ev}_\Gamma(W)=\delta^{\ell(\Gamma;W)-r}$, where $\ell(\Gamma;W)$ is the number of loops after inserting the $r$ cup--caps of $W$ at time $0$. By \eqref{eq:poisson} for the two time intervals $[-T,0]$ and $[0,T]$, and since the normalization of the cups is common to the numerator and the denominator of \eqref{eq:finitestate},
\[
 \omega_{N,T}(b)=\frac{\int\delta^{-n(\Gamma)}\,\mathrm{ev}_\Gamma(b)\,d\Gamma}{\int\delta^{-n(\Gamma)}\,\delta^{\ell(\Gamma)}\,d\Gamma},
\]
with both integrals over the configurations on $[-T,T]$. The finite volume loop law is the probability measure
\begin{equation}
 d\nu_{N,T}(\Gamma)=\frac{\delta^{-n(\Gamma)}\,\delta^{\ell(\Gamma)}}{Z_{N,T}}\,d\Gamma,\qquad Z_{N,T}=\int\delta^{-n(\Gamma)}\,\delta^{\ell(\Gamma)}\,d\Gamma,
 \label{eq:looplaw}
\end{equation}
and with it
\begin{equation}
 \omega_{N,T}(b)=\mathbb E_{\nu_{N,T}}\Bigl[\frac{\mathrm{ev}_\Gamma(b)}{\delta^{\ell(\Gamma)}}\Bigr],\qquad
 \omega_{N,T}(W)=\mathbb E_{\nu_{N,T}}[f_W],\qquad
 f_W(\Gamma)=\delta^{\ell(\Gamma;W)-\ell(\Gamma)-r}\in[\delta^{-2r},1].
 \label{eq:estimator}
\end{equation}
The range of $f_W$ holds because each inserted cup--cap either joins two loops into one or splits one loop into two, and so changes the number of loops by $\pm1$. With words $W_j$ of lengths $r_j$ inserted at the times $t_j$, the same computation gives
\[
 \mathcal S_{N,T}(W_1,t_1;\dots;W_s,t_s)=\mathbb E_{\nu_{N,T}}\bigl[\delta^{\ell(\Gamma;W_1,\dots,W_s)-\ell(\Gamma)-\sum_jr_j}\bigr],
\]
where $\ell(\Gamma;W_1,\dots,W_s)$ counts the loops after all the insertions.

\paragraph{The random-cluster model.} Configurations live on the space-time strip $I_N\times[-T,T]$, which we regard as a subset of $\Z\times\mathbb R$, the site $x\in I_N$ giving the first coordinate and the time $t\in[-T,T]$ the second. The primal line of the $Y$-site $x$ is the vertical line $\{x\}\times\mathbb R$ midway between its strands $2x-1$ and $2x$, so it passes through the cups of $\xi_N$. The dual line $\{x+\frac12\}\times\mathbb R$ lies midway between the sites $x$ and $x+1$ (Figure~\ref{fig:loopcluster}). Space is measured in $Y$-sites, and time is the imaginary time parameter $t$ of $e^{-th_N}$, where $h_N=-S_N$ is the finite volume Hamiltonian of Section~\ref{sec:TLJchain}. A bond $i$ is odd if $i=2x-1$, so that it joins the two strands $2x-1,2x$ of the site $x$, and even if $i=2x$, so that it joins the strands $2x$ and $2x+1$ of the neighbouring sites $x$ and $x+1$. An event of $\Gamma$ at time $t$ on the odd bond inside the site $x$ is a cut of the primal line $x$ at $t$. An event on the even bond between the sites $x$ and $x+1$ is a bridge, a horizontal segment joining the primal lines $x$ and $x+1$ at height $t$. The primal configuration $\omega$ is the union of the primal lines with the cuts removed and the bridges added. Its connected components are the primal clusters. Dually, a bridge cuts the dual line $x+\frac12$, and a cut of the line $x$ at $t$ carries a dual bridge joining the dual lines $x\pm\frac12$. This gives the dual configuration $\omega^*$ and the dual clusters. Primal and dual paths never cross. By checkerboard shading \cite{AN94}, the loops of $\Gamma$ are the boundaries of thin neighbourhoods of the primal clusters, and each loop separates a primal cluster from a dual one.

We now count the loops. For a configuration with $c$ cuts and $b$ bridges, regard the primal configuration as a planar graph whose vertices are its $N+c$ uncut vertical segments and whose edges are its $b$ bridges. Let $c_0$ be the number of its connected components, the primal clusters. Each loop of $\Gamma$ is exactly one boundary curve of a thin neighbourhood of exactly one primal cluster, so $\ell=\ell(\Gamma)$ is the total number of these boundary curves. A cluster with $v$ segments and $e$ bridges is a connected planar graph with $v$ vertices and $e$ edges, so by Euler's formula it has $e-v+1$ bounded faces. Its thin neighbourhood is a disc with one hole for each bounded face, and has $1+(e-v+1)=2-v+e$ boundary curves. Summing over the $c_0$ clusters gives $\ell=2c_0-(N+c)+b$. Let $c_0^*$ be the number of dual clusters, where the dual clusters meeting the boundary of $\Delta_{N,T}=I_N\times[-T,T]$ count as one. Each bounded face of a primal cluster contains exactly one of the other dual clusters, so $c_0^*-1=c_0-(N+c)+b$ and
\[
 \ell=c_0+c_0^*-1 .
\]
Each primal cluster has one outer boundary loop, and each dual cluster other than the outer one has one further loop. With $n(\Gamma)=c+b$ and $q=\delta^2$ the weight of the loop law \eqref{eq:looplaw} is
\[
 \delta^{-n(\Gamma)}\delta^{\ell}=\delta^{-1}\,(\delta^{-1})^{c}\,(\delta^{-1})^{b}\,\sqrt q^{\,c_0}\,\sqrt q^{\,c_0^*}.
\]
Cuts and bridges carry the same weight $\delta^{-1}$, and the loop weight is $\sqrt q$ per primal cluster and $\sqrt q$ per dual cluster. The shift by one strand exchanges primal and dual lines, odd and even bonds, and cuts and bridges, and away from the boundary it preserves this weight. This is the self-duality of the model. The results of \cite{DCLM} are stated for a weight on the primal clusters alone. Eliminating $c_0^*$ gives
\[
 \delta^{-n(\Gamma)}\delta^{\ell}=\delta^{-N}\,(\delta^{-2})^{c}\,1^{b}\,q^{c_0},
\]
so that all the weight is on the primal clusters, with cuts of weight $\delta^{-2}$ and bridges of weight $1$. Eliminating $c_0$ instead puts the weight $q^{c_0^*}$ on the dual clusters and exchanges the weights of cuts and bridges. So $\nu_{N,T}$ is the image under $\omega\mapsto\Gamma$ of the random-cluster measure with cluster weight $q$, cut intensity $\lambda_{\mathrm c}=\delta^{-2}$ and bridge intensity $\lambda_{\mathrm b}=1$, on the rectangle $\Delta_{N,T}=I_N\times[-T,T]$ with free boundary conditions, in the sense now recalled. The ratio $\lambda_{\mathrm b}/\lambda_{\mathrm c}=q$ is the trace of the symmetric weight above after the dual count has been eliminated. It is the self-dual point of the model \cite[Section~5.1]{DCLM}. No Potts realization is needed, so nonintegral $q$ is allowed.

\paragraph{Finite volumes, boundary conditions and infinite volume measures.} A finite volume is a bounded region $\Delta\subset\Z\times\mathbb R$ that is a finite union of rectangles $([a,b]\cap\Z)\times[c,d]$. A configuration $\omega$ in $\Delta$ is given by cuts on the primal lines in $\Delta$ and bridges with both ends in $\Delta$, and is identified, as above, with the union of the uncut segments and the bridges. A boundary condition $\zeta$ is an equivalence relation on the points of the primal lines on the boundary of $\Delta$, recording which of them are connected outside $\Delta$. The random-cluster measure on $\Delta$ with boundary condition $\zeta$ is
\[
 d\phi^\zeta_\Delta(\omega)=\frac{q^{k^\zeta(\omega)}\,\lambda_{\mathrm c}^{c(\omega)}\,\lambda_{\mathrm b}^{b(\omega)}}{Z^\zeta_\Delta}\,d\omega,
\]
where $d\omega$ is the following measure on configurations. A configuration with $c$ cuts and $b$ bridges is determined by the times of its cuts on each primal line and of its bridges in each gap between adjacent lines. Listing these times in increasing order on each line and in each gap identifies the configurations with given numbers of cuts and bridges on each line and in each gap with a product of simplices in $\mathbb R^{c+b}$. The measure $d\omega$ is Lebesgue measure on each of these products, summed over the numbers of cuts and bridges, as for $d\Gamma$ above. It is a measure on the positions of the finitely many cuts and bridges, not on the lines themselves. Further, $c(\omega)$ and $b(\omega)$ are the numbers of cuts and bridges, $k^\zeta(\omega)$ is the number of clusters of $\omega$ when boundary points equivalent under $\zeta$ are identified (the cluster count), and $Z^\zeta_\Delta$ normalizes. Up to the normalization, $\lambda_{\mathrm c}^{c}\lambda_{\mathrm b}^{b}\,d\omega$ is the law of independent Poisson processes of cuts of intensity $\lambda_{\mathrm c}$ on each primal line and of bridges of intensity $\lambda_{\mathrm b}$ between adjacent lines, which is the form used in \cite{AN94}, \cite[Section~1.3]{DCLM}, \cite[Chapter~9]{Gri10}. The free boundary condition $\zeta=0$ identifies nothing, so that nothing is connected outside $\Delta$. The wired boundary condition $\zeta=1$ identifies all boundary points, as if the whole boundary were one connected wire, so that all clusters meeting the boundary count as a single cluster (Figure~\ref{fig:bc}). For integer $q$ the random-cluster model is the Fortuin--Kasteleyn representation of the $q$-state Potts model. The free and wired boundary conditions correspond to free and fixed boundary conditions for the Potts spins. In the wired case the boundary spins are all set equal to one fixed value, and by the symmetry of the Potts model the choice of this value does not matter \cite[Section~3]{AN94}. Here $q=\delta^2$ is in general not an integer, and only the random-cluster formulation is used.

\begin{figure}[ht]
\centering
\begin{tikzpicture}[scale=0.8,bdA/.style={thin},bdB/.style={line width=2.2pt}]
\foreach \s/\lab/\w in {0/{(a) free, $\zeta=0$}/bdA, 5.5/{(b) wired, $\zeta=1$}/bdB} {
\begin{scope}[xshift=\s cm]
 \draw[\w] (0,0) rectangle (4,3);
 \draw[very thick] (0.5,0) -- (0.5,1.2) (0.5,1.5) -- (0.5,3);
 \draw[very thick] (1.5,0) -- (1.5,0.8) (1.5,1.1) -- (1.5,2.2) (1.5,2.5) -- (1.5,3);
 \draw[very thick] (2.5,0) -- (2.5,1.7) (2.5,2.0) -- (2.5,3);
 \draw[very thick] (3.5,0) -- (3.5,0.6) (3.5,0.9) -- (3.5,3);
 \draw[very thick] (0.5,1.9) -- (1.5,1.9) (2.5,0.4) -- (3.5,0.4) (1.5,1.3) -- (2.5,1.3);
 \node at (2,-0.6) {\scriptsize \lab};
\end{scope}}
\end{tikzpicture}
\caption{A configuration in a rectangle $\Delta$ under the free and the wired boundary conditions. In (a) the clusters meeting the boundary are counted separately in $k^0(\omega)$. In (b) all boundary points are identified, drawn as a thick boundary, so all clusters meeting the boundary form one cluster in $k^1(\omega)$.}
\label{fig:bc}
\end{figure}
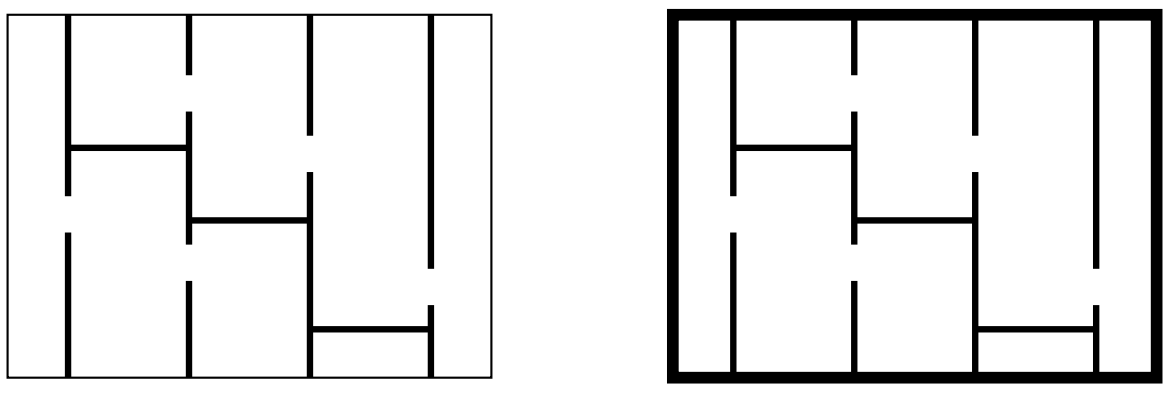 Configurations are ordered by inclusion as subsets of the plane. We write $\omega\le\omega'$ if $\omega\subseteq\omega'$. Removing a cut fills a gap in a primal line and so enlarges the subset, and adding a bridge enlarges it too. Hence $\omega\le\omega'$ exactly when every cut of $\omega'$ is a cut of $\omega$ and every bridge of $\omega$ is a bridge of $\omega'$, so a larger configuration has fewer cuts and more bridges. This is the order of \cite[Section~3 and (II.1)]{AN94}, whose $A$-bonds and $B$-bonds are our cuts and bridges. An event is increasing if it is preserved on passing to a larger configuration, and decreasing if its complement is increasing. Primal connection events are increasing and dual connection events are decreasing. For $q\ge1$ the model has the following three properties \cite[Chapter~9]{Gri10}, the continuum versions of those in \cite{FKG71} and of the classical comparison and Markov properties.
\begin{enumerate}[leftmargin=*,label=(\roman*)]
\item Positive association, the Fortuin--Kasteleyn--Ginibre (FKG) inequality \cite{FKG71}, \cite[Appendix~II]{AN94}. Increasing events are positively correlated under every $\phi^\zeta_\Delta$, and so are decreasing events.
\item Comparison of boundary conditions. We have $\phi^0_\Delta\le\phi^\zeta_\Delta\le\phi^1_\Delta$ in the sense of stochastic domination, where $\mu\le\mu'$ means that $\mu(E)\le\mu'(E)$ for every increasing event $E$, equivalently $\int f\,d\mu\le\int f\,d\mu'$ for every bounded increasing function $f$ \cite[(II.3), (II.5)]{AN94}. Thus increasing events are least likely under free and most likely under wired boundary conditions, and decreasing events the other way round.
\item The domain Markov property. For finite volumes $\Delta'\subseteq\Delta$, conditionally on the configuration in $\Delta\setminus\Delta'$, the configuration in $\Delta'$ has law $\phi^\xi_{\Delta'}$, where $\xi$ is the boundary condition induced by $\zeta$ and the connections in $\Delta\setminus\Delta'$.
\end{enumerate}
If $\Delta\subseteq\Delta'$, then by (iii) the restriction of $\phi^0_{\Delta'}$ to the configurations in $\Delta$ is an average of measures $\phi^\xi_\Delta$, each of which dominates $\phi^0_\Delta$ by (ii). So $\phi^0_{\Delta'}(E)\ge\phi^0_\Delta(E)$ for every increasing event $E$ depending on the configuration in $\Delta$, that is, $\phi^0_\Delta$ increases stochastically as $\Delta$ increases. In the same way $\phi^1_\Delta$ decreases stochastically. For an increasing event $E$ depending on the configuration in a bounded region, $\phi^0_\Delta(E)$ is therefore increasing in $\Delta$ and bounded by $1$, and it converges as $\Delta\uparrow\Z\times\mathbb R$. These limits determine a probability measure $\phi^0$, the free infinite volume measure, and the decreasing limits determine the wired infinite volume measure $\phi^1$ \cite[Chapter~9]{Gri10}, see also \cite[proof of Theorem~4.2]{AN94}. For a region $U$, let $\mathcal F_U$ be the $\sigma$-algebra on the space of configurations generated by the numbers of cuts and of bridges in the Borel subsets of $U$. An event or a function is $\mathcal F_U$-measurable exactly when it depends only on the cuts and bridges lying in $U$.

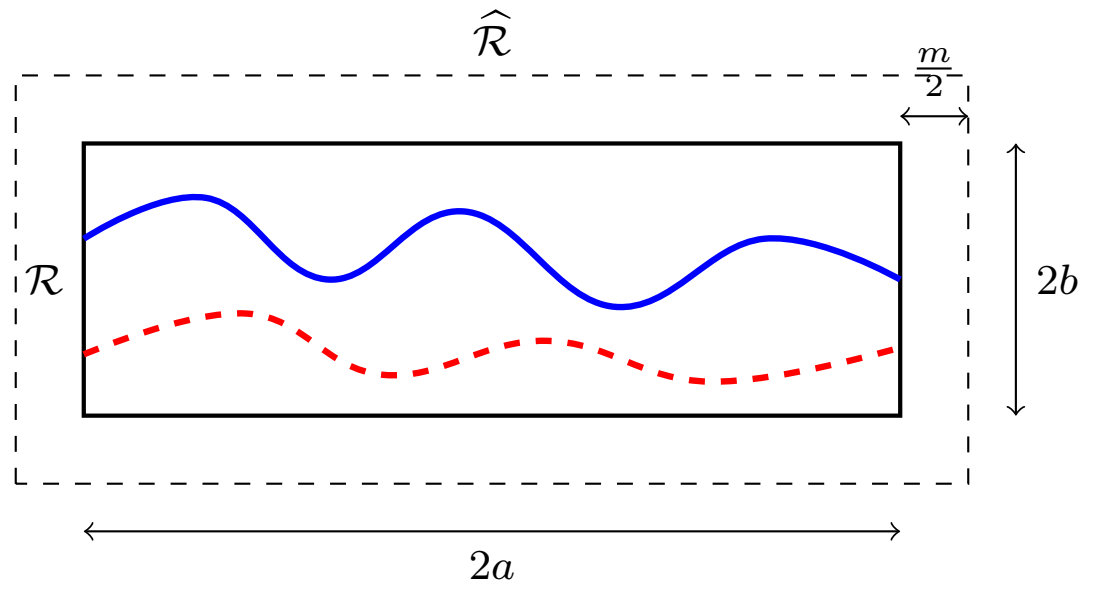
\begin{figure}[ht]
\centering
\begin{tikzpicture}[scale=0.9]
 \draw[dashed] (-0.5,-0.5) rectangle (6.5,2.5);
 \draw[thick] (0,0) rectangle (6,2);
 \draw[very thick,blue] plot[smooth,tension=0.7] coordinates {(0,1.3) (0.9,1.6) (1.8,1.0) (2.8,1.5) (3.9,0.8) (5.0,1.3) (6,1.0)};
 \draw[very thick,red,dashed] plot[smooth,tension=0.7] coordinates {(0,0.45) (1.2,0.75) (2.2,0.3) (3.4,0.55) (4.6,0.25) (6,0.5)};
 \node[above] at (3,2.5) {\scriptsize $\widehat{\mathcal R}$};
 \node[left] at (0,1.0) {\scriptsize $\mathcal R$};
 \draw[<->] (0,-0.85) -- (6,-0.85); \node[below] at (3,-0.85) {\scriptsize $2a$};
 \draw[<->] (6.85,0) -- (6.85,2); \node[right] at (6.85,1) {\scriptsize $2b$};
 \draw[<->] (6,2.2) -- (6.5,2.2); \node[above] at (6.25,2.2) {\scriptsize $\frac m2$};
\end{tikzpicture}
\caption{A rectangle $\mathcal R$ with $a>b$, so $m=b$ and the aspect ratio is $a/b$, inside its enlargement $\widehat{\mathcal R}$ (dashed). A primal crossing (solid) and a dual crossing (dashed) join the two shorter sides. Here the longer side lies in the space direction.}
\label{fig:crossing}
\end{figure}

For $z\in\Z\times\mathbb R$ and $r>0$ let $Q_r(z)$ be the space-time square of radius $r$ centred at $z$, that is the translate by $z$ of $([-r,r]\cap\Z)\times[-r,r]$. For a rectangle $\mathcal R$ with sides $2a$ in space and $2b$ in time, put $m=\min(a,b)$. The enlarged rectangle $\widehat{\mathcal R}$ has the same centre and sides $2a+m$ and $2b+m$, so it surrounds $\mathcal R$ by a margin equal to a quarter of its shorter side. A primal path is a continuous curve contained in $\omega$, that is, running along uncut segments of primal lines and along bridges, and a dual path is a continuous curve contained in $\omega^*$. A primal crossing of $\mathcal R$ is a primal path in $\mathcal R$ joining its two shorter sides, and a dual crossing is a dual path in $\mathcal R$ joining them. Let $C(\mathcal R)$ be the event that $\mathcal R$ has a primal crossing, and $C^*(\mathcal R)$ the event that it has a dual crossing. The aspect ratio of $\mathcal R$ is the ratio of its longer side to its shorter side (Figure~\ref{fig:crossing}).

\paragraph{The idea of Lemma~\ref{lem:critical-loops}.} Part (1) says that at the self-dual point the influence of the boundary dies out. Whatever boundary condition is imposed far away, the configuration in a fixed bounded region has, in the limit, the same law. Part (2) expresses scale invariance. The probability that a rectangle is crossed depends, up to constants, only on its shape and not on its size, so it stays bounded away from $0$ at all scales. Away from the self-dual point this fails, and crossings of large rectangles become either very likely or very unlikely. The bound in (2) is stated conditionally on everything outside a region $\Delta'$ containing the enlarged rectangle. Whatever happens outside acts on $\Delta'$ only as a boundary condition, and the worst boundary conditions are the free one for primal crossings and the wired one for dual crossings. So the crossing has a fresh chance of at least $c(\rho)$ whatever the outside information.

\begin{lemma}[The critical random-cluster measure]
\label{lem:critical-loops}
Let $q=\delta^2\in[2,4)$ and $(\lambda_{\mathrm c},\lambda_{\mathrm b})=(\delta^{-2},1)$. The constants below depend on $k$.
\begin{enumerate}[leftmargin=*]
\item The free and wired infinite volume measures coincide. The measure $\phi=\phi^0=\phi^1$ is invariant under integer space translations and real time translations. For every bounded $B\subset\Z\times\mathbb R$ and $\varepsilon>0$ there is $R$ such that, for every finite volume $\Delta$ containing the $R$-neighbourhood of $B$ and every boundary condition $\zeta$, the restrictions of $\phi^\zeta_\Delta$ and $\phi$ to $\mathcal F_B$ differ by at most $\varepsilon$ in total variation.
\item (Russo--Seymour--Welsh (RSW) property.) For every $\rho>1$ there are $c(\rho)>0$ and $s(\rho)>0$ with the following property. Let $\mathcal R$ be a rectangle of aspect ratio $\rho$, with its longer side in either the space or the time direction, and with shorter side at least $s(\rho)$. Let $\Delta'\subseteq\Delta$ be finite volumes with $\widehat{\mathcal R}\subseteq\Delta'$, and let $\zeta$ be a boundary condition on $\Delta$. Then, almost surely,
\[
 \phi^\zeta_\Delta\bigl(C(\mathcal R)\mid\mathcal F_{\Delta\setminus\Delta'}\bigr)\ge c(\rho),\qquad
 \phi^\zeta_\Delta\bigl(C^*(\mathcal R)\mid\mathcal F_{\Delta\setminus\Delta'}\bigr)\ge c(\rho).
\]
\end{enumerate}
\end{lemma}
\begin{proof}
For (1), in the notation of \cite{DCLM} the cut and bridge intensities are $\lambda$ and $\mu$, and for $q\in[1,4]$ and $\mu/\lambda=q$, Duminil-Copin, Li and Manolescu prove that $\phi^0=\phi^1$ and that there is no infinite cluster \cite[Theorem~1.5]{DCLM}. Here $\lambda_{\mathrm b}/\lambda_{\mathrm c}=q$. The limit $\phi^0$ does not depend on the increasing sequence of finite volumes, by monotonicity, so it is invariant under the translations that preserve the model.

For the total variation bound, let $Q$ be the $R$-neighbourhood of $B$ and $\Delta\supseteq Q$. By (iii) the restriction of $\phi^\zeta_\Delta$ to $\mathcal F_Q$ is an average of measures $\phi^\xi_Q$, and by (ii) each $\phi^\xi_Q$, as well as the restriction of $\phi$, lies stochastically between $\phi^0_Q$ and $\phi^1_Q$. Let $\mu$ be one of these measures. By Strassen's theorem \cite[Theorem~11]{Str65} there is a coupling of $\mu$ and $\phi^1_Q$ preserving the order, that is, a probability measure on pairs of configurations with marginals $\mu$ and $\phi^1_Q$ that is carried by the set of pairs $(\omega,\omega')$ with $\omega\le\omega'$. We write $(\omega^\mu,\omega^1)$ for a pair distributed according to this coupling, so that $\omega^\mu\le\omega^1$ holds with probability one, see also \cite[proof of Lemma~II.2]{AN94}. Write $N_{\mathrm c}$ and $N_{\mathrm b}$ for the numbers of cuts and bridges in $B$. If $\omega^\mu$ and $\omega^1$ differ in $B$, then $N_{\mathrm c}(\omega^\mu)-N_{\mathrm c}(\omega^1)+N_{\mathrm b}(\omega^1)-N_{\mathrm b}(\omega^\mu)\ge1$, and each difference is nonnegative. Since $N_{\mathrm c}$ is decreasing and $N_{\mathrm b}$ increasing,
\[
 \mathbb P(\omega^\mu|_B\ne\omega^1|_B)\le\eta_R:=\mathbb E_{\phi^0_Q}[N_{\mathrm c}]-\mathbb E_{\phi^1_Q}[N_{\mathrm c}]+\mathbb E_{\phi^1_Q}[N_{\mathrm b}]-\mathbb E_{\phi^0_Q}[N_{\mathrm b}].
\]
These counts are uniformly integrable. Adding a cut changes the cluster count $k^\zeta$ in the definition of $\phi^\zeta_\Delta$ by $0$ or $1$, and adding a bridge changes it by $0$ or $-1$. So under every $\phi^\zeta_\Delta$ the conditional intensity of cuts lies in $[\lambda_{\mathrm c},q\lambda_{\mathrm c}]$ and that of bridges in $[\lambda_{\mathrm b}/q,\lambda_{\mathrm b}]$. Here $q\lambda_{\mathrm c}=\lambda_{\mathrm b}=1$, so $N_{\mathrm c}$ and $N_{\mathrm b}$ are stochastically dominated by Poisson variables with mean the total length of the primal lines and of the bridge gaps in $B$, uniformly in $\Delta$ and $\zeta$. As $R\to\infty$ the expectations therefore converge monotonically to the expectations under $\phi^0$ and $\phi^1$, and these coincide. So $\eta_R\to0$. Comparing $\phi^\xi_Q$ and $\phi$ through $\phi^1_Q$ bounds their total variation distance on $\mathcal F_B$ by $2\eta_R$, which proves (1).

For (2), Duminil-Copin, Li and Manolescu prove in \cite[Section~5.2]{DCLM}, for $q\in[1,4]$ and $\mu/\lambda=q$, that the quantum random-cluster model, the model $\phi^\zeta_\Delta$ on $\Z\times\mathbb R$ above, has the RSW property in the form (BXP) of \cite[Lemma~3.1]{DCLM}. For every $\rho>1$ and $\nu_0>0$ there is $\delta_1>0$ such that, for all large $n$, primal horizontal and vertical crossings of the rectangles of sides $2\rho n$ and $2n$ have probability at least $\delta_1$ under the free boundary condition on the rectangle enlarged by the margin $\nu_0n$, and the dual crossings have probability at least $\delta_1$ under the wired boundary condition there. They obtain uniform crossing bounds on deformed isoradial square lattices and transfer them to the quantum random-cluster model \cite[Propositions~5.1 and~5.3]{DCLM}. We take $\nu_0=\frac12$, which is the margin of $\widehat{\mathcal R}$. Their intensities have ratio $q$ but a different common scale. Multiplying both intensities by $\alpha$ dilates the configuration vertically by $1/\alpha$ \cite[Section~1.3]{DCLM}, which changes aspect ratios by a fixed factor. Since (BXP) holds for every $\rho>1$, these bounds hold for crossings in the long direction. Consider a crossing in the short direction of a rectangle $[-a,a]\times[-b,b]$ with $b>a$, say. The subrectangle $[-a,a]\times[-\frac a2,\frac a2]$ has aspect ratio $2$, and its enlargement lies in $\widehat{\mathcal R}$. Every crossing of it in the long direction is a crossing of $\mathcal R$ in the short direction. So the bounds hold in both directions, at $(\lambda_{\mathrm c},\lambda_{\mathrm b})=(\delta^{-2},1)$ with a constant $c(\rho)>0$. The model is invariant under integer space translations and real time translations, so the bounds hold for all rectangles. Now let $\mathcal R,\Delta',\Delta,\zeta$ be as in (2). By (iii), conditionally on $\mathcal F_{\Delta\setminus\Delta'}$ and then on the configuration in $\Delta'\setminus\widehat{\mathcal R}$, the configuration in $\widehat{\mathcal R}$ has law $\phi^\xi_{\widehat{\mathcal R}}$ for an induced boundary condition $\xi$. The event $C(\mathcal R)$ is increasing, so by (ii) $\phi^\xi_{\widehat{\mathcal R}}(C(\mathcal R))\ge\phi^0_{\widehat{\mathcal R}}(C(\mathcal R))\ge c(\rho)$. The event $C^*(\mathcal R)$ is decreasing, so the least favourable boundary condition for it is the wired one, and $\phi^\xi_{\widehat{\mathcal R}}(C^*(\mathcal R))\ge\phi^1_{\widehat{\mathcal R}}(C^*(\mathcal R))\ge c(\rho)$. Averaging over the configuration in $\Delta'\setminus\widehat{\mathcal R}$ gives (2).
\end{proof}

We define the infinite volume loop law $\nu$ as the image of $\phi$ under $\omega\mapsto\Gamma$. Since $\nu_{N,T}$ is the image of $\phi^0_{\Delta_{N,T}}$, Lemma~\ref{lem:critical-loops}(1) makes $\nu$ the limit of $\nu_{N,T}$ on events depending on the configuration in a bounded region.

\subsection{Separating loops}

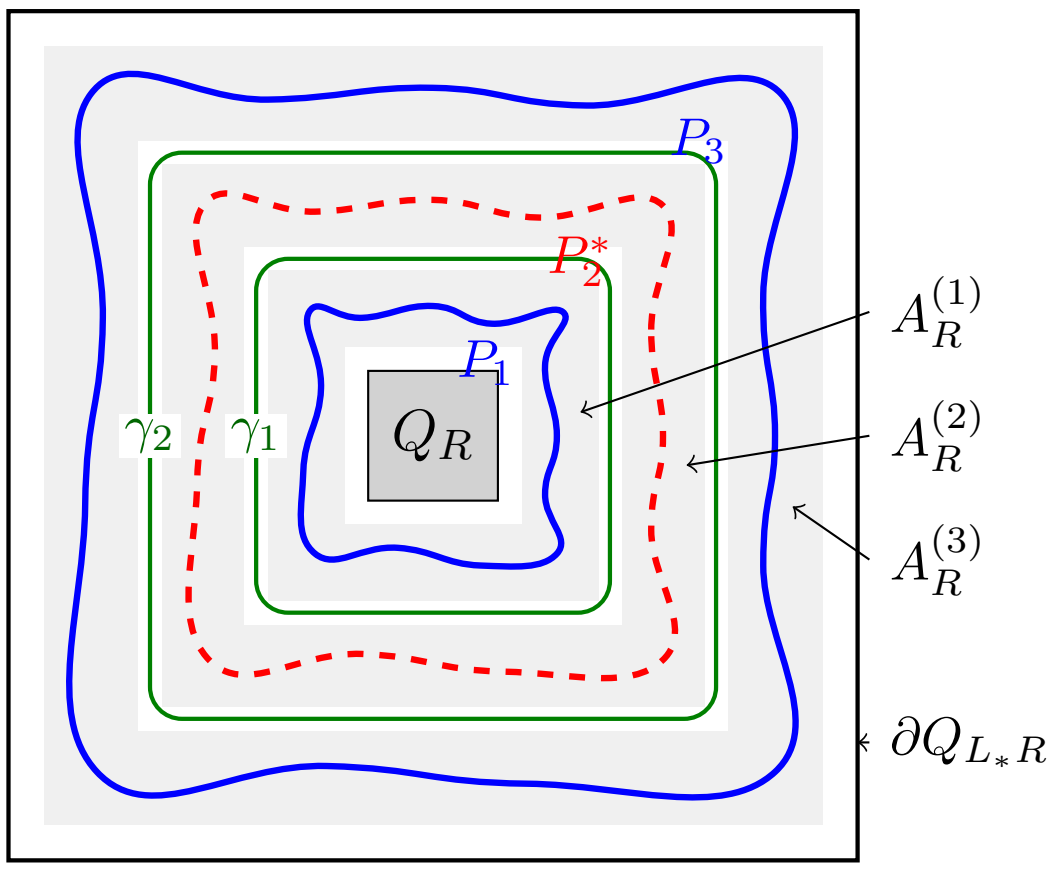
\begin{figure}[ht]
\centering
\begin{tikzpicture}[scale=0.78]
 \fill[gray!12] (-3.3,-3.3) rectangle (3.3,3.3);
 \fill[white] (-2.5,-2.5) rectangle (2.5,2.5);
 \fill[gray!12] (-2.3,-2.3) rectangle (2.3,2.3);
 \fill[white] (-1.6,-1.6) rectangle (1.6,1.6);
 \fill[gray!12] (-1.4,-1.4) rectangle (1.4,1.4);
 \fill[white] (-0.75,-0.75) rectangle (0.75,0.75);
 \fill[gray!35] (-0.55,-0.55) rectangle (0.55,0.55);
 \node at (0,0) {$Q_R$};
 \draw[thick] (-3.6,-3.6) rectangle (3.6,3.6);
 \draw[thin] (-0.55,-0.55) rectangle (0.55,0.55);
 \draw[very thick,blue] plot[smooth cycle,tension=0.8] coordinates {(1.05,0) (0.95,0.6) (1.1,1.05) (0.5,0.95) (0,1.1) (-0.6,1.0) (-1.05,1.05) (-0.95,0.4) (-1.1,-0.3) (-1.0,-1.0) (-0.3,-0.95) (0.4,-1.1) (1.05,-1.0) (0.95,-0.5)};
 \draw[very thick,red,dashed] plot[smooth cycle,tension=0.8] coordinates {(1.95,0) (1.85,0.9) (1.95,1.95) (0.9,1.85) (0,2.0) (-1.0,1.9) (-1.95,1.95) (-1.85,0.7) (-2.0,-0.4) (-1.9,-1.9) (-0.6,-1.85) (0.6,-2.0) (1.95,-1.9) (1.85,-0.8)};
 \draw[very thick,blue] plot[smooth cycle,tension=0.8] coordinates {(2.9,0) (2.8,1.2) (2.95,2.9) (1.3,2.8) (0,2.95) (-1.4,2.85) (-2.9,2.9) (-2.8,1.0) (-2.95,-0.6) (-2.85,-2.85) (-0.9,-2.8) (0.8,-2.95) (2.9,-2.85) (2.8,-1.1)};
 \draw[thick,green!50!black,rounded corners=6pt] (-1.5,-1.5) rectangle (1.5,1.5);
 \draw[thick,green!50!black,rounded corners=6pt] (-2.4,-2.4) rectangle (2.4,2.4);
 \node[blue] at (0.45,0.62) {\small $P_1$};
 \node[red] at (1.25,1.48) {\small $P_2^*$};
 \node[blue] at (2.25,2.5) {\small $P_3$};
 \node[green!40!black,fill=white,inner sep=1pt] at (-1.5,0) {\small $\gamma_1$};
 \node[green!40!black,fill=white,inner sep=1pt] at (-2.4,0) {\small $\gamma_2$};
 \node[right] at (3.7,1.05) {\small $A_R^{(1)}$};
 \draw[->] (3.7,1.05) -- (1.25,0.2);
 \node[right] at (3.7,0) {\small $A_R^{(2)}$};
 \draw[->] (3.7,0) -- (2.15,-0.25);
 \node[right] at (3.7,-1.05) {\small $A_R^{(3)}$};
 \draw[->] (3.7,-1.05) -- (3.05,-0.6);
 \node[right] at (3.7,-2.6) {\small $\partial Q_{L_*R}$};
 \draw[->] (3.7,-2.6) -- (3.6,-2.6);
\end{tikzpicture}
\caption{Schematic picture of Lemma~\ref{lem:separating}, not to scale. The shaded bands are the subannuli $A_R^{(1)},A_R^{(2)},A_R^{(3)}$ of $A_R=Q_{L_*R}\setminus Q_R$. On $G_R$ there are primal circuits (solid) in the first and third and a dual circuit (dashed) in the second. The loop $\gamma_1$ bounding the primal cluster of the inner circuit from outside, and the loop $\gamma_2$ bounding the dual cluster of the dual circuit from outside, both surround $Q_R$, lie in $A_R$, and are distinct.}
\label{fig:separating}
\end{figure}

Fix a centre $z\in\Z\times\mathbb R$ and abbreviate $Q_r=Q_r(z)$. For $R>0$ put $L_*=180$, $r_i=2\cdot6^{i-1}R$ for $i=1,2,3$, and
\[
 A_R=Q_{L_*R}\setminus Q_R,\qquad A_R^{(i)}=Q_{5r_i/2}\setminus Q_{r_i/2},\qquad C_R^{(i)}=Q_{2r_i}\setminus Q_{r_i}.
\]
So $A_R^{(1)}=Q_{5R}\setminus Q_R$, $A_R^{(2)}=Q_{30R}\setminus Q_{6R}$ and $A_R^{(3)}=Q_{180R}\setminus Q_{36R}$ are disjoint subannuli of $A_R$, ordered from inside to outside. Each contains the core annulus $C_R^{(i)}$ with a buffer of width $r_i/2$ on each side. A primal circuit around $Q_R$ is a closed primal path (Section~\ref{sec:loops}) winding once around $Q_R$, and a dual circuit around $Q_R$ is a closed dual path winding once around $Q_R$. Let $P_1$ and $P_3$ be the events that there is a primal circuit around $Q_R$ in $C_R^{(1)}$ and in $C_R^{(3)}$, and $P_2^*$ the event that there is a dual circuit around $Q_R$ in $C_R^{(2)}$. Put $G_R=P_1\cap P_2^*\cap P_3$. See Figure~\ref{fig:separating}.

\paragraph{The idea of Lemma~\ref{lem:separating}.} Four rectangles placed around a square annulus, each crossed the long way, link up into a closed curve around the inner square (Figure~\ref{fig:circuit}). By Lemma~\ref{lem:critical-loops}(2) each crossing has probability at least $c(4)$, and by positive association all four occur together with probability at least $c(4)^4$. A loop of $\Gamma$ is the interface between a primal and a dual cluster, so to obtain loops around $Q_R$ one needs circuits of both kinds. A primal circuit enclosed by a dual circuit gives one loop between them, and a further primal circuit outside gives a second loop. Since every bound holds conditionally on the outside, the three circuits can be placed in three disjoint annuli and their probabilities multiplied.

\begin{lemma}[Uniform conditional separating loops]
\label{lem:separating}
Put $c_*=c(4)^4$ and $p_*=c_*^3$, with $c(4)$ from Lemma~\ref{lem:critical-loops}(2), and $R_*=s(4)/2$. Let $R\ge R_*$, let $\Delta$ be a finite volume containing $Q_{L_*R}$, and let $\zeta$ be a boundary condition on $\Delta$. Then, almost surely,
\begin{equation}
 \phi^\zeta_\Delta(G_R\mid\mathcal F_{\Delta\setminus A_R})\ge p_* .
 \label{eq:annulus}
\end{equation}
The event $G_R$ depends only on the configuration in $A_R$. On $G_R$ the loops of $\Gamma$ include two distinct loops contained in $A_R$ and surrounding $Q_R$. The bound \eqref{eq:annulus} also holds for the infinite volume measure $\phi$.
\end{lemma}
\begin{proof}
Fix $i$ and write $r=r_i$. Consider the four rectangles
\[
 [-2r,2r]\times[r,2r],\quad[-2r,2r]\times[-2r,-r],\quad[-2r,-r]\times[-2r,2r],\quad[r,2r]\times[-2r,2r]
\]
relative to $z$, each of aspect ratio $4$ and shorter side $r\ge2R\ge s(4)$. Primal crossings of the first two in the long direction are horizontal and of the last two vertical. Their union contains a primal circuit in $C_R^{(i)}$ around $Q_r\supseteq Q_R$, because each horizontal crossing meets both vertical crossings. The enlargement of the first rectangle is $[-\frac94r,\frac94r]\times[\frac34r,\frac94r]$, and all four enlargements lie in the interior of $A_R^{(i)}$ (Figure~\ref{fig:circuit}). By (iii), conditionally on $\mathcal F_{\Delta\setminus A_R^{(i)}}$ the configuration in $A_R^{(i)}$ has law $\phi^\xi_{A_R^{(i)}}$ for an induced boundary condition $\xi$. By Lemma~\ref{lem:critical-loops}(2) with $\Delta'=A_R^{(i)}$, each of the four crossings has conditional probability at least $c(4)$. They are increasing events, so the FKG inequality (i) for $\phi^\xi_{A_R^{(i)}}$ bounds the conditional probability of all four, and hence of a primal circuit, below by $c(4)^4=c_*$. For $i=2$ the same argument with dual crossings, which are decreasing events and hence also positively correlated, gives the dual circuit. Thus, almost surely,
\[
 \phi^\zeta_\Delta\bigl(P_1\mid\mathcal F_{\Delta\setminus A_R^{(1)}}\bigr)\ge c_*,\qquad
 \phi^\zeta_\Delta\bigl(P_2^*\mid\mathcal F_{\Delta\setminus A_R^{(2)}}\bigr)\ge c_*,\qquad
 \phi^\zeta_\Delta\bigl(P_3\mid\mathcal F_{\Delta\setminus A_R^{(3)}}\bigr)\ge c_*.
\]
Only their combination \eqref{eq:annulus} is used later.

A primal or dual path in $C_R^{(i)}$ consists of segments of lines, bridges and cuts in $C_R^{(i)}$. So $P_1$, $P_2^*$ and $P_3$ depend only on the configurations in $C_R^{(1)}$, $C_R^{(2)}$ and $C_R^{(3)}$, and $G_R$ only on the configuration in $A_R$. The events $P_1$ and $P_3$ are measurable with respect to $\mathcal F_{\Delta\setminus A_R^{(2)}}$, and $\mathcal F_{\Delta\setminus A_R}\subseteq\mathcal F_{\Delta\setminus A_R^{(2)}}$. Iterating conditional expectations, first over $A_R^{(2)}$ and then over $A_R^{(1)}$,
\[
 \phi^\zeta_\Delta(P_1\cap P_2^*\cap P_3\mid\mathcal F_{\Delta\setminus A_R})
 \ge c_*\,\phi^\zeta_\Delta(P_1\cap P_3\mid\mathcal F_{\Delta\setminus A_R})
 \ge c_*^2\,\phi^\zeta_\Delta(P_3\mid\mathcal F_{\Delta\setminus A_R})\ge c_*^3 .
\]

On $G_R$, let $\gamma_1$ be the outer boundary of a thin neighbourhood of the primal cluster containing the inner circuit, and $\gamma_2$ the outer boundary of a thin neighbourhood of the dual cluster containing the dual circuit. These are loops of $\Gamma$ surrounding $Q_R$ (Figure~\ref{fig:gammas}). A primal cluster cannot cross a dual circuit, and a dual cluster cannot cross a primal circuit. So $\gamma_1$ lies inside the dual circuit and $\gamma_2$ between the dual circuit and the outer primal circuit. Both are contained in $A_R$, and they are distinct because the dual circuit separates them.

For $\phi$, let $F$ be an event depending on the configuration in a bounded region outside $A_R$. The bound $\phi^0_\Delta(G_R\cap F)\ge p_*\phi^0_\Delta(F)$ passes to the limit by Lemma~\ref{lem:critical-loops}(1). It extends to all $F\in\mathcal F_{A_R^c}$ by the monotone class theorem, which is \eqref{eq:annulus} for $\phi$.
\end{proof}

\begin{figure}[ht]
\centering
\begin{tikzpicture}[scale=0.85]
 \fill[gray!12] (-2.5,-2.5) rectangle (2.5,2.5);
 \fill[white] (-0.5,-0.5) rectangle (0.5,0.5);
 \draw[thin] (-2.5,-2.5) rectangle (2.5,2.5);
 \draw[thin] (-0.5,-0.5) rectangle (0.5,0.5);
 \draw[thick] (-1,-1) rectangle (1,1);
 \draw[thick] (-2,-2) rectangle (2,2);
 \draw[thick] (-1,1) -- (-1,2) (1,1) -- (1,2) (-1,-1) -- (-1,-2) (1,-1) -- (1,-2);
 \draw[thick] (-2,1) -- (-1,1) (1,1) -- (2,1) (-2,-1) -- (-1,-1) (1,-1) -- (2,-1);
 \draw[dashed] (-2.25,0.75) rectangle (2.25,2.25);
 \draw[very thick,blue] plot[smooth,tension=0.7] coordinates {(-2,1.5) (-1.2,1.7) (-0.4,1.3) (0.5,1.6) (1.3,1.35) (2,1.6)};
 \draw[very thick,blue] plot[smooth,tension=0.7] coordinates {(-2,-1.4) (-1.1,-1.65) (-0.3,-1.3) (0.6,-1.7) (1.4,-1.3) (2,-1.5)};
 \draw[very thick,blue] plot[smooth,tension=0.7] coordinates {(-1.5,-2) (-1.7,-1.1) (-1.3,-0.2) (-1.6,0.7) (-1.35,1.4) (-1.55,2)};
 \draw[very thick,blue] plot[smooth,tension=0.7] coordinates {(1.6,-2) (1.3,-1.2) (1.65,-0.3) (1.35,0.6) (1.7,1.3) (1.45,2)};
 \node at (0,0) {\small $Q_{r/2}$};
 \node at (0,0.75) {\small $Q_r$};
 \node[fill=white,inner sep=1pt] at (0,2.25) {\scriptsize $\widehat{\mathcal R}$};
 \node[right] at (2.6,1.5) {\small $\mathcal R=[-2r,2r]\times[r,2r]$};
 \node[right] at (2.6,-1.5) {\small $A_R^{(i)}=Q_{5r/2}\setminus Q_{r/2}$};
 \node[right] at (2.6,0) {\small $C_R^{(i)}=Q_{2r}\setminus Q_r$};
\end{tikzpicture}
\caption{The core annulus $C_R^{(i)}$ covered by four rectangles of aspect ratio $4$. Primal crossings of all four in the long direction contain a circuit around $Q_r$. The enlargement $\widehat{\mathcal R}$ of the top rectangle is dashed, and all four enlargements lie in the interior of the buffered annulus $A_R^{(i)}$ (shaded).}
\label{fig:circuit}
\end{figure}
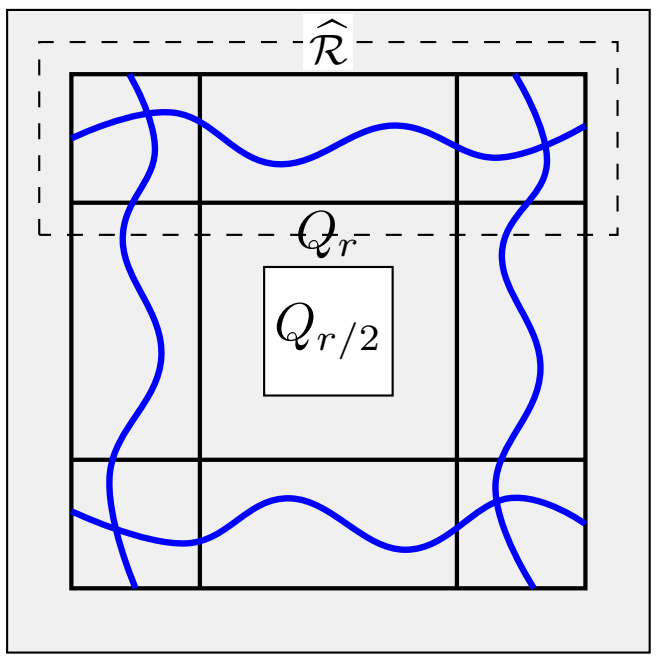

\begin{figure}[ht]
\centering
\begin{tikzpicture}[scale=0.9]
 \draw[very thick] (0,0.8) -- (9,0.8);
 \draw[very thick] (2,0.8) -- (2,1.7) (1.6,1.7) -- (2.6,1.7) (5.5,0.8) -- (5.5,1.4);
 \draw[very thick,red,dashed] (0,3.2) -- (9,3.2);
 \draw[very thick,red,dashed] (3.5,3.2) -- (3.5,4.1) (3.1,4.1) -- (4.0,4.1) (7,3.2) -- (7,3.7) (7,3.2) -- (7,2.6);
 \draw[very thick] (0,5.6) -- (9,5.6);
 \draw[thick,green!50!black] (0,1.0) -- (1.85,1.0) -- (1.85,1.55) -- (1.45,1.55) -- (1.45,1.9) -- (2.75,1.9) -- (2.75,1.55) -- (2.15,1.55) -- (2.15,1.0) -- (5.35,1.0) -- (5.35,1.6) -- (5.65,1.6) -- (5.65,1.0) -- (9,1.0);
 \draw[thick,green!50!black] (0,3.4) -- (3.35,3.4) -- (3.35,3.95) -- (2.95,3.95) -- (2.95,4.3) -- (4.15,4.3) -- (4.15,3.95) -- (3.65,3.95) -- (3.65,3.4) -- (6.85,3.4) -- (6.85,3.9) -- (7.15,3.9) -- (7.15,3.4) -- (9,3.4);
 \node[right] at (9,0.8) {\scriptsize inner primal circuit};
 \node[right] at (9,1.1) {\scriptsize $\gamma_1$};
 \node[right] at (9,3.1) {\scriptsize dual circuit};
 \node[right] at (9,3.45) {\scriptsize $\gamma_2$};
 \node[right] at (9,5.6) {\scriptsize outer primal circuit};
 \node at (4.5,0.2) {\scriptsize towards $Q_R$};
 \node at (4.5,6.1) {\scriptsize towards $\partial Q_{L_*R}$};
\end{tikzpicture}
\caption{The proof of Lemma~\ref{lem:separating} near the circuits, with the annulus $A_R$ cut open along a radius, so that the left and right edges are identified and the radial direction points up. The primal cluster of the inner circuit (solid) cannot cross the dual circuit (dashed), and the outer boundary $\gamma_1$ of a thin neighbourhood of this cluster lies between them. The outer boundary $\gamma_2$ of a thin neighbourhood of the dual cluster of the dual circuit lies between the dual circuit and the outer primal circuit. The dual circuit separates $\gamma_1$ from $\gamma_2$.}
\label{fig:gammas}
\end{figure}
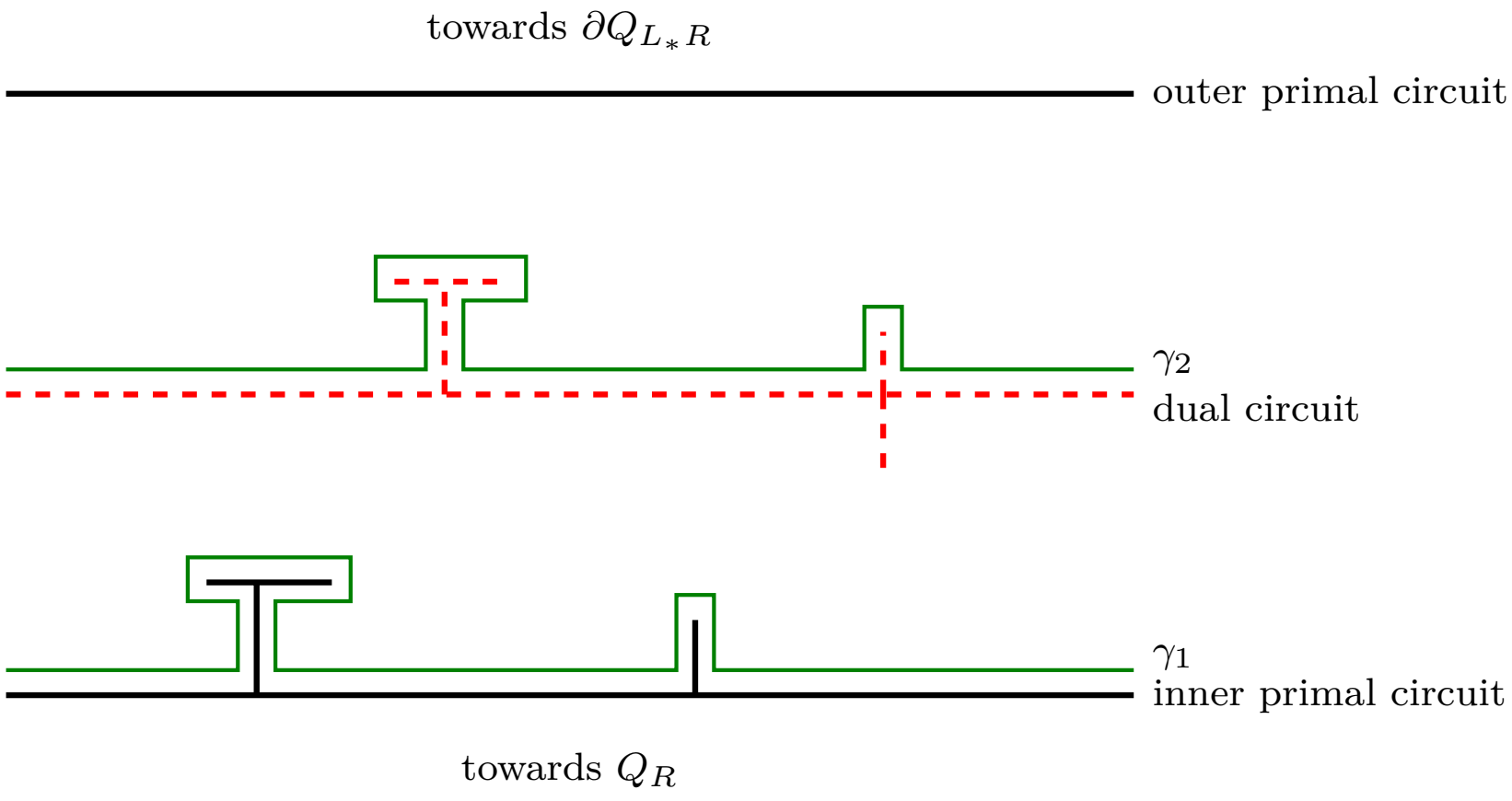

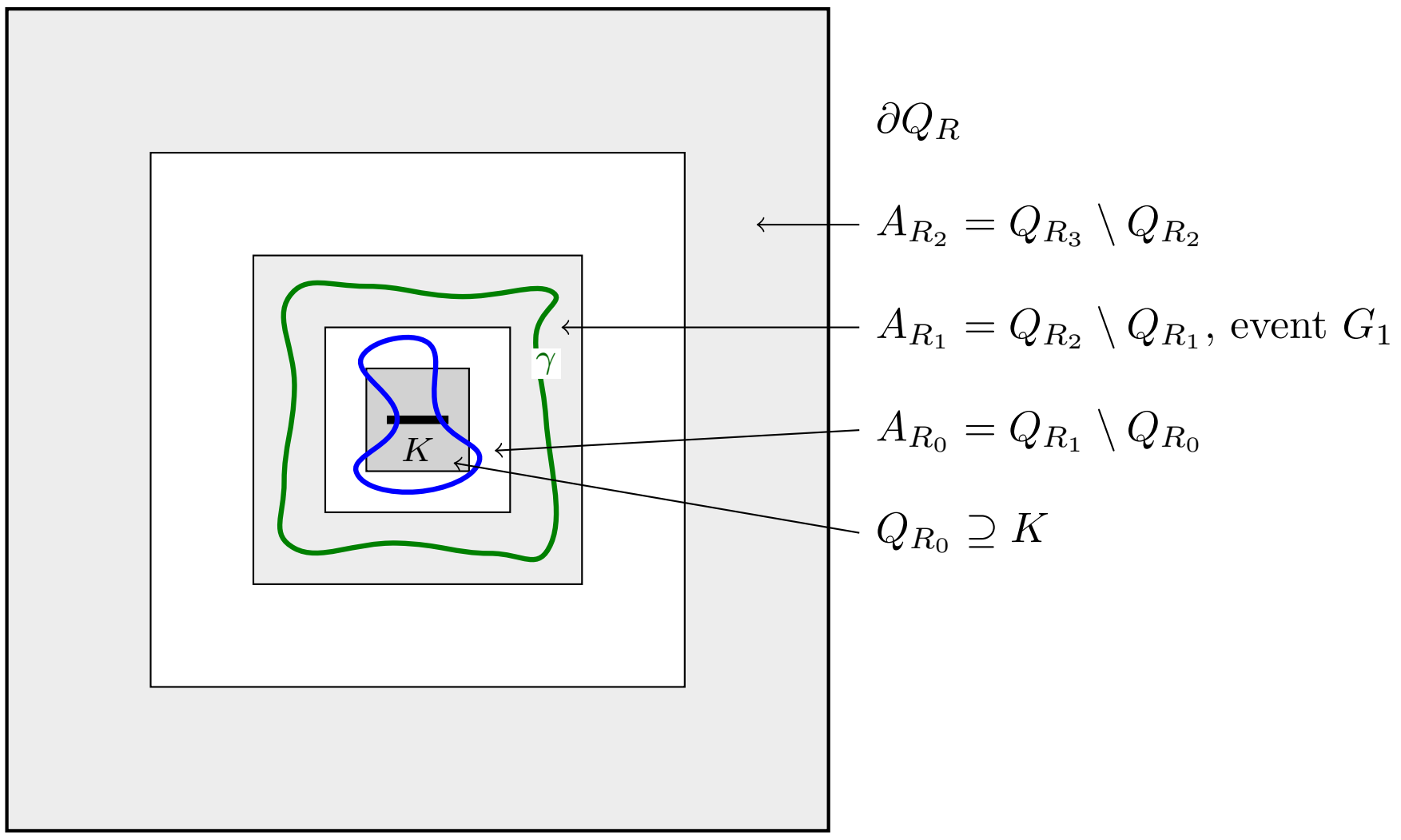
\begin{figure}[ht]
\centering
\begin{tikzpicture}[scale=0.85]
 \fill[gray!14] (-4,-4) rectangle (4,4);
 \fill[white] (-2.6,-2.6) rectangle (2.6,2.6);
 \fill[gray!14] (-1.6,-1.6) rectangle (1.6,1.6);
 \fill[white] (-0.9,-0.9) rectangle (0.9,0.9);
 \fill[gray!35] (-0.5,-0.5) rectangle (0.5,0.5);
 \draw (-0.5,-0.5) rectangle (0.5,0.5);
 \draw (-0.9,-0.9) rectangle (0.9,0.9);
 \draw (-1.6,-1.6) rectangle (1.6,1.6);
 \draw (-2.6,-2.6) rectangle (2.6,2.6);
 \draw[thick] (-4,-4) rectangle (4,4);
 \draw[line width=2pt] (-0.3,0) -- (0.3,0);
 \node[below] at (0,-0.02) {\scriptsize $K$};
 \draw[very thick,blue] plot[smooth cycle,tension=0.8] coordinates {(-0.2,0) (-0.55,0.6) (0.1,0.75) (0.2,0.05) (0.6,-0.4) (0,-0.7) (-0.6,-0.5)};
 \draw[very thick,green!50!black] plot[smooth cycle,tension=0.8] coordinates {(1.25,0) (1.15,0.8) (1.3,1.25) (0.4,1.2) (-0.5,1.3) (-1.25,1.2) (-1.2,0.3) (-1.3,-0.6) (-1.2,-1.25) (-0.2,-1.2) (0.7,-1.3) (1.3,-1.2)};
 \node[green!40!black,fill=white,inner sep=1pt] at (1.25,0.55) {\scriptsize $\gamma$};
 \node[right] at (4.3,2.9) {\small $\partial Q_R$};
 \node[right] at (4.3,1.9) {\small $A_{R_2}=Q_{R_3}\setminus Q_{R_2}$};
 \draw[->] (4.3,1.9) -- (3.3,1.9);
 \node[right] at (4.3,0.9) {\small $A_{R_1}=Q_{R_2}\setminus Q_{R_1}$, event $G_1$};
 \draw[->] (4.3,0.9) -- (1.4,0.9);
 \node[right] at (4.3,-0.1) {\small $A_{R_0}=Q_{R_1}\setminus Q_{R_0}$};
 \draw[->] (4.3,-0.1) -- (0.75,-0.3);
 \node[right] at (4.3,-1.1) {\small $Q_{R_0}\supseteq K$};
 \draw[->] (4.3,-1.1) -- (0.35,-0.42);
\end{tikzpicture}
\caption{Schematic picture of Lemma~\ref{lem:looptail}, not to scale. Between $K\subseteq Q_{R_0}$ and $\partial Q_R$ lie the annuli $A_{R_j}=Q_{R_{j+1}}\setminus Q_{R_j}$, $R_j=L_*^jR_0$, alternately shaded. Each is a copy of Figure~\ref{fig:separating} at scale $R_j$. Here the event $G_1$ occurs, and its loop $\gamma$ surrounds $K$. Since loops are disjoint, every loop meeting $K$ (solid, blue) stays inside $\gamma$. Each $G_j$ has conditional probability at least $p_*$ given the configuration outside $A_{R_j}$, so all of them fail with probability at most $(1-p_*)^m$.}
\label{fig:nested}
\end{figure}

\paragraph{The idea of Lemma~\ref{lem:looptail}.} A single annulus contains a pair of separating loops only with probability at least $p_*$, which does not tend to $1$. Between $K$ and distance $R$ there are about $\log R/\log L_*$ disjoint annuli at geometrically growing scales, each a copy of Figure~\ref{fig:separating} at its own scale (Figure~\ref{fig:nested}). Each gives a fresh chance of at least $p_*$, whatever happens outside it. So all of them fail with probability at most $(1-p_*)^{\log R/\log L_*}$, which is a negative power of $R$. Loops do not cross, so a single separating loop around $K$ confines every loop meeting $K$ to its inside. Mixing has a similar explanation (Remark~\ref{rem:mixing}).

\begin{lemma}[Loop tails and mixing]
\label{lem:looptail}
Let $K\subset\Z\times\mathbb R$ be bounded, and choose $z$ and $R_0\ge R_*$ with $K\subseteq Q_{R_0}(z)$. Put $\beta_0=-\log(1-p_*)/\log L_*$ and $C_K=R_0^{\beta_0}/(1-p_*)$. For every $R\ge R_0$,
\begin{equation}
 \mathbb P(\text{some loop meeting $K$ leaves $Q_R(z)$})\le C_KR^{-\beta_0},
 \label{eq:looptail}
\end{equation}
for $\mathbb P=\nu$, and for $\mathbb P=\nu_{N,T}$ whenever $Q_R(z)\subseteq\Delta_{N,T}$. Under $\nu$ every loop meeting $K$ is finite almost surely. The measure $\phi$, and hence the loop law $\nu$, is mixing under space translations by one primal column and under time translations.
\end{lemma}
\begin{proof}
The events below are events of the random-cluster configuration, and we write $\mathbb P$ also for $\phi$ and for $\phi^0_{\Delta_{N,T}}$, whose images under $\omega\mapsto\Gamma$ are $\nu$ and $\nu_{N,T}$. Let $m$ be the largest integer with $L_*^mR_0\le R$. Put $R_j=L_*^jR_0$ for $0\le j<m$. The annuli $A_{R_j}=Q_{R_{j+1}}\setminus Q_{R_j}$ have disjoint interiors and lie in $Q_R\setminus K$. Let $G_j=G_{R_j}$ be the event of Lemma~\ref{lem:separating} in $A_{R_j}$. By that lemma, $G_j$ depends only on the configuration in $A_{R_j}$, so $G_0^c,\dots,G_{m-2}^c$ are measurable with respect to the $\sigma$-algebra $\mathcal F_{A_{R_{m-1}}^c}$ of the configuration outside $A_{R_{m-1}}$. By the tower property of conditional expectations,
\[
 \mathbb P\Bigl(\bigcap_{j<m}G_j^c\Bigr)=\mathbb E\Bigl[\mathbf 1_{\bigcap_{j<m-1}G_j^c}\,\mathbb P\bigl(G_{m-1}^c\mid\mathcal F_{A_{R_{m-1}}^c}\bigr)\Bigr].
\]
By \eqref{eq:annulus}, $\mathbb P(G_{m-1}\mid\mathcal F_{A_{R_{m-1}}^c})\ge p_*$ almost surely, so the conditional probability inside the expectation is at most $1-p_*$, and
\[
 \mathbb P\Bigl(\bigcap_{j<m}G_j^c\Bigr)\le(1-p_*)\,\mathbb P\Bigl(\bigcap_{j<m-1}G_j^c\Bigr).
\]
By induction $\mathbb P(\bigcap_jG_j^c)\le(1-p_*)^m$. No independence of the annuli is used. By the choice of $m$, $m>\log(R/R_0)/\log L_*-1$, and since $1-p_*<1$ this gives
\[
 (1-p_*)^m\le(1-p_*)^{-1}(1-p_*)^{\log(R/R_0)/\log L_*}=(1-p_*)^{-1}(R/R_0)^{-\beta_0}=C_KR^{-\beta_0}.
\] On $G_j$ there is a loop $\gamma\subset A_{R_j}$ surrounding $Q_{R_j}\supseteq K$. Loops are disjoint, so a loop meeting $K$ and leaving $Q_R$ would meet both the inside and the outside of $\gamma$, hence equal $\gamma$, which does not meet $K$. This proves \eqref{eq:looptail}, and letting $R\to\infty$ gives finiteness.

For mixing, let $\tau$ be a translation by $n$ primal columns or by time $n$, and let $F,G$ be bounded functions of the configuration in a bounded region $B$. Given $\varepsilon>0$ choose $R$ as in Lemma~\ref{lem:critical-loops}(1) and let $Q$ be the $R$-neighbourhood of $B$. For $|n|$ large, $\tau B$ lies outside $Q$, so $G\circ\tau$ is $\mathcal F_{Q^c}$-measurable. For a large finite volume $\Delta$, (iii) and Lemma~\ref{lem:critical-loops}(1) give
\[
 \bigl|\mathbb E_{\phi^0_\Delta}[F\,G\circ\tau]-\mathbb E_\phi[F]\,\mathbb E_{\phi^0_\Delta}[G\circ\tau]\bigr|
 =\bigl|\mathbb E_{\phi^0_\Delta}\bigl[(\mathbb E_{\phi^\xi_Q}[F]-\mathbb E_\phi[F])\,G\circ\tau\bigr]\bigr|\le\varepsilon\|F\|_\infty\|G\|_\infty .
\]
Letting $\Delta\uparrow\Z\times\mathbb R$ by Lemma~\ref{lem:critical-loops}(1), and using $\mathbb E_\phi[G\circ\tau]=\mathbb E_\phi[G]$, gives $|\mathbb E_\phi[F\,G\circ\tau]-\mathbb E_\phi[F]\mathbb E_\phi[G]|\le\varepsilon\|F\|_\infty\|G\|_\infty$ for all large $|n|$. Bounded measurable functions are $L^1(\phi)$-limits of such local functions, and translation invariance keeps the approximation uniform in $n$. The map $\omega\mapsto\Gamma$ commutes with these translations, so mixing passes to $\nu$.
\end{proof}

\begin{remark}[Why there is mixing]
\label{rem:mixing}
By the domain Markov property, the configuration far from a box $B$ influences the configuration in $B$ only through the boundary condition it induces on a large box around $B$. By comparison of boundary conditions, every such measure lies between the free and the wired ones, and at criticality these have the same limit, $\phi^0=\phi^1$ \cite[Theorem~1.5]{DCLM}, so the boundary condition is forgotten. Geometrically, the dual and primal circuits of Lemma~\ref{lem:separating} around $B$, present at every scale with probability bounded below, act as free and wired boundary conditions that screen $B$ from the far region. The screening comes scale by scale, and the argument gives no exponential rate.
\end{remark}

\subsection{The thermodynamic limit}

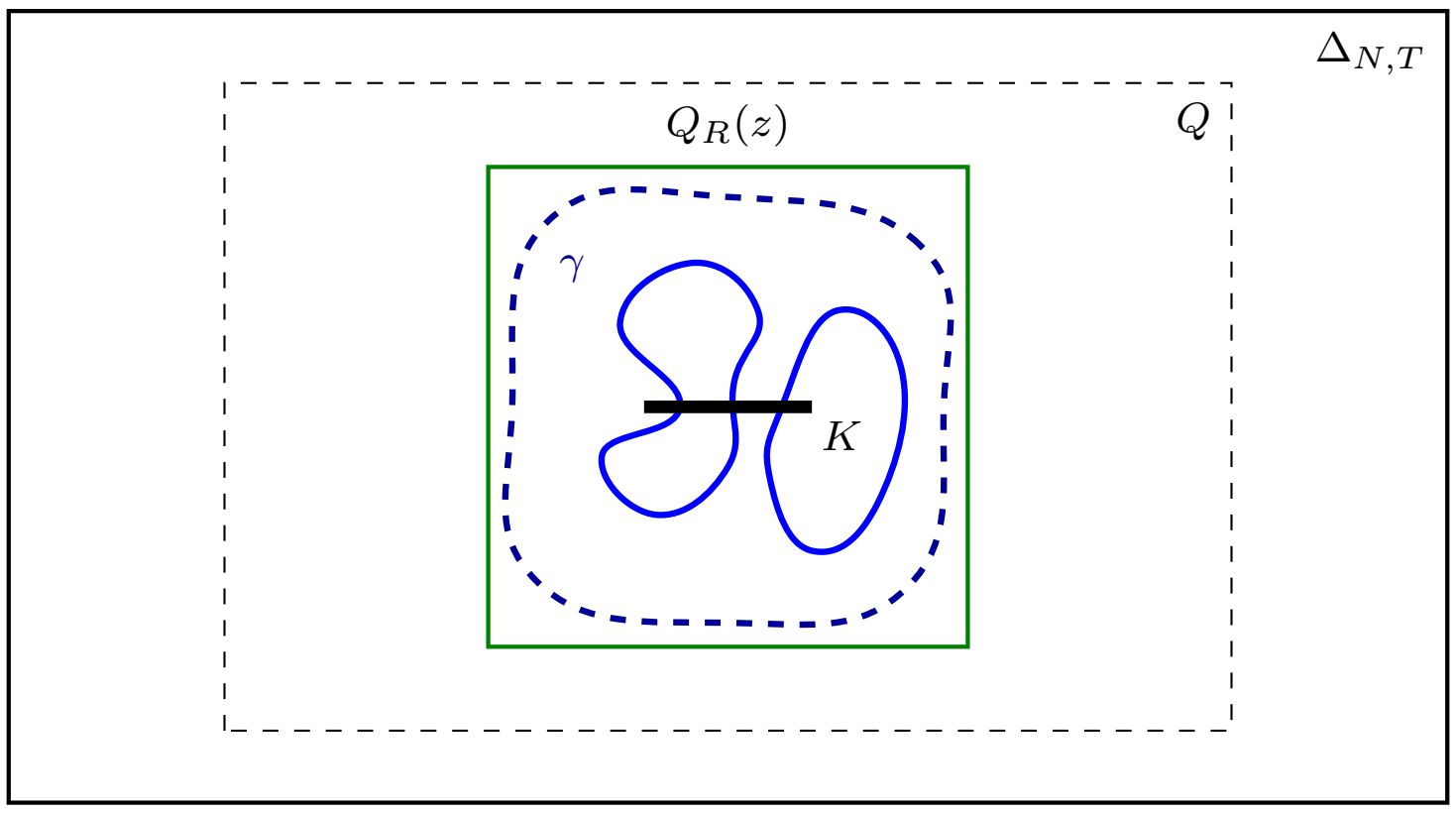
\begin{figure}[ht]
\centering
\begin{tikzpicture}[scale=0.8]
 \draw[thick] (-6,-3.3) rectangle (6,3.3);
 \node[below left] at (6,3.3) {\scriptsize $\Delta_{N,T}$};
 \draw[dashed] (-4.2,-2.7) rectangle (4.2,2.7);
 \node[below left] at (4.2,2.7) {\scriptsize $Q$};
 \draw[thick,green!50!black] (-2,-2) rectangle (2,2);
 \node[above] at (0,2) {\scriptsize $Q_R(z)$};
 \draw[very thick,blue!60!black,dashed] plot[smooth cycle,tension=0.8] coordinates {(1.8,0) (1.6,1.35) (0,1.75) (-1.5,1.55) (-1.8,0) (-1.6,-1.45) (0,-1.8) (1.5,-1.5)};
 \node[blue!60!black] at (-1.3,1.15) {\scriptsize $\gamma$};
 \draw[very thick,blue] plot[smooth cycle,tension=0.8] coordinates {(-0.4,0) (-0.9,0.7) (-0.3,1.2) (0.25,0.8) (0.05,0.2) (0.0,-0.5) (-0.6,-0.9) (-1.05,-0.4)};
 \draw[very thick,blue] plot[smooth cycle,tension=0.8] coordinates {(0.45,0) (0.9,0.8) (1.45,0.3) (1.25,-0.8) (0.7,-1.2) (0.35,-0.6)};
 \draw[line width=2.4pt] (-0.7,0) -- (0.7,0);
 \node[below] at (0.95,0.05) {\scriptsize $K$};
\end{tikzpicture}
\caption{The truncation in the proof of Proposition~\ref{prop:ground}. The word $W$ is inserted on $K=I\times\{0\}$. On $\mathcal E_R$ the loops meeting $K$ (solid) lie in $Q_R(z)$, so $f_W=f^{(R)}_W$ is determined by the configuration in $Q_R(z)$. A loop $\gamma$ (dashed) surrounding $K$ inside $Q_R(z)$ confines them, and by Lemma~\ref{lem:looptail} $\mathcal E_R$ fails with probability at most $C_KR^{-\beta_0}$. Once $\Delta_{N,T}$ contains the finite volume $Q$, Lemma~\ref{lem:critical-loops}(1) makes the laws of the configuration in $Q_R(z)$ under $\nu_{N,T}$ and $\nu$ $\varepsilon$-close.}
\label{fig:truncation}
\end{figure}

\begin{proposition}[The state $\omega_k$]
\label{prop:ground}
The limit $\omega_k(a)=\lim_{N\to\infty}\omega_N(a)$ exists for every local $a$. For a word $W$ it is given by $\omega_k(W)=\mathbb E_\nu[f_W]$, and $\omega_{N,T}(W)\to\mathbb E_\nu[f_W]$ as $N$ and $T$ tend to infinity jointly. The state $\omega_k$ is invariant under $\varsigma^2$ and is a ground state for the interaction $-\sum_ie_i$.
\end{proposition}
\begin{proof}
Let $W=e_{i_1}\cdots e_{i_r}$ be a word of length $r$ (Section~\ref{sec:loops}) supported in an interval $I$, and let $K=I\times\{0\}$ be the set of strand positions at which it is inserted. Choose $z$ and $R_0$ as in Lemma~\ref{lem:looptail}. Let $\mathcal E_R$ be the event that every loop meeting $K$ is contained in $Q_R(z)$. On $\mathcal E_R$, $\ell(\Gamma;W)-\ell(\Gamma)$ is the change in the number of loops through $K$, and it is determined by how these loops connect the strands at $K$. Each loop is traced by following the configuration along it, so on $\mathcal E_R$ these loops are determined by the configuration in $Q_R(z)$ (Figure~\ref{fig:truncation}). Hence there is a function $f^{(R)}_W$ of the configuration in $Q_R(z)$, with values in $[\delta^{-2r},1]$, which equals $f_W$ on $\mathcal E_R$. It is computed by following the loops through $K$ inside $Q_R(z)$, and we set $f^{(R)}_W=1$ if one of them reaches the boundary of $Q_R(z)$. Both $f_W$ and $f^{(R)}_W$ take values in $[\delta^{-2r},1]$, so $|f_W-f^{(R)}_W|\le1$, and $f_W-f^{(R)}_W$ vanishes on $\mathcal E_R$. The complement of $\mathcal E_R$ is the event in \eqref{eq:looptail} that some loop meeting $K$ leaves $Q_R(z)$. Hence $|\mathbb E_{\mathbb P}[f_W]-\mathbb E_{\mathbb P}[f^{(R)}_W]|\le\mathbb P(\mathcal E_R^c)$, and \eqref{eq:looptail} gives
\[
 \bigl|\mathbb E_{\mathbb P}[f_W]-\mathbb E_{\mathbb P}[f^{(R)}_W]\bigr|\le C_KR^{-\beta_0}\qquad(\mathbb P=\nu\text{ or }\nu_{N,T},\ Q_R(z)\subseteq\Delta_{N,T}).
\]
Given $\varepsilon>0$, choose $R$ with $C_KR^{-\beta_0}<\varepsilon$, and then, by Lemma~\ref{lem:critical-loops}(1) with $B=Q_R(z)$, a finite volume $Q\supseteq Q_R(z)$ such that $|\mathbb E_{\nu_{N,T}}[f^{(R)}_W]-\mathbb E_\nu[f^{(R)}_W]|<\varepsilon$ whenever $\Delta_{N,T}\supseteq Q$. With \eqref{eq:estimator} this gives $|\omega_{N,T}(W)-\mathbb E_\nu[f_W]|<3\varepsilon$ for all such $N,T$. Letting $T\to\infty$ at fixed $N$ gives $|\omega_N(W)-\mathbb E_\nu[f_W]|\le3\varepsilon$ for all large $N$. So $\omega_N(W)\to\mathbb E_\nu[f_W]$, and by linearity $\omega_N(a)$ converges for every local $a$. Positivity and normalization pass to the limit, and $\varsigma^2$-invariance follows from the invariance of $\nu$ under translation by one primal column.

A local $a$ preserves $V_N$, so
\[
 \omega_N(a^*[h_N,a])=\langle a\Omega_N,(h_N-\epsilon_N)a\Omega_N\rangle\ge0.
\]
For large $N$, $[h_N,a]$ does not depend on $N$, and the limit is the ground state inequality.
\end{proof}

The state is selected by the cup boundary vector and the vacuum block. No classification of infinite volume ground states is claimed for $k\ge3$. For $k=2$ the ground state is unique by \cite{ArakiMatsui85}.

\section{Purity}
\label{sec:purity}

This section shows that $\omega_k$ is pure (Proposition~\ref{prop:purity}), by proving imaginary time clustering from the loop representation and applying the Borchers--Arveson theorem. Let $(\pi,\HH,\Omega)$ be the GNS representation of $\omega_k$. The dynamics of the finite-range interaction, $\alpha_t(a)=\lim_Ne^{ith_N}ae^{-ith_N}$ for local $a$, exists as a norm limit \cite[Section~6.2]{BR2}, the proof using only that elements of $B_k$ localized in disjoint intervals commute. Since $\omega_k$ is a ground state, $\alpha_t$ is implemented in $\pi$ by $e^{it\mathcal L}$ with $\mathcal L\ge0$ and $\mathcal L\Omega=0$ \cite[Section~5.3]{BR2}.

\begin{lemma}[Imaginary time clustering]
\label{lem:clustering}
For local $a,b\in B_k$,
\[
 \lim_{s\to\infty}\langle\pi(a)\Omega,e^{-s\mathcal L}\pi(b)\Omega\rangle=\omega_k(a^*)\,\omega_k(b).
\]
\end{lemma}
\begin{proof}
By linearity we may take $a$ and $b$ to be words, of lengths $r_a$ and $r_b$. Put $G(s)=\langle\pi(a)\Omega,e^{-s\mathcal L}\pi(b)\Omega\rangle$ for $s\ge0$. The proof has three steps.

\emph{Step 1, finite volume.} The two-insertion Schwinger function on the strip $I_N\times[-T,T+s]$, with $b$ inserted at time $0$ and $a^*$ at time $s$, is
\[
 G_{N,T}(s)=\frac{\langle\xi_N,e^{-Th_N}a^*e^{-sh_N}be^{-Th_N}\xi_N\rangle}{\langle\xi_N,e^{-(2T+s)h_N}\xi_N\rangle}.
\]
By \eqref{eq:PF}, since $\langle\Omega_N,\xi_N\rangle\ne0$,
\[
 \lim_{T\to\infty}G_{N,T}(s)=\langle a\Omega_N,e^{-s(h_N-\epsilon_N)}b\Omega_N\rangle=:G_N(s).
\]

\emph{Step 2, $G_N(s)\to G(s)$.} By polarization it suffices to take $a=b$. Let $m_{N,a}$ be the spectral measure of the operator $h_N-\epsilon_N$ on $V_N$ in the vector $a\Omega_N$, a positive measure on $[0,\infty)$ of total mass $\omega_N(a^*a)$. Its Fourier transform is
\[
 \int e^{it\lambda}\,dm_{N,a}(\lambda)=\langle a\Omega_N,e^{it(h_N-\epsilon_N)}a\Omega_N\rangle=\omega_N\bigl(a^*\alpha^{(N)}_t(a)\bigr),\qquad\alpha^{(N)}_t(a)=e^{ith_N}ae^{-ith_N},
\]
because $h_N\Omega_N=\epsilon_N\Omega_N$. As $N\to\infty$, $\alpha^{(N)}_t(a)\to\alpha_t(a)$ in norm and $\omega_N\to\omega_k$ on local elements (Proposition~\ref{prop:ground}), so the Fourier transform converges to
\[
 \omega_k(a^*\alpha_t(a))=\langle\pi(a)\Omega,e^{it\mathcal L}\pi(a)\Omega\rangle=\int e^{it\lambda}\,dm_a(\lambda),
\]
where $m_a$ is the spectral measure of $\mathcal L$ in $\pi(a)\Omega$, a positive measure on $[0,\infty)$ because $\mathcal L\ge0$. The total masses converge as well, so by L\'evy's continuity theorem $m_{N,a}\to m_a$ weakly. The function $\lambda\mapsto e^{-s\lambda}$ is bounded and continuous on $[0,\infty)$, which carries all these measures, so $G_N(s)=\int e^{-s\lambda}dm_{N,a}\to\int e^{-s\lambda}dm_a=G(s)$.

\emph{Step 3, the loop representation and $s\to\infty$.} As in \eqref{eq:estimator},
\[
 G_{N,T}(s)=\mathbb E\bigl[f^{(2)}_s\bigr],\qquad f^{(2)}_s(\Gamma)=\delta^{\ell(\Gamma;b,a^*)-\ell(\Gamma)-r_a-r_b}\in[\delta^{-2(r_a+r_b)},1],
\]
the expectation being under the loop law of the strip $I_N\times[-T,T+s]$, and $\ell(\Gamma;b,a^*)$ counting the loops after inserting $b$ at time $0$ and $a^*$ at time $s$. The truncation argument in the proof of Proposition~\ref{prop:ground}, with the two insertion sets $K_0$ at time $0$ and $K_s$ at time $s$, shows that $G_{N,T}(s)\to\mathbb E_\nu[f^{(2)}_s]$ as $N$ and $T$ tend to infinity jointly, and in particular $\lim_N\lim_{T}G_{N,T}(s)=\mathbb E_\nu[f^{(2)}_s]$. With Steps 1 and 2, $G(s)=\mathbb E_\nu[f^{(2)}_s]$.

Now let $E_s$ be the event that every loop meeting $K_0$ is contained in $Q_{s/3}(z_0)$ and every loop meeting $K_s$ is contained in $Q_{s/3}(z_s)$, where $z_0$ and $z_s=z_0+(0,s)$ are centres for $K_0$ and $K_s$. On $E_s$ no loop meets both $K_0$ and $K_s$. The surgery at time $0$ changes only loops through $K_0$, and the new loops it creates consist of arcs of these loops and arcs near $K_0$, so they are disjoint from $K_s$. Hence the surgery at time $s$ acts on the same loops before and after the surgery at time $0$. The changes in the number of loops add, and
\[
 f^{(2)}_s=f_b\cdot(f_{a^*}\circ\tau_s)\qquad\text{on }E_s,
\]
where $f_{a^*}\circ\tau_s$ is the one-insertion estimator of $a^*$ inserted at time $s$. Let $K'$ be the insertion set of $a^*$ at time $0$, so that $K_s$ is its translate by $s$. By \eqref{eq:looptail} and time invariance of $\nu$, $\nu(E_s^c)\le(C_{K_0}+C_{K'})(s/3)^{-\beta_0}$. All estimators take values in $[0,1]$, so
\[
 \bigl|G(s)-\mathbb E_\nu[f_b\cdot(f_{a^*}\circ\tau_s)]\bigr|\le2\nu(E_s^c)\to0 .
\]
By time mixing (Lemma~\ref{lem:looptail}), $\mathbb E_\nu[f_b\cdot(f_{a^*}\circ\tau_s)]\to\mathbb E_\nu[f_b]\,\mathbb E_\nu[f_{a^*}]=\omega_k(b)\,\omega_k(a^*)$, using Proposition~\ref{prop:ground} for the last equality.
\end{proof}

\begin{proposition}
\label{prop:purity}
$\omega_k$ is pure.
\end{proposition}
\begin{proof}
As $s\to\infty$, $e^{-s\mathcal L}$ converges strongly to the projection $P_0$ onto $\ker\mathcal L$. On the dense set $\pi(B_{k,\mathrm{loc}})\Omega$, Lemma~\ref{lem:clustering} gives $P_0=|\Omega\rangle\langle\Omega|$. Since $\Omega\in\ker\mathcal L$ is cyclic for $\pi(B_k)$, the ground state subspace $\ker\mathcal L$ is generating, so $\pi$ with the implementation $e^{it\mathcal L}$ is a ground state representation in the sense of \cite{NR24}. By the Borchers--Arveson theorem in the form \cite[Proposition~2.7]{NR24} the positive implementation is then minimal, and in particular $e^{it\mathcal L}\in\pi(B_k)''$, see also \cite[Section~5.3]{BR2}. Every $x\in\pi(B_k)'$ therefore commutes with $P_0$, so $x\Omega=c\Omega$, and cyclicity gives $x=c\one$.
\end{proof}

\section{Affine modules and sewing}
\label{sec:affine}

This section and the next prove the remaining assertion $H_{\omega_k}\cong\one$ of Theorem~\ref{thm:TLJ}. Here we provide the two ingredients that turn charged coefficients into quantities controlled by loops. Section~\ref{sec:affine}.1 recalls the Hilbert affine Temperley--Lieb module $\mathscr V_Z$ of a centre object $Z$ \cite{DGG} and proves the affine estimate (Lemma~\ref{lem:affine-estimate}). A plain affine diagram acts, relative to its number of loops, with norm at most $\delta^{(m+n)/2}\gamma_Z^{\lfloor r/2\rfloor}$, where $r$ is the number of essential circles and $\gamma_Z=(\mu_Z/\delta)^2<1$ for $Z\not\cong\one$ by the explicit values of $\mu_Z$ (Proposition~\ref{prop:muZ}). Section~\ref{sec:affine}.2 defines the charged coefficients $\langle\xi\otimes\Omega,U_Z^n(\eta\otimes\Omega)\rangle$ and proves the sewing lemma (Lemma~\ref{lem:sewing}). For each loop configuration the charged evaluation is a matrix coefficient of an affine diagram whose essential circles are the loops separating the two insertions. Together these give the bound \eqref{eq:F-bound} on the charged estimator. Section~\ref{sec:decay} combines it with the separating loops of Lemma~\ref{lem:separating} to prove decay (Proposition~\ref{prop:decay}), and deduces $H_{\omega_k}\cong\one$ by rigidity (Lemma~\ref{lem:rigidity}).

\subsection{The affine representation of a centre object}

\paragraph{Affine diagrams.} An affine diagram $D:(m,\epsilon)\to(n,\eta)$ is drawn in the annulus between an inner circle with $2m$ marked points and an outer circle with $2n$ marked points. It consists of disjoint strings, each either joining two marked points or closed, and finitely many disjoint boxes, each labelled by a morphism of $\C_k$ between tensor powers of $X$, whose legs are ends of strings. Strings never cross. The complementary regions are shaded alternately, and $\epsilon,\eta\in\{+,-\}$ record the shading next to a distinguished marked point on the inner and on the outer circle. Diagrams are taken up to isotopy fixing both circles pointwise, so a string that winds around the annulus is not isotopic to one that does not. We call $m$ the affine degree, and the $2m$ marked points carry the strands of $Y^{\otimes m}=X^{\otimes2m}$. A plain affine Temperley--Lieb diagram is an affine diagram without boxes. It has three kinds of strings. Through strings join the inner circle to the outer circle. Arcs join two points of the same circle. The remaining strings are closed loops. A closed loop either bounds a disc in the annulus, and is then called contractible, or winds once around the inner circle, and is then called essential (Figure~\ref{fig:affine}).

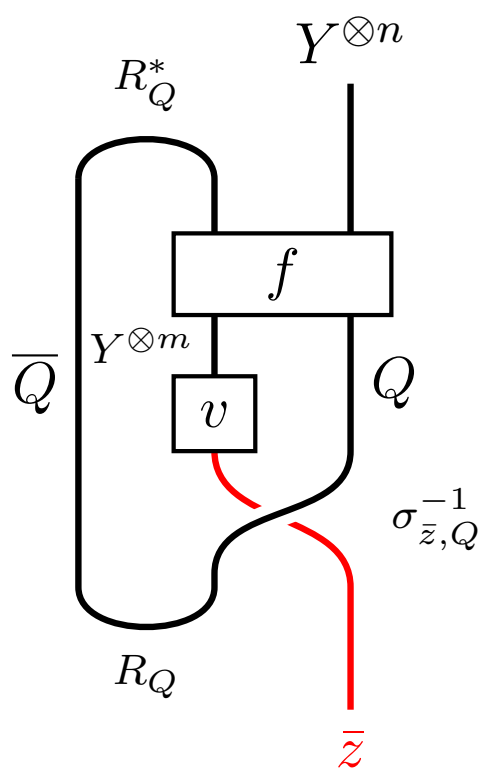
\begin{figure}[ht]
\centering
\begin{tikzpicture}[scale=0.9]
 \draw[very thick,red] (2,0) -- (2,0.9) to[out=90,in=-90] (1,1.9);
 \draw[very thick,white,line width=5pt] (1,1.0) to[out=90,in=-90] (2,1.9);
 \draw[very thick] (0,3.9) -- (0,0.9) to[out=-90,in=-90] (1,0.9) -- (1,1.0) to[out=90,in=-90] (2,1.9) -- (2,2.9);
 \draw[very thick] (1,2.45) -- (1,2.9);
 \draw[very thick] (1,3.5) -- (1,3.9) to[out=90,in=90] (0,3.9);
 \draw[very thick] (2,3.5) -- (2,4.6);
 \draw[thick,fill=white] (0.7,1.9) rectangle (1.3,2.45);
 \node at (1,2.17) {$v$};
 \draw[thick,fill=white] (0.7,2.9) rectangle (2.3,3.5);
 \node at (1.5,3.2) {$f$};
 \node[below,red] at (2,0) {$\bar z$};
 \node[left] at (0,2.4) {$\overline Q$};
 \node[right] at (2,2.4) {$Q$};
 \node[left] at (1,2.68) {\scriptsize $Y^{\otimes m}$};
 \node[above] at (2,4.6) {$Y^{\otimes n}$};
 \node[below] at (0.5,0.55) {\scriptsize $R_Q$};
 \node[above] at (0.5,4.25) {\scriptsize $R_Q^*$};
 \node[right] at (2.15,1.4) {\scriptsize $\sigma^{-1}_{\bar z,Q}$};
\end{tikzpicture}
\caption{The action \eqref{eq:affineaction}, read from bottom to top. The cup $R_Q$ creates $\overline QQ$, the strand $Q$ crosses $\bar z$ by the inverse half-braiding, $v:\bar z\to Y^{\otimes m}$ is applied, then the cut diagram $f:Y^{\otimes m}\otimes Q\to Q\otimes Y^{\otimes n}$, and the cap $R_Q^*$ closes $Q$ with $\overline Q$. Closing $Q$ around this way is the wrapping of the cut diagram around the annulus.}
\label{fig:affineaction}
\end{figure}

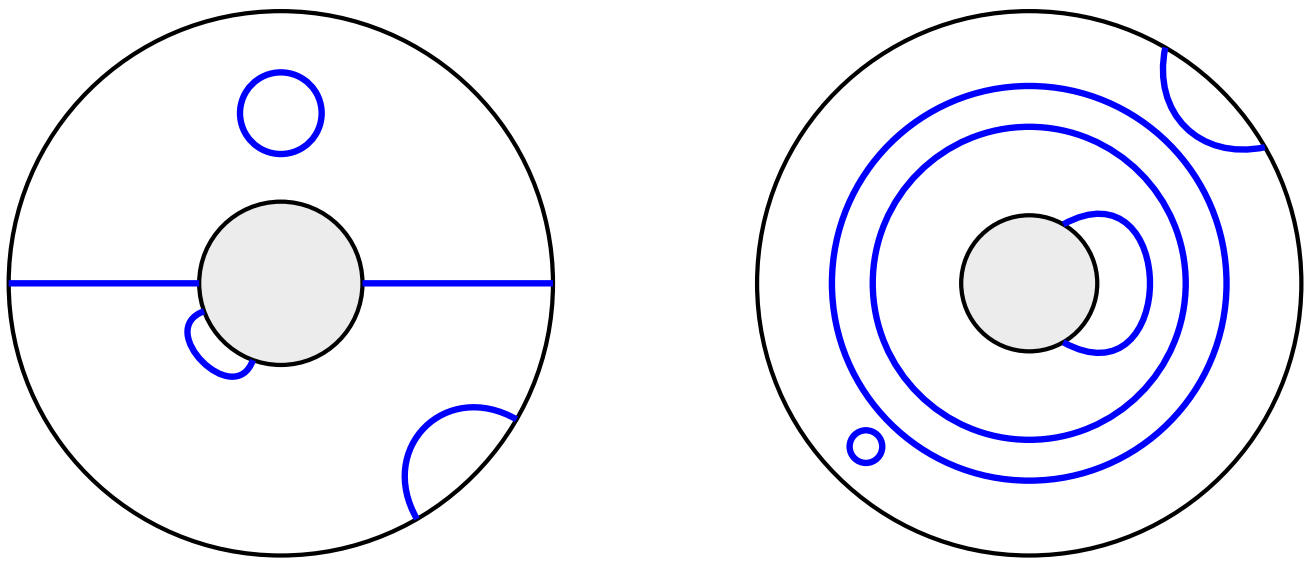
\begin{figure}[ht]
\centering
\begin{tikzpicture}[scale=0.9]
\begin{scope}
 \draw[thick] (0,0) circle (2);
 \draw[thick,fill=gray!15] (0,0) circle (0.6);
 \draw[very thick,blue] (0.6,0) -- (2,0);
 \draw[very thick,blue] (-0.6,0) -- (-2,0);
 \draw[very thick,blue] (0,1.25) circle (0.3);
 \draw[very thick,blue] (1,-1.732) .. controls (0.7,-1.21) and (1.21,-0.7) .. (1.732,-1);
 \draw[very thick,blue] (-0.564,-0.205) .. controls (-0.94,-0.34) and (-0.34,-0.94) .. (-0.205,-0.564);
 \node at (0,-2.5) {(a) through strings};
\end{scope}
\begin{scope}[xshift=5.5cm]
 \draw[thick] (0,0) circle (2);
 \draw[thick,fill=gray!15] (0,0) circle (0.5);
 \draw[very thick,blue] (0,0) circle (1.15);
 \draw[very thick,blue] (0,0) circle (1.45);
 \draw[very thick,blue] (0.25,0.433) .. controls (1.1,0.9) and (1.1,-0.9) .. (0.25,-0.433);
 \draw[very thick,blue] (1.732,1) .. controls (1.25,0.9) and (0.9,1.25) .. (1,1.732);
 \draw[very thick,blue] (-1.2,-1.2) circle (0.12);
 \node at (0,-2.5) {(b) essential circles};
\end{scope}
\end{tikzpicture}
\caption{Plain affine Temperley--Lieb diagrams in the annulus between the two holes. All diagrams are planar, so strings never cross. In (a) the two radial strings are through strings, there is an arc on each circle, and the closed loop is contractible. In (b) there are no through strings. The arcs on the inner and outer circles cut off discs, and the diagram factors through a zero-boundary space. Each of the two essential circles contributes a single-circle map of norm $\mu_Z$. The small loop in (b) is contractible.}
\label{fig:affine}
\end{figure}

For $Z\in\mathcal Z(\C_k)$ write $z$ for its underlying object, with half-braiding $\sigma_{z,Q}:z\otimes Q\to Q\otimes z$. Let $\mathscr V_Z$ be the Hilbert affine module of the dual object $\overline Z$, with
\begin{equation}
 \mathscr V_Z(m,+)=\Hom_{\C_k}(\bar z,Y^{\otimes m}),\qquad\langle v,w\rangle=\mathrm{Tr}_{\bar z}(v^*w).
 \label{eq:affinespace}
\end{equation}
Equation \eqref{eq:affinespace} defines the spaces $\mathscr V_Z(m,+)$, for boundaries whose distinguished region is unshaded. Boundaries whose distinguished region is shaded carry a second family of spaces $\mathscr V_Z(m,-)$, defined in \cite{DGG} through the dual shading. Together they form the shaded affine module $\mathscr V_Z$ of \cite{DGG}, in which an affine diagram $D:(m,\epsilon)\to(n,\eta)$ maps $\mathscr V_Z(m,\epsilon)$ to $\mathscr V_Z(n,\eta)$. We use the $-$-spaces only for $m=0$, through the zero-boundary space $\mathscr V_Z(0,-)$. Cutting an affine diagram $D:(m,+)\to(n,+)$ along a seam, a radial arc from the inner to the outer circle that avoids the boxes and meets the strings transversally in an even number $2r$ of points, turns the annulus into a rectangle and $D$ into a morphism $f:Y^{\otimes m}\otimes Q\to Q\otimes Y^{\otimes n}$ of $\C_k$ with $Q=Y^{\otimes r}$. With a standard coevaluation $R_Q$, $D$ acts by (Figure~\ref{fig:affineaction})
\begin{equation}
 \rho_Z(D)v=(R_Q^*\otimes\id_{Y^{\otimes n}})(\id_{\overline Q}\otimes f)(\id_{\overline Q}\otimes v\otimes\id_Q)(\id_{\overline Q}\otimes\sigma^{-1}_{\bar z,Q})(R_Q\otimes\id_{\bar z}).
 \label{eq:affineaction}
\end{equation}
Naturality of the half-braiding gives independence of the cut, and its tensor compatibility gives composition. The conjugate equations give cancellation of cups and caps. Unitarity gives $\langle v,\rho_Z(D)w\rangle=\langle\rho_Z(D^*)v,w\rangle$ for the inner products \eqref{eq:affinespace} and their counterparts on the $-$-spaces, where $D^*$ is the reflection of $D$ exchanging the two circles.

 Only $Q\in\C_k$ occurs. The full shaded Hilbert module, including the negative zero space, is that of Das--Ghosh--Gupta, whose main theorem identifies the centre of the even part with the locally finite Hilbert affine representations of finite support \cite[Theorem~5.13]{DGG}. In finite depth the support condition is automatic.

\paragraph{Reduction to plain diagrams.} The category $\C_k$ is the even part of the Temperley--Lieb--Jones category, whose morphism spaces between tensor powers of $X$ are quotients of the spaces spanned by Temperley--Lieb diagrams, the noncrossing pairings of the boundary points \cite{GHJ}. So every box label of an affine diagram $D$ is a linear combination of Temperley--Lieb diagrams. Substituting these expansions for the boxes writes $D$ as a linear combination of plain affine Temperley--Lieb diagrams $D'$, and since \eqref{eq:affineaction} is linear in $f$, $\rho_Z(D)$ is the same linear combination of the operators $\rho_Z(D')$. The module $\mathscr V_Z$ is therefore determined by the action of plain affine Temperley--Lieb diagrams.

\paragraph{Zero spaces and single-circle maps.} The zero spaces are $\mathscr V_Z(0,+)=\Hom_{\C_k}(\bar z,\one)$ and $\mathscr V_Z(0,-)$. The single-circle maps are the two operators $C_+:\mathscr V_Z(0,+)\to\mathscr V_Z(0,-)$ and $C_-:\mathscr V_Z(0,-)\to\mathscr V_Z(0,+)$ given by the plain affine diagrams with no marked points and one essential circle, across which the shading changes. The reflection exchanging the two circles carries each of these diagrams to the other, so $C_-=C_+^*$. The composite $C_-C_+$ is given by two parallel essential $X$-circles, that is an essential $Y$-circle, acting on $\Hom_{\C_k}(\bar z,\one)$ through the half-braiding. If $\one\not\subseteq z$, then $\mathscr V_Z(0,+)=0$. Each single-circle map has $\mathscr V_Z(0,+)$ as its domain or its target, so both vanish, and we put $\mu_Z=0$. If $\one\subseteq z$, Proposition~\ref{prop:muZ} below shows that $\mathscr V_Z(0,+)$ is one-dimensional and that the essential $Y$-circle acts on it by a scalar $\mu_Z^2$ with $\mu_Z\ge0$. Then
\[
 \|C_+\|=\|C_-\|=\|C_-C_+\|^{1/2}=\mu_Z .
\]
The number $\mu_Z$ is the parameter of the Hilbert module $\mathscr V_Z$ in the classification of Jones and Reznikoff \cite{JR06}. We do not use the classification. The values of $\mu_Z$ are computed from the $S$-matrix of $\mathrm{SU}(2)_k$.

\begin{proposition}[Values of $\mu_Z$]
\label{prop:muZ}
Let $Z\in\mathcal Z(\C_k)$ be simple with $\one\subseteq z$. Then $\Hom_{\C_k}(\bar z,\one)$ is one-dimensional, and the essential $Y$-circle acts on it by $\mu_Z^2$, where $\mu_Z=|2\cos\frac{(a+1)\pi}{k+2}|$ for a label $a$ determined by $Z$ as follows, or $\mu_Z=0$.
\begin{enumerate}[leftmargin=*]
\item For odd $k$, $Z=a\boxtimes\bar a$ with $a\in\{0,2,\dots,k-1\}$.
\item For even $k=2n$, either $Z$ is the class $[a,a]$ with $0\le a\le k$ and $a\ne n$, or $Z=(n,n)_+$ and $\mu_Z=0$.
\end{enumerate}
In particular $0\le\mu_Z\le\delta$ for every simple $Z$, and
\begin{equation}
 0\le\mu_Z<\delta\qquad\text{for every simple }Z\not\cong\one .
 \label{eq:muZstrict}
\end{equation}
\end{proposition}

\begin{proof}
For odd $k$ the simple current of $\mathrm{SU}(2)_k$ has odd spin, so the even part $\C_k$, with simple objects the even labels $a\in\{0,2,\dots,k-1\}$, is modular and $\mathcal Z(\C_k)\simeq\C_k\boxtimes\C_k^{\mathrm{rev}}$. This factorisation is \cite[Corollary~3.10]{BEK01}, with the trivial inclusion, see also \cite{Mue03}. The subfactor setting of \cite{BEK01} is no restriction, since a unitary fusion category always acts on the hyperfinite II$_1$ factor by Popa, see the discussion in \cite{CHPJP22} and \cite[Section~1.3]{EvansJones2025}. The simple objects are $Z=a\boxtimes\bar b$, with underlying object $z=a\otimes b$. All labels are self-dual, so $\one\subseteq z$ only for $a=b$, and then $\Hom_{\C_k}(\bar z,\one)$ is one-dimensional. A $c$-loop around the charge passes over the strand $a$ and under the strand $\bar a$. It forms a Hopf link with $a$ and acts by $S_{ca}/S_{0a}$, where $S_{ab}=\sqrt{2/(k+2)}\,\sin\frac{(a+1)(b+1)\pi}{k+2}$. Since $Y=X\otimes X\cong0\oplus2$ and $S_{2a}/S_{0a}=(S_{1a}/S_{0a})^2-1$, the essential $Y$-circle acts by $(S_{1a}/S_{0a})^2$, so
\[
 \mu_{a\boxtimes\bar a}=\Bigl|\frac{S_{1a}}{S_{0a}}\Bigr|=\Bigl|2\cos\frac{(a+1)\pi}{k+2}\Bigr| .
\]

For even $k=2n$ the simple current $J=k$ lies in $\C_k$ and braids trivially with all of $\C_k$. So $\C_k$ is degenerate and $\mathcal Z(\C_k)$ does not factor. Evans and Kawahigashi describe it through the asymptotic inclusion \cite[Section~4]{EK98b}. Its simple objects are the classes $[a,b]$ of pairs of labels $0\le a,b\le k$ with $a+b$ even, modulo $(a,b)\sim(k-a,k-b)$. The only fixed pair is $(n,n)$, and it gives two objects $(n,n)_+$ and $(n,n)_-$. There are $k^2/4+k/2+2$ simple objects \cite[Corollary~4.3]{EK98b}. The underlying object of $[a,b]$ is $a\otimes b$. By \cite[Theorem~4.2]{EK98b} the underlying objects of $(n,n)_+$ and $(n,n)_-$ are the sums of the labels in $n\otimes n$ congruent to $0$ and to $2$ modulo $4$. In $\mathrm{SU}(2)_k$ the product $a\otimes b$ contains $0$ only for $a=b$, and then once. So $\one\subseteq z$ only for $[a,a]$ and for $(n,n)_+$, and then $\Hom_{\C_k}(\bar z,\one)$ is one-dimensional. For a class $[a,b]$ with $(a,b)\ne(n,n)$ the corresponding centre object is obtained from the pair $(a,b)$ in the quantum double construction of Evans and Kawahigashi \cite[Section~4]{EK98b}. In the Longo--Rehren description the two chiral factors carry the braiding and the opposite braiding \cite[Sections~2.2 and~4]{BEK01}. Equivalently, the half-braiding of $a\otimes b$ with an object $Q$ passes $Q$ over the $a$-strand and under the $b$-strand,
\[
 \sigma_{a\otimes b,Q}=\bigl(\mathsf c^{-1}_{Q,a}\otimes\id_b\bigr)\circ\bigl(\id_a\otimes\mathsf c_{b,Q}\bigr),
\]
where $\mathsf c$ is the braiding of $\mathrm{SU}(2)_k$. The convention for the order of the tensor factors and for which crossing is drawn as over plays no role. What matters is that the two factors use opposite braidings. We do not use \cite[Corollary~4.8]{BEK01} here, since it assumes a nondegenerate braiding, and $\C_k$ is degenerate for even $k$. The general half-braiding and relative braiding of \cite[Sections~2.2 and~4]{BEK01} need no such assumption, and the identification of the centre objects with the classes $[a,b]$ is that of \cite[Section~4]{EK98b}. Thus for $[a,a]$ the half-braiding with the even object $Y$ passes $Y$ over one $a$-strand and under the other. So the computation for odd $k$ applies unchanged and gives $\mu_{[a,a]}=|2\cos\frac{(a+1)\pi}{k+2}|$, the same for $a$ and $k-a$. The zero space of $(n,n)_+$ is the zero space of $n\otimes n$, so $\mu_{(n,n)_+}=|2\cos\frac{\pi}{2}|=0$.

In both cases $|2\cos\frac{(a+1)\pi}{k+2}|\le2\cos\frac{\pi}{k+2}=\delta$, with equality only for $a\in\{0,k\}$. For odd $k$ the label $a$ is even and at most $k-1$, so equality means $a=0$, that is $Z=\one$. For even $k$ the classes $[0,0]$ and $[k,k]$ coincide and are $\one$. This proves \eqref{eq:muZstrict}.
\end{proof}

\begin{lemma}[Affine diagram estimate]
\label{lem:affine-estimate}
Let $D:(m,\epsilon)\to(n,\eta)$ be a plain affine Temperley--Lieb diagram. Let $c(D)$ be the number of its contractible closed loops and $r(D)$ the number of its essential closed loops, and let $\ell(D)$ be the number of loops after filling both holes with fixed noncrossing pairings. For $v\in\mathscr V_Z(n,\eta)$, $w\in\mathscr V_Z(m,\epsilon)$ and $\gamma_Z=(\mu_Z/\delta)^2$, we have
\begin{equation}
 \frac{|\langle v,\rho_Z(D)w\rangle|}{\delta^{\ell(D)}}\le\delta^{(m+n)/2}\|v\|\|w\|\,\gamma_Z^{\lfloor r(D)/2\rfloor}.
 \label{eq:affine-estimate}
\end{equation}
Here $0^0=1$.
\end{lemma}
\begin{proof}
Contractible circles contribute $\delta^{c(D)}$. An essential circle separates the inner from the outer circle, so through strings and essential circles do not occur together. If there are $2h>0$ through strings there are no essential circles, rotation is unitary, and each of the $m+n-2h$ boundary arcs has norm $\sqrt\delta$, so $\|\rho_Z(D)\|\le\delta^{c(D)+(m+n)/2}$. Suppose there are no through strings. An arc with both ends on the same circle cuts off a disc from the annulus, and essential circles avoid these discs. So, up to contractible circles, $D=D_{\mathrm{out}}\circ C_r\circ\dots\circ C_1\circ D_{\mathrm{in}}$, where $D_{\mathrm{in}}:(m,\epsilon)\to(0,\epsilon')$ consists of the arcs on the inner circle, $D_{\mathrm{out}}:(0,\eta')\to(n,\eta)$ of the arcs on the outer circle, and $C_1,\dots,C_r$, $r=r(D)$, are single-circle maps. Thus $D$ factors through the zero-boundary spaces, and $\|\rho_Z(D)\|\le\delta^{c(D)+(m+n)/2}\mu_Z^{r(D)}$. Since $\ell(D)\ge c(D)+r(D)$ and $\delta\ge1$,
\[
 \frac{\|\rho_Z(D)\|}{\delta^{\ell(D)}}\le\delta^{(m+n)/2}\Bigl(\frac{\mu_Z}{\delta}\Bigr)^{r(D)}\le\delta^{(m+n)/2}\gamma_Z^{\lfloor r(D)/2\rfloor},
\]
using $0\le\mu_Z/\delta\le1$ (Proposition~\ref{prop:muZ}). In the case with through strings the same computation, with $r(D)=0$, gives the bound with final factor $1$. The two cases are shown in Figure~\ref{fig:affine}.
\end{proof}

\begin{remark}[$k=2$]
\label{rem:k2affine}
At $k=2$ the four objects of $\mathcal Z(\Vect(\Z_2))$ have $z\in\{\one,\psi\}$. For $m$ and $f$, $z=\psi$, so $\mathscr V_Z(0,+)=0$ and $\mu_Z=0$. For $\one$ and $e$, $z=\one$ and the half-braiding with $Q$ is a character $\chi$ of the parity of $Q$, trivial for $\one$ and the sign for $e$. For $e$ this half-braiding corresponds to the half-line parity $\rho_e=\Lambda_-|_B$ of \eqref{eq:sector-representatives}. Two parallel essential $X$-circles form an essential $Y$-circle, which acts on $\mathscr V_Z(0,+)$ by $\sum_{a\subset Y}d_a\chi(a)=1+\chi(\psi)$. Thus $\mu_\one^2=2=\delta^2$ and $\mu_e=0$. For every nontrivial $Z$, therefore, both single-circle maps vanish. For $m$ and $f$ one of them has zero domain and the other zero target, and for $e$ they are zero because $\mu_e=0$. Now let $D$ be a plain affine diagram with $r(D)\ge1$. It has no through strings, and the proof of Lemma~\ref{lem:affine-estimate} writes $\rho_Z(D)$ as a composite containing $r(D)$ single-circle maps. So $\rho_Z(D)=0$. In the sewing of Lemma~\ref{lem:sewing} below, the essential circles of $D_\Gamma$ are the loops of $\Gamma$ separating the two insertions, so at $k=2$ a single separating loop annihilates the charged evaluation.
\end{remark}

\subsection{Charged coefficients}

For $Z\in\mathcal Z(\C_k)$ the local correspondence spaces are
\[
 \mathsf Y_Z(I)=\Hom_{\C_k}\bigl(Y^{\otimes|I|},Y^{\otimes|I|}\otimes z\bigr),\qquad
 a\xi b=(a\otimes\id_z)\xi b,\qquad\langle\xi,\eta\rangle_{B_k}=\xi^*\eta.
\]
Let $A=Y^{\otimes|I|}$ and $Q=Y^{\otimes r}$. Let $J$ be the interval of the $r$ sites immediately to the right of $I$, and $J'$ that of the $r$ sites immediately to the left. Adding $Q$ on the right embeds $\mathsf Y_Z(I)$ into $\mathsf Y_Z(I\cup J)$ by
\begin{equation}
 \iota^R_Q(\xi)=(\id_A\otimes\sigma_{z,Q})(\xi\otimes\id_Q),
 \label{eq:rightembed}
\end{equation}
and adding $Q$ on the left embeds $\mathsf Y_Z(I)$ into $\mathsf Y_Z(J'\cup I)$ by $\xi\mapsto\id_Q\otimes\xi$. The half-braiding identities and unitarity make these embeddings compatible and isometric, and completion gives the DHR bimodule $\mathsf Y_Z$ \cite[Section~6.2]{JNPW}. Relabelling an interval by its translate by one site defines a unitary $s_Z$ of $\mathsf Y_Z$ with $s_Z(a\xi b)=\varsigma^2(a)s_Z(\xi)\varsigma^2(b)$ and $\langle s_Z\xi,s_Z\eta\rangle=\varsigma^2(\langle\xi,\eta\rangle)$, which maps $\mathsf Y_Z(I)$ onto $\mathsf Y_Z(I+1)$. Let $U$ be the unitary on $\HH$ with $U\pi(a)\Omega=\pi(\varsigma^2(a))\Omega$, which exists because $\omega_k$ is $\varsigma^2$-invariant, and put
\[
 \mathcal K_Z=\mathsf Y_Z\boxtimes_{B_k}\HH,\qquad\pi_Z(b)(\xi\otimes v)=b\xi\otimes v,\qquad U_Z(\xi\otimes v)=s_Z(\xi)\otimes Uv.
\]
Then $U_Z$ is unitary, $U_Z\pi_Z(b)U_Z^*=\pi_Z(\varsigma^2(b))$, vectors $\xi\otimes\Omega$ with local $\xi$ are dense, and
\begin{equation}
 \langle\xi\otimes\Omega,U_Z^n(\eta\otimes\Omega)\rangle=\omega_k\bigl(\langle\xi,s_Z^n\eta\rangle_{B_k}\bigr).
 \label{eq:charged-coeff}
\end{equation}

For $A=Y^{\otimes p}$ the endpoint map (Figure~\ref{fig:endpoint}(a)) and its inverse are
\begin{align*}
 J_A(\xi)&=(\id_{\overline AA}\otimes\overline R_z^*)(\id_{\overline A}\otimes\xi\otimes\id_{\bar z})(R_A\otimes\id_{\bar z})\in\Hom_{\C_k}(\bar z,\overline AA)=\mathscr V_Z(2p,+),\\
 J_A^{-1}(v)&=(\overline R_A^*\otimes\id_{Az})(\id_A\otimes v\otimes\id_z)(\id_A\otimes R_z),
\end{align*}
using $\overline AA\cong Y^{\otimes2p}$. The conjugate equations give $\langle J_A\xi,J_A\eta\rangle=\mathrm{Tr}_A(\xi^*\eta)$ and $\|J_A\xi\|\le\delta^p\|\xi\|$.

\begin{figure}[ht]
\centering
\begin{tikzpicture}[scale=0.85]
 \draw[very thick,red] (3.15,1.45) .. controls (3.3,2.3) and (0.6,1.9) .. (0.45,2.55);
 \node[red] at (2.15,1.65) {\small $z$};
 \draw[very thick,white,line width=5pt] (1.6,0) -- (1.6,4);
 \draw[very thick] (1.6,0) -- (1.6,4);
 \draw[very thick,white,line width=5pt] (2.9,1.45) -- (2.9,4);
 \draw[very thick] (2.9,1.45) -- (2.9,4);
 \draw[very thick] (2.9,0) -- (2.9,0.85);
 \draw[very thick] (0.2,0) -- (0.2,2.55);
 \draw[very thick] (0.2,3.15) -- (0.2,4);
 \draw[thick,fill=white] (2.55,0.85) rectangle (3.4,1.45);
 \node at (2.97,1.15) {$\eta$};
 \draw[thick,fill=white] (-0.15,2.55) rectangle (0.75,3.15);
 \node at (0.3,2.85) {$\xi^*$};
 \node[below] at (0.2,0) {$A$};
 \node[below] at (1.6,0) {$M$};
 \node[below] at (2.9,0) {$B'$};
 \node[above] at (0.2,4) {$A$};
 \node[above] at (1.6,4) {$M$};
 \node[above] at (2.9,4) {$B'$};
 \node[right] at (3.5,2.1) {\small $\sigma^{-1}_{z,MB'}$};
 \draw[decorate,decoration={brace,mirror,amplitude=4pt}] (1.25,-0.55) -- (1.95,-0.55);
 \node[below] at (1.6,-0.7) {\small $g$ sites};
\end{tikzpicture}
\caption{The charged insertion $b_{\xi,\eta}$ of \eqref{eq:charged-insertion}, read from bottom to top. The strands $A$, $M$, $B'$ belong to the sites of $I$, of the gap, and of $J$. The box $\eta$ emits the charge line $z$, which passes under the strands of $B'$ and $M$ by the inverse half-braiding and is absorbed by $\xi^*$.}
\label{fig:gap}
\end{figure}
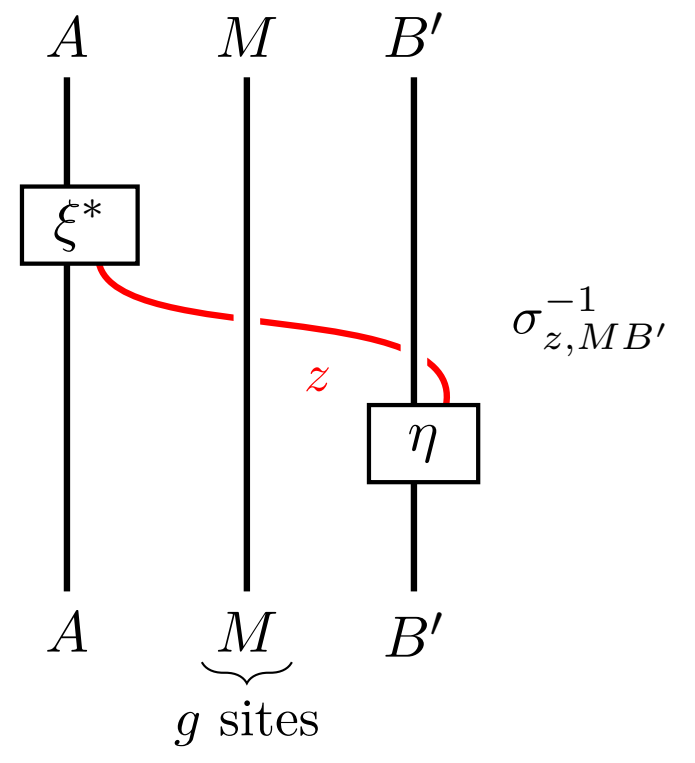

\paragraph{The charged insertion.} Let $\xi\in\mathsf Y_Z(I)$ and $\eta\in\mathsf Y_Z(J)$, where $I$ and $J$ are intervals of $p$ and $p'$ sites and $I$ lies to the left of $J$. Let $g\ge0$ be the number of sites strictly between them, and put
\[
 A=Y^{\otimes p},\qquad M=Y^{\otimes g},\qquad B'=Y^{\otimes p'},
\]
so that $AMB'=Y^{\otimes|K|}$ for the interval $K$ from the left end of $I$ to the right end of $J$. The gap word $M$ consists of the $Y$-strands of the sites between $I$ and $J$. To compare $\xi$ and $\eta$ in $\mathsf Y_Z(K)$, embed $\xi$ by adding $MB'$ on the right, by \eqref{eq:rightembed}, and $\eta$ by adding $AM$ on the left. Their inner product is then
\begin{equation}
 b_{\xi,\eta}=\langle\iota^R_{MB'}(\xi),\id_{AM}\otimes\eta\rangle_{B_k}=(\xi^*\otimes\id_{MB'})(\id_A\otimes\sigma^{-1}_{z,MB'})(\id_{AM}\otimes\eta)\in\End_{\C_k}(AMB')=B_k(K),
 \label{eq:charged-insertion}
\end{equation}
using $\sigma^*=\sigma^{-1}$. Read from bottom to top, $\eta$ creates the charged leg $z$ at the right, the inverse half-braiding transports it to the left across the strands of $M$ and $B'$, and $\xi^*$ absorbs it at the left (Figure~\ref{fig:gap}). For the coefficient \eqref{eq:charged-coeff} the relevant element is $\langle\xi,s_Z^n\eta\rangle_{B_k}=b_{\xi,s_Z^n\eta}$, where $s_Z^n\eta\in\mathsf Y_Z(J+n)$ is the translate of $\eta$ by $n$ sites, for $n$ so large that $J+n$ lies to the right of $I$. Its gap word is $Y^{\otimes g_n}$ with $g_n=g+n$.

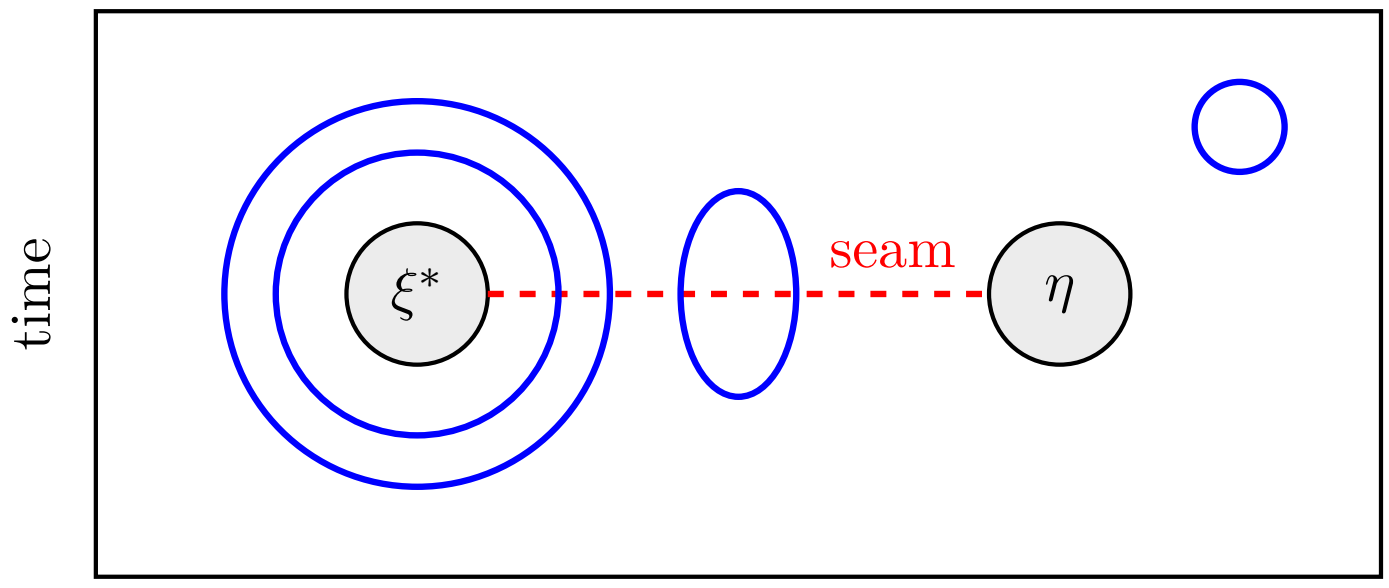
\begin{figure}[ht]
\centering
\begin{tikzpicture}[scale=0.85]
 \draw[thick] (-5,-2.2) rectangle (5,2.2);
 \draw[thick,fill=gray!15] (-2.5,0) circle (0.55);
 \draw[thick,fill=gray!15] (2.5,0) circle (0.55);
 \node at (-2.5,0) {$\xi^*$};
 \node at (2.5,0) {$\eta$};
 \draw[very thick,red,dashed] (-1.95,0) -- (1.95,0);
 \node[red,above] at (1.2,0.05) {seam};
 \draw[very thick,blue] (-2.5,0) circle (1.1);
 \draw[very thick,blue] (-2.5,0) circle (1.5);
 \draw[very thick,blue] (0,0) ellipse (0.45 and 0.8);
 \draw[very thick,blue] (3.9,1.3) circle (0.35);
 \node at (0,-2.6) {space};
 \node[rotate=90] at (-5.5,0) {time};
\end{tikzpicture}
\caption{The sewing of Lemma~\ref{lem:sewing}. Removing small discs around the two charged insertions leaves, on the sphere, an annulus in which the charge line becomes the seam. The two loops surrounding the disc of $\xi^*$ separate the two holes and are essential in the annulus. The loop crossing the seam without surrounding a disc, and the loop on the right, are contractible in the annulus and contribute only through $\delta^{c(D)}$.}
\label{fig:sewing}
\end{figure}

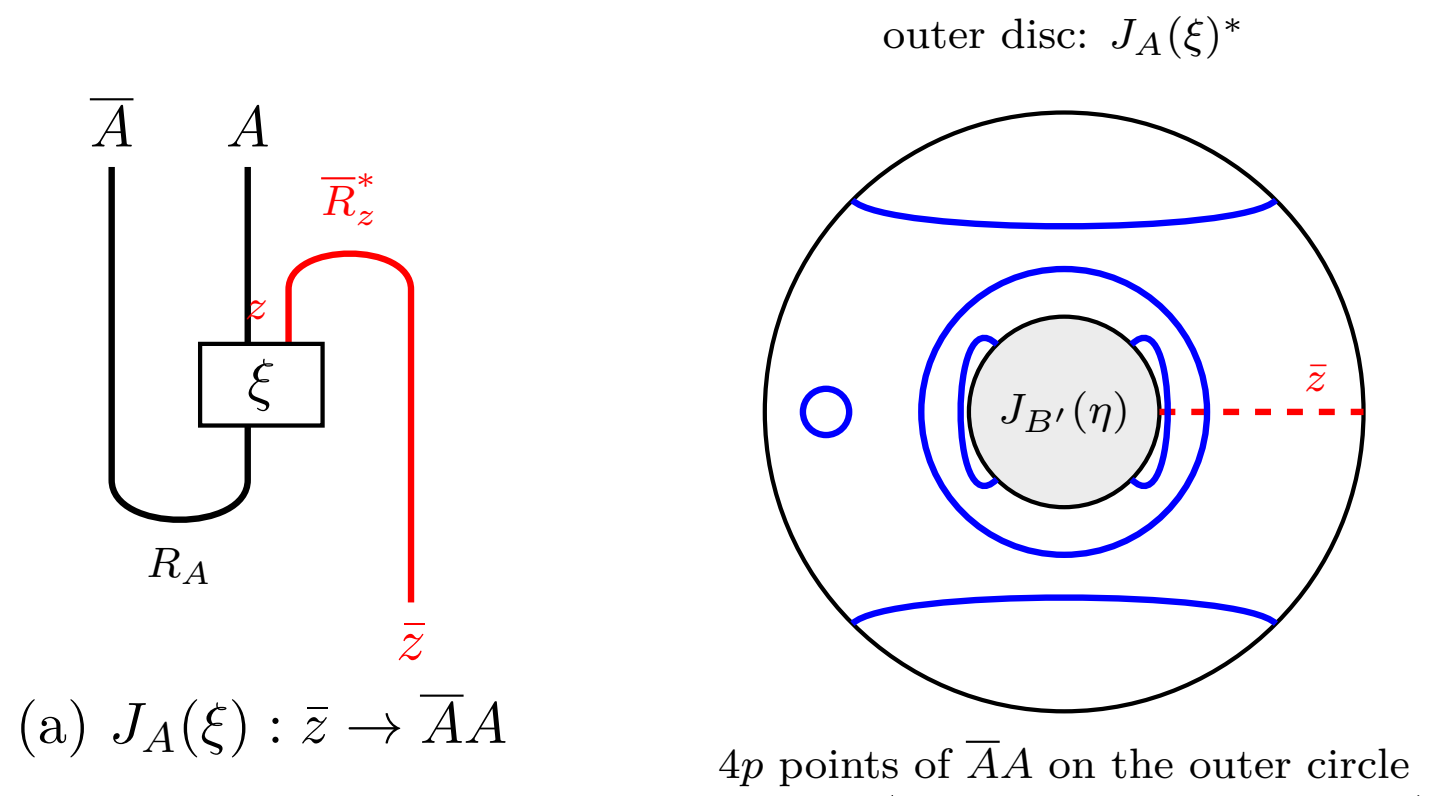
\begin{figure}[ht]
\centering
\begin{tikzpicture}[scale=0.9]
\begin{scope}
 \draw[very thick] (0,3.2) -- (0,0.9) to[out=-90,in=-90] (1,0.9) -- (1,1.3);
 \draw[very thick] (1,1.9) -- (1,3.2);
 \draw[very thick,red] (1.3,1.9) -- (1.3,2.3) to[out=90,in=90] (2.2,2.3) -- (2.2,0);
 \draw[thick,fill=white] (0.65,1.3) rectangle (1.55,1.9);
 \node at (1.1,1.6) {$\xi$};
 \node[above] at (0,3.2) {$\overline A$};
 \node[above] at (1,3.2) {$A$};
 \node[below,red] at (2.2,0) {$\bar z$};
 \node[red,left] at (1.3,2.15) {\scriptsize $z$};
 \node[below] at (0.5,0.55) {\scriptsize $R_A$};
 \node[red,above] at (1.75,2.62) {\scriptsize $\overline R_z^*$};
 \node at (1.1,-0.9) {(a) $J_A(\xi):\bar z\to\overline AA$};
\end{scope}
\begin{scope}[xshift=7cm,yshift=1.4cm]
 \draw[thick] (0,0) circle (2.2);
 \draw[thick,fill=gray!15] (0,0) circle (0.7);
 \node at (0,0) {\scriptsize $J_{B'}(\eta)$};
 \draw[very thick,red,dashed] (0.7,0) -- (2.2,0);
 \node[red,above] at (1.85,0) {\scriptsize $\bar z$};
 \draw[very thick,blue] (0,0) circle (1.05);
 \draw[very thick,blue] (0.495,0.495) .. controls (0.85,0.85) and (0.85,-0.85) .. (0.495,-0.495);
 \draw[very thick,blue] (-0.495,0.495) .. controls (-0.85,0.85) and (-0.85,-0.85) .. (-0.495,-0.495);
 \draw[very thick,blue] (1.556,1.556) .. controls (1.3,1.3) and (-1.3,1.3) .. (-1.556,1.556);
 \draw[very thick,blue] (1.556,-1.556) .. controls (1.3,-1.3) and (-1.3,-1.3) .. (-1.556,-1.556);
 \draw[very thick,blue] (-1.75,0) circle (0.17);
 \node at (0,2.75) {\scriptsize outer disc: $J_A(\xi)^*$};
 \node[below] at (0,-2.25) {\scriptsize $4p$ points of $\overline AA$ on the outer circle};
 \node at (0,-3.05) {(b) $\mathrm{Tr}_{\bar z}\bigl(J_A(\xi)^*\rho_Z(D_\Gamma)J_{B'}(\eta)\bigr)$};
\end{scope}
\end{tikzpicture}
\caption{The endpoint bends and the trace in Lemma~\ref{lem:sewing}. (a) The endpoint map. Here $\xi\in\mathsf Y_Z(I)=\Hom_{\C_k}(A,A\otimes z)$, where $z$ is the underlying object of $Z$ (Section~\ref{sec:affine}), so the box $\xi$ has the input $A$ and the outputs $A$ and $z$. The cup $R_A$ bends the input $A$-string of $\xi$ into $\overline AA$, and the cap $\overline R_z^*$ joins the charge leg $z$ to the input $\bar z$. Likewise $J_{B'}(\eta):\bar z\to\overline{B'}B'$. (b) On the sphere, the disc around $\eta$ carries $J_{B'}(\eta)$, with $4p'$ points of $\overline{B'}B'$ on the inner circle, and the complementary disc carries $J_A(\xi)^*$. The annulus carries $D_\Gamma$, here with arcs on both circles, one essential circle and one contractible loop. The charge line $\bar z$ runs along the seam and crosses the strands by the inverse half-braiding of \eqref{eq:affineaction}. Closing $\bar z$ gives $\mathrm{Tr}_{\bar z}$, which is $\langle J_A(\xi),\rho_Z(D_\Gamma)J_{B'}(\eta)\rangle$ by \eqref{eq:affinespace}.}
\label{fig:endpoint}
\end{figure}

\begin{lemma}[Sewing]
\label{lem:sewing}
Let $\xi,\eta$ and $b_{\xi,\eta}$ be as above, with $K\subseteq I_N$. For each configuration $\Gamma$ on $\Delta_{N,T}$ there is a plain affine Temperley--Lieb diagram $D_\Gamma:(2p',+)\to(2p,+)$ such that the charged evaluation \eqref{eq:evaluation} is (Figure~\ref{fig:endpoint}(b))
\[
 \mathrm{ev}_\Gamma(b_{\xi,\eta})=\langle J_A(\xi),\rho_Z(D_\Gamma)J_{B'}(\eta)\rangle .
\]
Filling the holes of $D_\Gamma$ with the pairings of $\id_A$ and $\id_{B'}$ recovers the diagram of $\Gamma$, so $\ell(D_\Gamma)=\ell(\Gamma)$. The essential circles of $D_\Gamma$ are the loops of $\Gamma$ that separate the two insertions $\xi^*$ and $\eta$.
\end{lemma}
\begin{proof}
Insert $b_{\xi,\eta}$ into the closed loop diagram of $\Gamma$ at time $0$, and regard the resulting closed planar diagram as drawn on the sphere. Remove disjoint small discs around the boxes $\xi^*$ and $\eta$. The complement is an annulus, and the charge line running through the gap becomes a distinguished seam in it (Figure~\ref{fig:sewing}). A loop of $\Gamma$ is essential in this annulus exactly when it surrounds one of the two discs and not the other, that is when it separates the two insertions. Cutting along the seam, which leaves the Temperley--Lieb strands it crosses unchanged, produces a plain affine diagram $D_\Gamma:(2p',+)\to(2p,+)$, with $4p$ and $4p'$ boundary points. Bending the $A$- and $B'$-strings with the standard solutions $R_A,R_{B'}$ and the charge leg with $\overline R_z^*$ produces exactly $J_A(\xi)$ and $J_{B'}(\eta)$ (Figure~\ref{fig:endpoint}(a)). Moving a charge bend through the seam uses
\[
 (\sigma_{z,Q}\otimes\id_{\bar z})(\id_z\otimes\sigma_{\bar z,Q})(\overline R_z\otimes\id_Q)=\id_Q\otimes\overline R_z,\qquad Q\in\C_k,
\]
and with the conjugate equations for $Q$ this gives $(\sigma^{-1}_{z,Q})^\vee=\sigma^{-1}_{\bar z,\overline Q}$, the crossing in \eqref{eq:affineaction}. Gluing the output to the reflected first endpoint gives $\mathrm{Tr}_{\bar z}(J_A(\xi)^*\rho_Z(D_\Gamma)J_{B'}(\eta))$ (Figure~\ref{fig:endpoint}(b)). The argument uses only three facts. The endpoint maps can be undone, $J_A^{-1}J_A=\id$, by the conjugate equations. The displayed identity moves a charge bend through the seam. The half-braiding is natural and compatible with tensor products. None of these depends on how the strands of $\Gamma$ wind around the two insertions, so the identity holds for every configuration $\Gamma$. Erasing the charged insertion restores the identity on $AMB'$, and the holes are filled by the rotated identity pairings.
\end{proof}

The uncharged evaluation is $\mathrm{ev}_\Gamma(1)=\delta^{\ell(\Gamma)}$, and the weight $\delta^{-n(\Gamma)}$ and the boundary cups are common to the charged and uncharged evaluations. So \eqref{eq:estimator} gives
\begin{equation}
 \omega_{N,T}(b_{\xi,\eta})=\mathbb E_{\nu_{N,T}}[F_{\xi,\eta}],\qquad F_{\xi,\eta}(\Gamma)=\frac{\mathrm{ev}_\Gamma(b_{\xi,\eta})}{\mathrm{ev}_\Gamma(1)}=\frac{\langle J_A(\xi),\rho_Z(D_\Gamma)J_{B'}(\eta)\rangle}{\delta^{\ell(\Gamma)}},
 \label{eq:charged-expectation}
\end{equation}
an identity of evaluations in the unitary quotient that does not depend on a Temperley--Lieb expansion of $b_{\xi,\eta}$. Let $r(\Gamma)$ be the number of loops of $\Gamma$ separating the two insertions. By Lemma~\ref{lem:sewing}, $r(\Gamma)=r(D_\Gamma)$ and $\ell(D_\Gamma)=\ell(\Gamma)$, so Lemma~\ref{lem:affine-estimate}, with $m=2p'$, $n=2p$ and the bounds on $\|J_A\xi\|$ and $\|J_{B'}\eta\|$, gives
\begin{equation}
 |F_{\xi,\eta}(\Gamma)|\le\delta^{2(p+p')}\|\xi\|\|\eta\|\,\gamma_Z^{\lfloor r(\Gamma)/2\rfloor}.
 \label{eq:F-bound}
\end{equation}
The constant does not depend on the separation of the insertions, the number of events, or the finite volume cutoffs.

\section{Charged decay and the algebra object}
\label{sec:decay}

\begin{proposition}[Decay of charged coefficients]
\label{prop:decay}
Let $Z\in\mathcal Z(\C_k)$ be simple with $Z\not\cong\one$, and let $U_Z$ be the unitary on $\mathcal K_Z=\mathsf Y_Z\boxtimes_{B_k}\HH$ given by $U_Z(\xi\otimes v)=s_Z(\xi)\otimes Uv$, which implements the translation $\varsigma^2$ in the representation $\pi_Z$ (Section~\ref{sec:affine}). Then $U_Z^n\to0$ in the weak operator topology as $|n|\to\infty$.
\end{proposition}
\begin{proof}
By \eqref{eq:muZstrict}, $\gamma_Z<1$. Let $\xi\in\mathsf Y_Z(I)$ and $\eta\in\mathsf Y_Z(J)$ be local, supported on $p$ and $p'$ sites. By \eqref{eq:charged-coeff} and the discussion after \eqref{eq:charged-insertion},
\[
 \langle\xi\otimes\Omega,U_Z^n(\eta\otimes\Omega)\rangle=\omega_k(b_{\xi,s_Z^n\eta})
\]
for all large $n>0$, and $\|s_Z^n\eta\|=\|\eta\|$. Choose a centre $z$ on $I\times\{0\}$ and $R_0\ge R_*$ such that $Q_{R_0}(z)$ contains the small disc around the insertion $\xi^*$. For $n$ large let $a(n)$ be the largest integer with $L_*^{a(n)}R_0<n/2$, so that $a(n)\ge\log(n/(2R_0))/\log L_*-1$, and $Q_{L_*^{a(n)}R_0}(z)$ does not meet the insertion of $s_Z^n\eta$. Put $R_j=L_*^jR_0$ and let $G_j=G_{R_j}$ be the event of Lemma~\ref{lem:separating} in the annulus $A_{R_j}=Q_{R_{j+1}}(z)\setminus Q_{R_j}(z)$, for $0\le j<a(n)$. On $G_j$ there are two distinct loops in $A_{R_j}$ surrounding the insertion of $\xi^*$ and not that of $s_Z^n\eta$, and loops in different annuli are distinct. So
\[
 r(\Gamma)\ge2\sum_{j<a(n)}\mathbf 1_{G_j}(\Gamma),\qquad\gamma_Z^{\lfloor r(\Gamma)/2\rfloor}\le\prod_{j<a(n)}\gamma_Z^{\mathbf 1_{G_j}(\Gamma)} .
\]
For a single annulus, \eqref{eq:annulus} gives $\mathbb E[\gamma_Z^{\mathbf 1_{G_j}}\mid\mathcal F_{A_{R_j}^c}]=1-(1-\gamma_Z)\,\mathbb P(G_j\mid\mathcal F_{A_{R_j}^c})\le1-p_*(1-\gamma_Z)$. Conditioning successively on the complements of the annuli, from the outermost inwards, as in the proof of Lemma~\ref{lem:looptail} and without any independence of the annuli,
\[
 \mathbb E_{\nu_{N,T}}\bigl[\gamma_Z^{\lfloor r(\Gamma)/2\rfloor}\bigr]\le\bigl(1-p_*(1-\gamma_Z)\bigr)^{a(n)}\le C_Z(1+n)^{-\varepsilon_Z},\qquad\varepsilon_Z=-\frac{\log\bigl(1-p_*(1-\gamma_Z)\bigr)}{\log L_*}>0,
\]
whenever $\Delta_{N,T}$ contains $Q_{R_{a(n)}}(z)$ and the support of $b_{\xi,s_Z^n\eta}$. With \eqref{eq:charged-expectation} and \eqref{eq:F-bound}, $|\omega_{N,T}(b_{\xi,s_Z^n\eta})|\le\delta^{2(p+p')}\|\xi\|\|\eta\|C_Z(1+n)^{-\varepsilon_Z}$ for all such $N,T$. Proposition~\ref{prop:ground} gives $\omega_k(b)=\lim_N\lim_{T}\omega_{N,T}(b)$ for local $b$, so
\[
 |\langle\xi\otimes\Omega,U_Z^n(\eta\otimes\Omega)\rangle|\le C'_{\xi,\eta}(1+n)^{-\varepsilon_Z}
\]
for large $n>0$. Negative $n$ follow from unitarity. Indeed $\langle v,U_Z^{-n}w\rangle=\overline{\langle w,U_Z^nv\rangle}$, so $|\langle v,U_Z^{-n}w\rangle|=|\langle w,U_Z^nv\rangle|$. Density of the vectors $\xi\otimes\Omega$ and unitarity of $U_Z$ give the claim.
\end{proof}

\begin{remark}[The role of the annular structure]
\label{rem:annular}
The annular Temperley--Lieb structure enters the paper only through Proposition~\ref{prop:decay}. The sewing of Lemma~\ref{lem:sewing} turns a charged coefficient, whose strings may wind around the insertions and which involves an unbounded number of events, into a matrix coefficient of the single operator $\rho_Z(D_\Gamma)$ in the Hilbert affine module $\mathscr V_Z$, which depends only on $Z$ by \cite[Theorem~5.13]{DGG}. The loops separating the two insertions become the essential circles of $D_\Gamma$. Each essential circle costs a factor $\mu_Z/\delta$, and $\mu_Z<\delta$ for every simple $Z\not\cong\one$ by the explicit values of Proposition~\ref{prop:muZ}. If $\one\not\subseteq z$, a single separating loop annihilates the coefficient. Lemma~\ref{lem:separating} supplies logarithmically many separating loops, and the contraction $\mu_Z/\delta<1$ per loop converts them into decay. For $Z=\one$ separating loops cost nothing, and there is no decay. At $k=2$ the inequality is checked directly in Remark~\ref{rem:k2affine}, and Part~\ref{part:Ising} uses no annular structure. For $k\ge3$ the inequality is read off from the $S$-matrix of $\mathrm{SU}(2)_k$, through the factorisation of the double for odd $k$ \cite{BEK01}, and for even $k=2n$ through the description of the double with the fixed point $(n,n)$ resolved \cite{EK98b} and the Longo--Rehren half-braiding \cite[Sections~2.2 and~4]{BEK01}. The classification of Hilbert annular Temperley--Lieb modules of Jones and Reznikoff \cite{JR06} gives the same inequality for all $k$ at once, without the values. The rate $\varepsilon_Z$ obtained in this way is not optimal.
\end{remark}

The next step uses only purity, translation covariance of the bimodules and rigidity of $\mathcal Z(\C_k)$.

\begin{lemma}[Rigidity]
\label{lem:rigidity}
Let $\omega$ be a pure $\varsigma^2$-invariant state on $B_k$, with $\HH$, $U$, $\mathcal K_Z$ and $U_Z$ constructed from $\omega$ as above, and suppose $U_Z^n\to0$ weakly for every nontrivial simple $Z$. Then $H_\omega=\one$.
\end{lemma}
\begin{proof}
Let $\mathrm{Int}_Z=\Hom_{B_k}(\HH,\mathcal K_Z)$ be the space of intertwiners from $\pi$ to $\pi_Z$. For $S,S'\in\mathrm{Int}_Z$, $S^*S'\in\pi(B_k)'=\mathbb C\one$ by purity, which defines an inner product by $S^*S'=\langle S,S'\rangle\one$. We show that $\mathrm{Int}_Z$ is finite dimensional. Let $\Phi_1,\dots,\Phi_m$ be orthonormal, so $\Phi_i^*\Phi_j=\delta_{ij}\one$ and $\sum_i\Phi_i\Phi_i^*\le\one$. Let $R:\one\to\overline Z\otimes Z$ and $\overline R:\one\to Z\otimes\overline Z$ be standard solutions of the conjugate equations for $Z$ in $\mathrm{DHR}(B_k)\simeq\mathcal Z(\C_k)$, so that $(\overline R^*\otimes\id_Z)(\id_Z\otimes R)=\id_Z$, $(R^*\otimes\id_{\overline Z})(\id_{\overline Z}\otimes\overline R)=\id_{\overline Z}$ and $R^*R=\overline R^*\overline R=d_Z$ \cite{LR97}. As a morphism of DHR bimodules, $R$ is a bimodule map $B_k\to\mathsf Y_{\overline Z}\boxtimes_{B_k}\mathsf Y_Z$, and
\[
 R_\HH=R\boxtimes\id_\HH:\HH=B_k\boxtimes_{B_k}\HH\longrightarrow\mathsf Y_{\overline Z}\boxtimes_{B_k}\mathsf Y_Z\boxtimes_{B_k}\HH=\mathsf Y_{\overline Z}\boxtimes_{B_k}\mathcal K_Z
\]
is the induced intertwiner, with $R_\HH^*R_\HH=d_Z\one$. Put $S_i=(\id_{\mathsf Y_{\overline Z}}\boxtimes \Phi_i^*)R_\HH:\HH\to\mathcal K_{\overline Z}$, which is well defined because $\Phi_i^*$ is $B_k$-linear. Then $\sum_iS_i^*S_i=R_\HH^*(\id\boxtimes\sum_i\Phi_i\Phi_i^*)R_\HH\le d_Z\one$. The first conjugate equation gives $\Phi_i^*=(\overline R^*\boxtimes\id_\HH)(\id_{\mathsf Y_Z}\boxtimes S_i)$, so $1=\|\Phi_i\|\le\sqrt{d_Z}\|S_i\|$. Each $S_i^*S_i$ is a scalar by purity, at least $1/d_Z$, so $m/d_Z\le d_Z$ and $m\le d_Z^2$.

Translation acts unitarily on $\mathrm{Int}_Z$ by $\Theta_Z(S)=U_ZSU^*$. If $\mathrm{Int}_Z\ne0$, take an eigenvector $S$ with $\Theta_Z(S)=\lambda S$. Since $S^*S=c\one$ with $c>0$, $S\Omega\ne0$, and since $U\Omega=\Omega$, $U_ZS\Omega=\lambda S\Omega$. So $|\langle S\Omega,U_Z^nS\Omega\rangle|=\|S\Omega\|^2$ for all $n$, contradicting weak decay. Hence $\mathrm{Int}_Z=0$ for nontrivial $Z$, and by purity it is $\mathbb C$ for $Z=\one$. These are the multiplicities of the internal end, whose algebra structure is then that of $\one$.
\end{proof}

\begin{theorem}[Critical TLJ states]
\label{thm:TLJ}
For every $k\ge2$, $\omega_k$ is a pure, $\varsigma^2$-invariant ground state for $-\sum_ie_i$, and
\[
 H_{\omega_k}\cong\one\quad\text{in }\mathcal Z(\C_k).
\]
In particular $\omega_k$ is local, and $H_{\omega_k}$ is not Lagrangian.
\end{theorem}
\begin{proof}
Propositions~\ref{prop:ground}, \ref{prop:purity} and~\ref{prop:decay}, and Lemma~\ref{lem:rigidity}. The tensor unit is commutative, and it is not Lagrangian because $\mathcal Z(\C_k)\not\simeq\Vect$.
\end{proof}

At $k=2$ the tensor unit is moreover not Morita equivalent to a Lagrangian algebra (Section~\ref{sec:invariant}), so $\omega_2$ is categorically gapless.

\begin{remark}[Comparison with Part~\ref{part:Ising}]
\label{rem:comparison}
At $k=2$, $\omega_2=\psi_\sharp$ (Section~\ref{sec:TLJchain}), and Theorem~\ref{thm:TLJ} reproves $H_{\phi_\sharp}=\one$. Up to the choice of local vectors, the coefficients \eqref{eq:charged-coeff} are the order, disorder and fermion correlators
\[
 \langle\sigma_x^{(0)}\sigma_x^{(n)}\rangle,\qquad
 \Bigl\langle\prod_{0<j\le n}\sigma_z^{(j)}\Bigr\rangle,\qquad
 \langle i\,b_0a_n\rangle .
\]
At $\phi_\sharp$ the first decays like $n^{-1/4}$ \cite{Pfeuty70}, the second equals the first because $\kappa(\sigma_x^{(0)}\sigma_x^{(n)})=\prod_{j=1}^n\sigma_z^{(j)}$ and $\psi_\sharp\circ\kappa=\psi_\sharp$, and the third is $2/(\pi(2n-1))$ by Proposition~\ref{prop:covariance}. By Remark~\ref{rem:k2affine}, at $k=2$ the coefficient is bounded by a constant times the probability that no loop separates the two insertions.

In Part~\ref{part:Ising} the implementability criterion is necessary and sufficient for quasifree states and also yields $\one\oplus e$ and $\one\oplus m$ off criticality (Corollary~\ref{cor:offcritical-sector-calculation}), whereas decay is only sufficient. Away from the self-dual point, primal connections decay exponentially when $\lambda_{\mathrm b}/\lambda_{\mathrm c}<q$ \cite[Theorem~1.5]{DCLM}, and dual connections when $\lambda_{\mathrm b}/\lambda_{\mathrm c}>q$, since the shift by one strand exchanges primal and dual clusters together with the odd and even rates. A loop surrounding an insertion at radius $R$ borders a primal and a dual cluster of diameter at least of order $R$. Such loops are exponentially rare, and only finitely many loops separate distant insertions. So the argument gives no decay, consistently with the nonzero disorder or order correlator. The exclusion of sectors thus comes from separating loops at every scale.
\end{remark}

\section{Overview}
\label{sec:scope}

Theorem~\ref{thm:TLJ} extends the Hamiltonian half of Theorem~\ref{thm:main} to every unitary Temperley--Lieb--Jones chain below index four. The other results of Part~\ref{part:Ising} use the quasifree structure, and their extension is open. This concerns the transfer matrix state, the noncritical algebra objects, cone Haag duality, factoriality of half-line algebras, standardness, and the conjugacy $\alpha$.

For the transfer matrix state, the critical Ising transfer matrix \eqref{eq:transfer} is, up to a scalar,
\[
 T_K=\exp\Bigl(K\sum_je_{2j}\Bigr)\exp\Bigl(2K\sum_je_{2j-1}\Bigr)\exp\Bigl(K\sum_je_{2j}\Bigr),
\]
which makes sense for every $k$. Since $e^{2Ke}=1+u\,E$ with $u=(e^{2K}-1)/\delta$, the bulk of $T_K^M$ inserts cup--caps with equal weight on odd and even bonds. This is the anisotropic random-cluster model on the square lattice at its self-dual point. The rectangular lattices are isoradial, so \cite{DCLM} is the source of crossing estimates, but we have not carried out the remaining steps. Purity is the main obstacle. $T_K$ selects a different state from $-\sum_ie_i$ ($\phi_c\ne\phi_\sharp$ at $k=2$), so Section~\ref{sec:purity} does not apply. At $k=2$ purity of $\phi_c$ follows from \cite[Theorem~2 and Lemma~5.1]{ArakiEvans83}. For $k\ge3$ no dynamics with this state as ground state is known.

The Hamiltonian $-\sum_ie_i$ is the critical restricted solid-on-solid chain \cite{ABF84}, equivalently the antiferromagnetic anyonic chain of \cite{FTLTKWF07}. Its standard conformal description in the physics literature is the unitary minimal model of central charge $1-6/((k+1)(k+2))$. As in Section~\ref{sec:TL}, this identification is background and is not a consequence of $H_{\omega_k}=\one$, and we do not prove spectral gaplessness of $-\sum_ie_i$ for $k\ge3$.

The restriction to index below four is essential. For $\delta>2$, that is $q>4$, the Temperley--Lieb category is unitary but of infinite depth, so there is no fusion category in place of $\C_k$ and the framework of \cite{EvansJones2025} does not apply as it stands. At the self-dual point the quantum random-cluster model then has coexisting phases. The wired measure has an infinite cluster and free connections decay exponentially \cite[Theorem~1.5]{DCLM}. The crossing estimates used in Section~\ref{sec:loops} fail, and the decay argument gives nothing. Through the same loop representation, Aizenman, Duminil-Copin and Warzel prove breaking of translation symmetry for loop weight greater than two. Examples are the dimerization of the spin-$S$ antiferromagnetic chain with projection interaction, $S>1/2$, and N\'eel order of the XXZ chain at $\Delta>1$ \cite{ADW20}. Correspondingly, the planar Potts and random-cluster transitions are continuous for $1\le q\le4$ \cite{DCST17} and first order for $q>4$ \cite{DCGHMT21}. For $q>4$ the correlation length at the transition point is finite. So there is no diverging correlation length and no scale invariant conformal description. The weakly first order behaviour just above $q=4$ has been interpreted in the physics literature through complex fixed points \cite{GRZ18}.

The chain $-\sum_ie_i$ is no exception, since its ground states are governed by the same quantum random-cluster model. For $\delta=2\cos\gamma\le2$ the Coulomb gas central charge is $c=1-6\gamma^2/(\pi(\pi-\gamma))$ \cite{Nienhuis84}, which equals $1-6/((k+1)(k+2))$ at $\gamma=\pi/(k+2)$ and $1$ at $\delta=2$. For $\delta=2\cosh\lambda>2$ the Bethe ansatz for the XXZ chain at $\Delta=\cosh\lambda$ gives a gap closing only as $\lambda\to0$, like $\exp(-\pi^2/(2\lambda))$ \cite{dCG66}. Since $q-4=4\sinh^2\lambda$, this matches the two-dimensional correlation length. Continuing the formula for $c$ to $\gamma=i\lambda$ gives a complex central charge, which corresponds to the non-Hermitian fixed points of \cite{GRZ18} and not to a ground state of $-\sum_ie_i$. These statements concern the self-dual Hamiltonian $-\sum_ie_i$ only.

\appendix
\section{Standardness of the critical half-line subspaces}
\label{sec:standardness}

This appendix establishes a property additional to the duality theorem. A \emph{ray} $L$ is a half-line $\{k\ge k_0\}$ or $\{k\le k_0\}$ of indices. In \S\ref{sec:sharpstd} the indices are Majorana indices, and in \S\ref{sec:cstd} they are sites, each carrying both Majoranas. Let $Q_L$ be the projection onto vectors supported in $L$, in the complexified real Majorana space, and let $E$ be the basis projection of the state in that presentation. The pair $(E,Q_L)$ is in \emph{generic position} if
\begin{equation}
 \begin{gathered}
 \ran E\cap\ran Q_L=\ker E\cap\ran Q_L=0,\\
 \ran E\cap\ker Q_L=\ker E\cap\ker Q_L=0.
 \end{gathered}
 \label{eq:generic}
\end{equation}

\Needspace{8\baselineskip}
\begin{lemma}
\label{lem:standard}
If \eqref{eq:generic} holds, the closed real half-line field subspace $M_L$ in the Fock one-particle Hilbert space is standard. That is,
\[
 M_L\cap iM_L=0,\qquad\overline{M_L+iM_L}=\HH_1.
\]
Consequently the Fock vacuum is cyclic and separating for the half-line field algebra.
\end{lemma}
\begin{proof}
Let $J$ be the orthogonal complex structure on the real Majorana space. Its two spectral projections on the complexification are $E$ and $1-E$, in an order depending on the Fock convention. If a real vector $v$ and $Jv$ are both supported in $L$, the vectors $(v+iJv)/2$ and $(v-iJv)/2$ belong to the two intersections $\ran E\cap\ran Q_L$ and $\ker E\cap\ran Q_L$. They vanish, proving $M_L\cap iM_L=0$. Applying the same argument to $L^c$ shows that the real orthogonal complement of $M_L+iM_L$ is zero. The vacuum assertion is \cite[Proposition~2.5]{Foit83}.
\end{proof}

\subsection{The quantum chain state}
\label{sec:sharpstd}

By \eqref{eq:Omega}, the convolution operator with matrix entries $\Omega(j-k)$ has Fourier multiplier $i\sgn(p)$ on $(-\pi,\pi)$, using $\widehat f(p)=\sum_k f_ke^{ikp}$. Consequently the two projections $E$ and $1-E$ in Majorana coordinates are the multipliers $\mathbf1_{(0,\pi)}$ and $\mathbf1_{(-\pi,0)}$. Since \eqref{eq:generic} is unchanged on exchanging these projections, either choice gives the same argument. This Fourier transformation sends $Q_L$ for $L=\{k\ge0\}$ to the Hardy projection onto $H^2(\T)$, and for $L=\{k\ge k_0\}$ to the projection onto $e^{ik_0p}H^2(\T)$. Since $E$ commutes with translations, and \eqref{eq:generic} is symmetric under $Q_L\leftrightarrow1-Q_L=Q_{L^c}$, it suffices to treat $L=\{k\ge0\}$.

\begin{proposition}
\label{prop:sharpstandard}
For $\omega_\sharp$, the pair $(E,Q_L)$ is in generic position for every Majorana ray.
\end{proposition}
\begin{proof}
A nonzero Hardy function cannot have boundary values vanishing on a set of positive measure. A function in $\ran E\cap H^2$ vanishes on $(-\pi,0)$, and a function in $\ker E\cap H^2$ vanishes on $(0,\pi)$. Both intersections are zero. The orthogonal complement of $H^2$ is $\overline{H^2_0}$, and complex conjugation reduces the other two intersections to the same boundary uniqueness statement. By the reduction above this covers every ray.
\end{proof}

\subsection{The transfer matrix state}
\label{sec:cstd}

We use the Fourier representation on sites, with two components. For the site ray $L=\{j\ge1\}$, $Q_L$ is the Hardy projection onto $H^2(\T;\mathbb C^2)$. For $\{j\ge j_0\}$ it is its conjugate by multiplication by $e^{\pm i(j_0-1)\theta}$, the sign depending on the Fourier convention. As in \S\ref{sec:sharpstd}, translation invariance and the symmetry of \eqref{eq:generic} under $Q_L\leftrightarrow Q_{L^c}$ reduce every physical ray to $L=\{j\ge1\}$. The $(-1)$-eigenline of $V(\Theta_c)$ is spanned by
\[
 \begin{pmatrix}i\sin(\Theta_c/2)\\\cos(\Theta_c/2)\end{pmatrix}.
\]
We use \emph{meromorphic bounded type} to mean a quotient of bounded holomorphic functions on the disc, with nonzero denominator. Such functions form a field, include rational functions with poles in the disc, and are determined by their almost-everywhere nontangential boundary values. The same boundary class is obtained from quotients of two Hardy functions. For example, an outer multiplier can simultaneously convert two $H^2$ functions to bounded ones. Boundary uniqueness follows by cross-multiplication and Hardy uniqueness.

\begin{lemma}
\label{lem:boundedtype}
There is a nonzero section of $\ran E_c$ in $H^2(\T;\mathbb C^2)$ if and only if $e^{i\Theta_c}$ is the boundary function of a meromorphic bounded-type function.
\end{lemma}
\begin{proof}
For a nonzero Hardy section $(f_1,f_2)$, the eigenline condition is
\[
 \frac{f_1}{f_2}=\tau:=i\tan(\Theta_c/2)
 \quad\text{almost everywhere}.
\]
Its second component cannot vanish identically, since $\cos(\Theta_c/2)\ne0$ almost everywhere. Hence $\tau$ is of meromorphic bounded type. Conversely, write $\tau=p/q$ with bounded holomorphic $p,q$. Then $(p,q)$ is a nonzero Hardy section. Finally,
\[
 \tau=\frac{e^{i\Theta_c}-1}{e^{i\Theta_c}+1},\qquad
 e^{i\Theta_c}=\frac{1+\tau}{1-\tau},
\]
and the bounded-type field is preserved by these transformations.
\end{proof}

\begin{proposition}[The odd-order obstruction]
\label{prop:criticalstandard}
For $0\le z<1$, the pair $(E_c,Q_L)$ is in generic position for every physical ray. Its real half-line field subspace is standard.
\end{proposition}
\begin{proof}
Put $w=e^{i\theta}$. The square of the critical symbol is
\begin{equation}
 e^{2i\Theta_c}=-w\frac{1-z^2w}{1-z^2/w}
 =-\frac{w^2(1-z^2w)}{w-z^2}=:G(w).
 \label{eq:rationalsquare}
\end{equation}
If $e^{i\Theta_c}$ were the boundary function of a meromorphic function $h$ of bounded type, boundary uniqueness would imply $h^2=G$ on the disc. For $0<z<1$, $G$ has a simple pole at $z^2$, whereas every pole of a meromorphic square has even order. For $z=0$, $G(w)=-w$ has a simple zero, giving the same contradiction. Lemma~\ref{lem:boundedtype} therefore proves $\ran E_c\cap\ran Q_L=0$.

The complementary eigenline has component ratio $1/\tau=-i\cot(\Theta_c/2)$. Closure of the field of functions of bounded type under reciprocals shows that this ratio is also not of bounded type, giving $\ker E_c\cap\ran Q_L=0$. For a section in $\ker Q_L=\overline{H^2_0(\T;\mathbb C^2)}$, complex conjugation produces a Hardy section with ratio $-\tau$ or $-1/\tau$, since $\Theta_c$ is real. Both are impossible. This proves all four intersections for $L=\{j\ge1\}$, hence for every physical ray. Apply Lemma~\ref{lem:standard}.
\end{proof}

The odd-order argument proves standardness independently of the conjugacy in Section~\ref{sec:conjugacy}. It is not needed to invoke Foit's duality theorem, and factoriality (Proposition~\ref{prop:factor}) is not inferred from it.

\section*{Use of Artificial Intelligence}

Claude Opus 5.5 was used in assisting with  the proofs, drawing the diagrams and ChatGPT in auditing the manuscript.

\section*{Acknowledgements}

I would like to thank Corey Jones for suggesting the problem of computing the algebra object for the critical states  particularly for the Temperley-Lieb-Jones models and  for numerous discussions.  I am  grateful to Magdalena Musat and Mikael R{\o}rdam for hospitality at the University of Copenhagen during a visit in May 2026 and Denjoe O'Connor at DIAS Dublin during a visit in October-November 2026. This work was supported by the Leverhulme Trust through a Leverhulme Emeritus Fellowship, and in part by a grant of access to OpenAI models through the ChatGPT for Academic Researchers programme.

\bigskip
\noindent{\small David E.~Evans\, evansde@caerdydd.ac.uk\\ School of Computational and Mathematical Sciences, Cardiff University, Cardiff, Wales, UK}

\end{document}